\documentclass[11pt, english, reqno]{amsart}

    \usepackage{amsmath}
    \usepackage{enumerate}
    \usepackage{amsfonts}
    \usepackage{amssymb}
    \usepackage{amsthm}
    \usepackage{bbm}
    \usepackage{xcolor}
    \usepackage{tabularx}
    \usepackage{graphicx}
    \usepackage{xr}
    \usepackage{xspace}
    \usepackage{thmtools}
    \usepackage{thm-restate}
    \usepackage{fourier}
    \usepackage[multiple]{footmisc}
    \usepackage{microtype}
    \usepackage[font=small,format=plain,labelfont=sf,up,textfont=sf,up]{caption}
    \usepackage[tmargin=1in,bmargin=1in,lmargin=1.24in,rmargin=1.24in]{geometry}
    \usepackage{hyperref}

\definecolor{halfgray}{gray}{0.55} 
\definecolor{webgreen}{rgb}{0,.5,0}
\definecolor{webbrown}{rgb}{.6,0,0}
\definecolor{Maroon}{cmyk}{0, 0.87, 0.68, 0.32}
\definecolor{RoyalBlue}{cmyk}{1, 0.50, 0, 0}
\definecolor{Black}{cmyk}{0, 0, 0, 0}

\hypersetup{%
    colorlinks=true, linktocpage=true, pdfstartpage=3, pdfstartview=FitV,%
    breaklinks=true, pdfpagemode=UseNone, pageanchor=true, pdfpagemode=UseOutlines,%
    plainpages=false, bookmarksnumbered, bookmarksopen=true, bookmarksopenlevel=1,%
    hypertexnames=true, pdfhighlight=/O,
    urlcolor=webbrown, linkcolor=RoyalBlue, citecolor=webgreen, 
}

\declaretheoremstyle[ 
    spaceabove=0.5em, spacebelow=0.5em,
    headfont=\bfseries\sffamily,
    notefont=\normalfont, notebraces={[}{]},
    bodyfont=\itshape,
    postheadspace=0.5em,
    numbered=no,
                    ]{myStyle}

\theoremstyle{myStyle}
\declaretheorem[numberwithin=section]{theorem}
\declaretheorem[sibling=theorem, name=Lemma]{lemma}
\declaretheorem[sibling=theorem, name=Remark] {rk}
\declaretheorem[sibling=theorem, name=Definition]{defn}
\declaretheorem[sibling=theorem, name=Proposition]{prop}
\declaretheorem[sibling=theorem, name=Corollary]{corr}

\newtheorem{innercustomthm}{Theorem}
\newenvironment{customthm}[1]
  {\renewcommand\theinnercustomthm{#1}\innercustomthm}
  {\endinnercustomthm}

\newtheorem{innercustomass}{Assumption}
\newenvironment{customass}[1]
  {\renewcommand\theinnercustomass{#1}\innercustomass}
  {\endinnercustomass}

\newcommand{\bt}   {\begin{theorem}}
\newcommand{\et}   {\end  {theorem}}
\newcommand{\bl}   {\begin{lemma}}
\newcommand{\el}   {\end  {lemma}}
\newcommand{\bp}   {\begin{prop}}
\newcommand{\ep}   {\end  {prop}}
\newcommand{\bc}   {\begin{corr}}
\newcommand{\ec}   {\end  {corr}}
\newcommand{\bd}   {\begin{defn}}
\newcommand{\ed}   {\end  {defn}}
\newcommand{\ba}   {\begin{array}}
\newcommand{\ea}   {\end  {array}}
\newcommand{\bi}   {\begin{itemize}}
\newcommand{\ei}   {\end  {itemize}}
\newcommand{\be}   {\begin{enumerate}}
\newcommand{\ee}   {\end  {enumerate}}
\newcommand{\eq}[1]{\begin{equation}#1}
\newcommand{\en}   {\end{equation}}
\newcommand{\eqn}[1]    {\begin{equation} #1 \end{equation}}
\newcommand{\eqan}[1]   {\begin{align} #1 \end{align}}
\newcommand{\eqarray}   {\begin{eqnarray}}
\newcommand{\enarray}   {\end{eqnarray}}

\newcommand{\refeq}[1]{(\ref{e:#1})}
\newcommand{\lbeq}[1] {\label{e:#1}}

\numberwithin{equation}{section}

   \makeatletter 
    \renewcommand\subsection{\@startsection{subsection}{1}{\z@}
                                  {3.25ex \@plus 1ex \@minus .2ex}{.2\baselineskip}
                                  {\normalfont\sffamily}}
    \makeatother

    \newcommand{\vep}           {\varepsilon}

    \newcommand{\eps}           {\varepsilon}

    \renewcommand{\d}           {{\rm d}}

    \newcommand{\jel}           {J_{{\rm E}}(\lambda)}
    
    \newcommand{\cso}[1]        {c_{#1} (\lambda_0)}

    \newcommand{\sss}           { \scriptscriptstyle }

    \newcommand{\nn}            {\nonumber}
    \newcommand{\nnb}           {\nonumber \\}

    \newcommand{\white}         {\color{white}}

\DeclareMathAlphabet{\pazocal}{OMS}{zplm}{m}{n}

\newcommand{\Bcal}  {\pazocal{B}}
\newcommand{\Ccal}  {\pazocal{C}}

\newcommand{\Ecal}  {\pazocal{E}}
\newcommand{\Fcal}  {\pazocal{F}}

\newcommand{\Ncal}  {\pazocal{N}}

\newcommand{\Xcal}  {\pazocal{X}}

\newcommand{\R}     {\mathbb{R}}
\newcommand{\N}     {\mathbb{N}}
\renewcommand{\P}   {\mathbb{P}}

\newcommand{\E}     {\mathbb{E}}
\newcommand{\Z}     {\mathbb{Z}}
\newcommand{\1}     {\mathbbm{1}}

\newcommand{\Q}     {\mathbb{Q}}

\newcommand{\Zd}    {\mathbb{Z}^d}

\newcommand{\bb}    {\underline{b}}
\newcommand{\tb}    {\overline{b}}
\newcommand{\bbe}   {\underline{e}}
\newcommand{\te}    {\overline{e}}
\newcommand{\ole}   {\overline{e}}
\newcommand{\ule}   {\underline{e}}
\newcommand{\oly}   {\overline{y}}
\newcommand{\uly}   {\underline{y}}
\newcommand{\ulb}   {\underline{b}}
\newcommand{\olb}   {\overline{b}}

\newcommand{\Piic}      {\mathbb{P}_{\mathsf{ \sss IIC}}}
\newcommand{\Eiic}      {\mathbb{E}_{\mathsf{ \sss IIC}}}
\newcommand{\Pp}        {\mathbb{P}_p}
\newcommand{\Ppc}       {\mathbb{P}_{p_c}}

\newcommand{\prob}      {{\mathbb P}}

\newcommand{\Eotr}      {E^0_{\omega} \tau_{Q_r}}
\newcommand{\Epc}       {\mathbb{E}_{p_c}}
\newcommand{\expec}     {\mathbb{E}}

\newcommand{\indic}[2][]    {\1_{\{#2\}_{#1}}}

\newcommand{\taupc}         {\tau_{p_c}}

\newcommand{\taup}          {\tau_p}

\newcommand{\indi}          {\mathbbm{1}}

\newcommand{\gqr}           {G_{Q_r}}
\newcommand{\pqr}           {p_{Q_r}^x}

\renewcommand{\rho}         {\varrho}

\newcommand{\nin}       {\notin}
\newcommand{\ua}        {\nearrow}

\newcommand{\conn}      {\leftrightarrow}
\newcommand{\nconn}     {\nleftrightarrow}
\newcommand{\Conn}      {\Longleftrightarrow}

\newcommand{\on}        {{\text{\rm{ on }}}}

\newcommand{\convd}     {\stackrel{d}{\longrightarrow}}

\newcommand{\mconn}     {\stackrel{m}{\longleftrightarrow}}

\newcommand{\limp}          {\lim_{p \ua p_{c}}}

\newcommand{\Reff}      {R_{\rm \sss eff}}

\newcommand{\Br}        {B^{\sss (R)}}
\newcommand{\events}    {\mathfrak{F}}

\newcommand{\BER}       {Q_r}
\newcommand{\Bpiv}[1]   {B_{\rm piv}(#1)}

\newcommand{\VIR}       {|B_r\cap \Ccal(0)|}
\newcommand{\VI}[1]     {|B(#1) \cap \Ccal(0)|}
\newcommand{\VER}       {|Q_r\cap \Ccal(0)|}
\newcommand{\Npiv}[1]   {N_{\rm piv}(#1)}
\newcommand{\Vpiv}[1]   {V_{\rm piv}(#1)}

\newcommand{\Ct}        {{\tilde {\Ccal}^b(0)}}
\newcommand{\Bb}        {\mathsf{Bb}}

\newcommand{\smi}       {\mathsf{S}_m^{\sss \infty}}
\newcommand{\smr}       {\mathsf{S}_m^{\sss (r)}}

\newcommand{\smx}       {\mathsf{S}_m^x}

\newcommand{\sni}       {\mathsf{Z}_m^{\sss \infty}}
\newcommand{\snr}       {\mathsf{Z}_m^{\sss (r)}}
\newcommand{\snR}       {\mathsf{Z}_m^{\sss (R)}}
\newcommand{\snx}       {\mathsf{Z}_m^x}

\newcommand{\iic}       {\mathsf{IIC}}

\newcommand{\liR}       {\Lambda^{\sss \infty}_{\sss (R)}}

\newcommand{\lrR}       {\Lambda^{\sss (r)}_{\sss (R)}}
\newcommand{\lrx}       {\Lambda^{\sss (r)}_{\sss x}}
\newcommand{\bliR}      {(\Lambda^{\sss \infty}_{\sss (R)})^c}
\newcommand{\blir}      {(\Lambda^{\sss \infty}_{\sss (r)})^c}
\newcommand{\blrR}      {(\Lambda^{\sss (r)}_{\sss (R)})^c}

\newcommand{\Be}        {\Ecal_r}

\newcommand{\bbb}       {\mathsf{Bb}}
\newcommand{\bbr}       {U_r^{\mathsf{ \sss Bb}}}

\newcommand{\dqr}       {Q_{r}^{c}}
\newcommand{\tCcal}     {\tilde{\Ccal}}
\newcommand{\dBr}       {\partial B_r}
\newcommand{\Nbb}       {N_{\sss \mathsf{Bb}}}
\newcommand{\Ngood}     {N_{\rm good}}
\newcommand{\Nbad}      {N_{\rm bad}}
\newcommand{\Bear}      {\Ecal_{\vep^a r}}

\newcommand{\zri}       {\mathsf{Z}_r^{\sss \infty}}
\newcommand{\zrx}       {\mathsf{Z}_r^x}
\newcommand{\zni}       {\mathsf{Z}_n^{\sss \infty}}
\newcommand{\znx}       {\mathsf{Z}_n^x}

\newcommand{\tqr}       {\tau_{Q_{r}}}
\newcommand{\bqr}       {U_r}
\newcommand{\bqhr}      {U_{\vep r/2}}
\newcommand{\bqi}[1]    {U_{#1}}

\newcommand{\Fe}		{\Fcal_r}
\newcommand{\Beh}	{{\mathsf{\Ecal}_{\eta r}}}
\newcommand{\Zer}	       {{Z_{\sss \Bb} (\eta r, r)}}

\newcommand{\twa}       {{(2 \wedge \alpha)}}

\newcommand{\ReffB}     {\Reff (0, Q_{r}^c)}
\newcommand{\Reffbb}    {\Reff^{\mathsf{ \sss Bb}}}
\newcommand{\Reffmod}   {\Reff^{\mathsf{ \sss mod}}}
\newcommand{\Reffpn}    {R_{\rm \sss eff}(0,\ole_n)}

\newcommand{\EVBr}      {V(B_r)}
\newcommand{\Vr}        {V(\bqr)}

\newcommand{\Sbbr}      {\lvert \bbr \rvert}

\begin{document}

\title{Random walk on the high-dimensional IIC}
\author{Markus Heydenreich}
\address{Mathematisch Instituut, Universiteit Leiden, P.O.~Box 9512, 2300 RA Leiden, The~Netherlands;
 Centrum voor Wiskunde \& Informatica, P.O.~Box 94079, 1090~GB Amsterdam, The~Netherlands.}
\email{markus@math.leidenuniv.nl}
\author{Remco van der Hofstad}
\address{Department of Mathematics and
	Computer Science, Eindhoven University of Technology, P.O.\ Box 513,
	5600 MB Eindhoven, The Netherlands.}
\email{rhofstad@win.tue.nl}
\author{Tim Hulshof}
\address{Department of Mathematics, the University of British Columbia, Vancouver, Canada.}
\address{Pacific Institute for Mathematical Sciences, the University of British Columbia, Vancouver, Canada.}
\email{thulshof@math.ubc.ca}

\begin{abstract}We study the asymptotic behavior the exit times of random walk from Euclidean balls around the origin of the incipient infinite cluster in a manner inspired by \cite{KumMis08}. We do this by getting bounds on the effective resistance between the origin and the boundary of these Euclidean balls. We show that the geometric properties of long-range percolation clusters are significantly different from those of finite-range clusters. We also study the behavior of random walk on the backbone of the IIC and we prove that the Alexander-Orbach conjecture holds for the incipient infinite cluster in high dimensions, both for long-range percolation and for finite-range percolation.
\end{abstract}
\date{\today}
\maketitle
\vspace{1cm}
{\small

}

\vspace{1cm}
\hrule
\vspace{1em}

\section{Introduction and main results}
We study properties of random walk on the \emph{incipient infinite cluster} (IIC) of $\Zd$. The IIC is an infinite random subgraph of $\Zd$ that was proposed by physicists in the 1980's as an infinite analogue of a critical percolation cluster (see e.g. \cite{AleOrb82}, \cite{LeySta83}).

Percolation is a model that, given a graph $G$, generates random subgraphs of $G$ by independently retaining edges according to a Bernoulli process with parameter $p$, and removing them otherwise. When $G = \Zd$ (i.e., the $d$-dimensional integer lattice with $d \ge 2$), there is a non-trivial value $p_c$, the \emph{critical threshold}, so that when $p < p_c$ the model almost surely does not generate an infinite connected subgraph, whereas when $p>p_c$, the model does generate a unique infinite connected subgraph almost surely. We are interested in the case where $p = p_c$, i.e., \emph{critical} percolation.
It is widely believed that critical percolation clusters of $\Zd$ are almost surely \emph{finite} whenever $d \ge 2$. This has been proved when $d=2$ \cite{Kest80}, and when $d$ is `high enough' \cite{BarAiz91}, \cite{HarSla90a}. This paper focusses on the high-dimensional setting.

The asymptotic behavior of random walk reveals a lot about the structure of the graph it walks on. Of particular interest to us are the \emph{exit time} $\tau_A$ of the walk from a set $A$ (for particular choices of $A$), and the \emph{return probability} $p_n (x,x)$, i.e., the probability that a walk started at $x$ returns to $x$ after $n$ steps. The main focus of this paper is the scaling behavior of the exit time from balls as the radius of the balls increases.

Our motivation for studying random walk on the IIC rather than on critical clusters is that the IIC is an infinite graph that is constructed to locally `look' like a (very large) critical cluster. This way we can use random walk asymptotics to study geometric properties of critical percolation clusters without having to deal with finite-size effects. To generate the IIC we rely on the construction of the IIC-measure from $\Piic$ \cite{HofJar04}, \cite{HeyHofHul12a}. We will elaborate on this construction below (cf.\ \refeq{IICdef}).

Our main results consist of: (1) asymptotic bounds on the random walk return probability, and (2) bounds on the random walk exit time from the intersection between the IIC and Euclidean and intrinsic balls. Of course, the bounds we get on the exit time from a ball depend crucially on the metric that we use to define that ball. In this paper we consider two metrics for this purpose: the \emph{extrinsic distance metric} $\lvert \;\cdot\; \rvert$ (or Euclidean metric, or $L^2$-norm) and the \emph{intrinsic distance metric} $d_G$ (or graph metric). The intrinsic distance metric between two vertices of a graph is the shortest distance between them \emph{through} the graph. Hence, it is sensitive to the topology of the graph, but not to the topology of the space the graph may be embedded in. The extrinsic distance metric, on the other hand, is a metric of $\Zd$, so it is sensitive to this topology, but it ignores the topology of embedded graphs.
Using both metrics will bring to light some fundamental similarities and differences between the various percolation models that we study. In particular, we observe that properties of the random walk that have to do with the \emph{graph structure} (e.g.\ return probabilities, exit times from intrinsic metric balls) are universal for a broad class of models, whereas those properties that have to do with the \emph{spatial structure} (e.g.\ exit time from extrinsic metric balls) are shared only among models that have similar local properties (that is, similar edge probability distributions, see below).

\subsubsection*{Our contribution.} This paper establishes novel bounds on random walk exit times for the IIC in high dimensions. In more detail, our contributions are as follows:

\begin{enumerate}
\item[(1)] We identify the asymptotics of the exit times from Euclidean or extrinsic balls for random walks on the IIC for finite-range percolation.
\item[(2)] We generalize the results by Nachmias and Kozma \cite{KozNac09} on random walk exit times in the intrinsic distance, and show that these hold under the strong triangle condition rather than on upper bounds in $x$-space on the two-point function, which, in applications is a weaker condition.
\item[(3)] We identify the asymptotics of the exit times of Euclidean and intrinsic balls for random walks on the IIC backbone for finite-range percolation, and show that these obey similar scaling as the exit times of random walk on a random walk trace.
\item[(4)] We extend all the above results to the IIC for long-range percolation, and show that while the results in the intrinsic distance are unchanged, the results in the Euclidian or extrinsic distance depend sensitively on the long-range nature of the percolation model. 
\end{enumerate}

In Theorem \ref{th:loglimits} below we summarize some of our main results as they apply to three important percolation models. A precise definition of these models is given further along in this section. A few brief definitions and remarks are needed before we state this theorem: Let $\omega$ be a subgraph of the complete graph on $\Zd$,  $(\Zd, \Zd \times \Zd)$ with induced graph metric $d_\omega$. Given $r \in \R$ and a configuration $\omega$, define the sets
\begin{equation}
    Q_r (x) = \{ y \in \Zd \colon \lvert x -y\rvert \le r\} \qquad \text{ and } \qquad B_r (x ; \omega) = \{y \in \Zd \colon d_{\omega} (x, y) \le r\}.
\end{equation}
We call $Q_r (x)$ the \emph{extrinsic ball of radius $r$ around $x$} and $B_r (x;\omega)$ the \emph{intrinsic ball of radius $r$ around $x$}. Typically, we write $B_r(x)$ instead of $B_r(x;\omega)$.
Note that $Q_r(x)$ is a deterministic set while $B_r(x)$ is a random set (if $\omega$ is random).

Given a random walk on $\omega$ started at $0$, we write $\tau_{B_r}$ and $\tau_{Q_r}$ for the exit times of that random walk from $B_r$ and $Q_r$. The probability measure $P_\omega^0$ and the expectation $E_{\omega}^0$ below are for a random walk started at $0$ on a \emph{fixed} $\omega$, that is, they only consider the randomness of the walk.
\begin{customthm}{0}\label{th:loglimits}
    Let $\Piic$ be the IIC-measure for critical percolation on $\Zd$ with $d \ge 19$ for nearest-neighbor percolation {\sc (nnp)}, $d>6$ for finite-range spread-out percolation {\sc (frp)} and $d > 3 \twa$ for long-range spread-out percolation {\sc (lrp)} with decay exponent $\alpha \in (0,\infty)$. Consider a simple random walk on the IIC. Then, $\Piic$-a.s.,
   \begin{equation}\label{e:AlexanderOrbach}
      \lim_{n\to\infty}\frac{\log p_{2n}(0,0)}{\log n}=-\frac23\qquad\text{and}\qquad
      \lim_{n\to\infty}\frac{\log|W_n|}{\log n}=\frac23
    \end{equation}
    where $W_n$ is the range of random walk after $n$ steps. Furthermore, for $\Piic$-almost all $\omega$,
    \begin{equation}\label{e:Exittimethm}
        \lim_{r\to\infty}\frac{\log E_{\omega}^0 \tau_{B_r}}{\log r}=3 \quad\text{ for {\sc nnp}, {\sc frp} \& {\sc lrp}}\qquad \text{and}\qquad
        \lim_{r\to\infty}\frac{\log E_{\omega}^0 \tau_{Q_r}}{\log r}  =6  \quad\text{ for {\sc nnp} \& {\sc frp},}
    \end{equation}
    whereas, for {\sc lrp}, uniformly in $r$,
    \begin{equation}
        \int P_{\omega}^0 \left(\tau_{Q_r} > \theta r^{3(4 \wedge \alpha)/2}\right) \Piic(\d \omega) \to 0 \quad \text{ as } \theta \to \infty.
    \end{equation}
\end{customthm}
(Kozma and Nachmias proved \refeq{AlexanderOrbach} and the left-hand limit of \refeq{Exittimethm} for {\sc nnp} and {\sc frp} in \cite{KozNac09} -- we mention these here for completeness.)

Let us summarize the contributions that we make to the literature in this paper: we estimate exit times from the extrinsic ball $Q_r$ for all three percolation models, we bound exit times from the intrinsic ball $Q_r$ for long-range percolation, and we prove exit time estimates for random walk on the IIC backbone.

A paper by Kumagai and Misumi \cite{KumMis08} provides most of the tools that we need to prove the above theorem (as well as stronger results, given in Theorems  \ref{th:Eucexittimes} and \ref{th:Eucwalkdist} -- \ref{th:IIC-bb2}
 below). In this paper, they prove that bounds on the volume and effective resistance of a graph intersected with a ball of a given metric imply bounds on the random walk return probability and exit time from that ball. (Their results generalize results of Barlow, J\'arai, Kumagai and Slade \cite{BarJarKumSla08}, where such bounds are obtained specifically for balls in the intrinsic metric.) Most of the work in this paper goes into proving the required bounds on the volume and effective resistance (see Definition \ref{def:radiussets} and Theorem \ref{th:Goodballs} for a precise statement of what they are). We prove these bounds for both the extrinsic and the intrinsic distance metric, and we prove them for a broad class of percolation models.

\subsection{High-dimensional percolation}
\subsubsection*{The triangle condition.}
Our results apply to percolation models that satisfy the so-called `strong triangle condition'. Define the \emph{triangle diagram} $\triangle_p (0)$ by
\begin{equation}\label{e:trianglecond}
    \triangle_p (0) \equiv \sum_{x,y \in \Zd} \Pp (0 \conn x) \Pp (x \conn y) \Pp (y \conn 0),
\end{equation}
where $\Pp(a \conn b)$ is the probability that there is a path between $a, b \in \Zd$ in the subgraph of $\Zd$ generated by the percolation measure with parameter $p$. The \emph{triangle condition} is that $\triangle_p (0)$ is finite when $p \le p_c$. The \emph{strong triangle condition} is that $\triangle_{p} (0) = 1 + O(\beta)$ when $p \le p_c$, for a sufficiently small parameter $\beta$,
where $\beta=1/(2d)$ for nearest-neighbor percolation and $\beta=L^{-d}$ for spread out models. 
The strong triangle condition is known to hold for nearest neighbor-percolation when $d$ is sufficiently large, and in spread-out models when $d$ is larger than some \emph {upper critical dimension}  $d_c$ and $L$ is sufficiently large. 
From here on, when we say that a result holds `in high dimension', we mean that the dimension is high enough for the strong triangle condition to be satisfied with sufficiently small $\beta$.

\subsubsection*{Percolation models.}
In independent percolation models the probability that an edge $\{x,y\}$ is retained can be described by a transition kernel $D(x,y)$, and a parameter $p \in [0, \|D\|_{\infty}^{-1}]$. The results in this paper hold for models that satisfy the strong triangle condition with a sufficiently small $\beta$, and whose associated edge weight distribution $D(\;\cdot\;,\;\cdot\;)$ is invariant under the symmetries of $\Zd$.  Moreover, the results in our paper hold in the general regime of percolation models considered in \cite{HeyHofHul12a}. 
For the purpose of presentation, we further restrict our attention to the following three `standard' models.

The first model to consider is \emph{nearest-neighbor percolation}. The underlying graph of this model is $\Zd$ with edge set $E = \{\{x,y\}\colon x,y \in \Zd, |x-y| =1\}$. All edges are retained independently with probability $p/2d$ where $p \in [0,2d]$ and removed with probability $1-p/2d$ (note that we follow the convention of the high-dimensional percolation literature and choose $p$ in such a way that $D(x,y) = (1/2d)\indi_{\{|x-y|=1\}}$ is a normalized transition kernel). Hara and Slade proved that nearest-neighbor percolation on $\Zd$ satisfies the strong triangle condition when $d \ge 19$ \cite{HarSla94}, and Fitzner and the second author announced a proof that shows that this holds for $d \ge 15$ (see \cite{Fitz13}), although it is generally believed that $d>6$ is enough.

The second model is \emph{finite-range spread-out percolation}. Now the underlying graph is the complete graph with vertex set $\Zd$. The probability that an edge is retained is positive and the same for all edges up to length $L$, and $0$ for longer edges, i.e.,
\begin{equation}
    p D(x,y) = \frac{p}{(2 L +1)^d -1} \indi_{\{0< \| x-y\|_{\infty}  \le L\}}.
\end{equation}
The parameter $L$ is known as the \emph{spread-out parameter}, and it is typically chosen to be large for technical reasons. For this model it has been proved that the strong triangle condition is satisfied when $d > 6$ and $L$ is sufficiently large \cite{HarSla90a}.

Finally, we consider \emph{long-range spread-out percolation}. Again, we use the complete graph with vertex set $\Zd$. The probability that an edge is retained decays as a power-law with the (extrinsic) distance between its ends, that is, for an edge $\{x,y\}$,
\begin{equation}\label{e:DLRP}
    p D(x,y) = p \frac{\Ncal_L}{\max \{|x-y|/L , 1\}^{d+\alpha}},
\end{equation}
for $\alpha \in (0,\infty)$ and where $\Ncal_L$ is a normalizing constant. For long-range spread-out percolation $d$ is high enough when $d > 3 \twa$ \cite{HeyHofSak08}. The decay exponent $\alpha$ determines the decay of the edge retention probability as a function of the length of the edge. As can be seen from the definition, when $\alpha \le2$ the spatial variance $\sum_x \lvert x \rvert^2 D(0,x)$ becomes infinite, whereas the spatial variance is finite when $\alpha > 2$. As a result the long-range model behaves different for $\alpha < 2$ and $\alpha > 2$. From here on we take $\alpha \in (0,2)\cup(2,\infty)$, that is, we do not consider the case where $\alpha=2$. When $\alpha=2$ we get logarithmic corrections on many of the results that follow, and these make it cumbersome to read.

Unless we say otherwise, the results below hold for finite-range models. But we will state results in terms of the parameter $\alpha$ whenever the result also applies to long-range percolation. To make sense of these results for models that do not depend on the parameter $\alpha$ one should think of $\alpha$ as a redundant parameter that is always set to $\infty$.

\subsubsection*{The IIC-measure.}
We cannot construct an IIC-measure by simply conditioning the critical percolation measure on the event that the cluster of the origin is infinite because this is an event of measure $0$. But it is a well-known property of high-dimensional critical percolation that in a box of linear size $n$ there is a cluster whose size is of order $n$ with high probability \cite{Aize97}. In other words, large critical clusters are common. We can use this fact to condition the critical percolation measure on an event that implies that the origin is part of a cluster whose size is proportional to $n$ (e.g.\ the event that $0$ is connected to a point at distance at least $n$). Taking the limit $n \to \infty$ yields an IIC-measure. This needs to be proved, and these proofs are typically quite involved (for high-dimensional models one needs to use lace-expansion techniques). It turns out that several different limiting schemes can give the same IIC-measure. This has been proved for both two- and high-dimensional percolation and oriented percolation models (although most schemes have not been shown in all three settings) cf.\ \cite{HeyHofHul12a, HofHolSla02, HofJar04, Jara03b, Kest86a}.

The particular scheme that we use in the proofs of this paper relies on the expected cluster size, or \emph{susceptibility} of a percolation model, which is defined by
\begin{equation}
    \chi(p) \equiv \E_p[\lvert \Ccal(0) \rvert] = \sum_{x \in \Zd} \Pp(0 \conn x),
\end{equation}
where $\Ccal(0)$ is the connected component containing the origin and $\lvert \Ccal(0) \rvert$ denotes the number of vertices in $\Ccal(0)$.
The susceptibility is finite when $p<p_c$, but it diverges when $p$ approaches $p_c$ from below. With this in mind, \cite{HofJar04} proposes the limiting scheme
\begin{equation}\label{e:IICdef}
    \Piic(F)  \equiv \limp \frac{1}{\chi(p)} \sum_{x \in \Zd} \Pp(F \cap \{0 \conn x\}), 
    \qquad F \in \events_0,
\end{equation}
where $\events_0$ is the algebra of cylinder events. 
This extends to an IIC measure on the $\sigma$-algebra generated by $\events_0$.  
The second author and J\'arai \cite{HofJar04} proved that the limit exists for models where the \emph{two-point function}
\begin{equation}
    \tau_p (x) \equiv \Pp(0 \conn x)
\end{equation}
obeys the asymptotic relation $\taupc(x) = c|x|^{2-d}(1+o(1))$ as $x \to \infty$. This asymptotic relation holds for finite-range models (cf.\ \cite{HarHofSla03} and \cite{Hara08}), but does not hold for long-range models when $\alpha <2$. In a previous paper \cite{HeyHofHul12a} we prove that the same limit also exists under the weaker condition that the strong triangle condition holds, so the limit also holds for the long-range models we discuss in this paper.

An IIC configuration contains a special subgraph, the \emph{backbone}, that consists of all  vertices $x\in\Zd$ (and the edges between them) in the IIC with the property that there is a path of open edges from $0$ to $x$ and disjoint from this path there is another path from $x$ to $\infty$. Given a configuration $\omega$ we write $\bbb(\omega)$ for the backbone of $\omega$. We say an edge $e$ in $\bbb(\omega)$ is \emph{backbone-pivotal} if $e$ is open, and if closing $e$ would disconnect $0$ and $\infty$. It has been proved that the backbone is \emph{essentially unique} \cite{HofJar04}. This means that any two infinite self-avoiding paths started at $0$ share an infinite number of edges. The set of edges shared by all infinite self-avoiding paths is exactly the set of backbone-pivotal edges.

For a more in-depth discussion of the construction of IIC-measures we refer the reader to \cite{HeyHofHul12a}.

\subsubsection*{Random walk.}
In this paper random walks and the associated spaces are defined as follows:
Let $\Omega = \{0,1\}^{E(\Zd)}$, where $E(\Zd)$ is the set of edges (typically, our percolation models require $E(\Zd) = \Zd \times \Zd$), so that $\Omega$ is the state space of percolation configurations. Consider the probability space $(\Omega, \Fcal, \Piic)$ that describes the family of random graphs  $\Gamma_{\mathsf{ \sss IIC}} (\omega) = \left(\left(\iic(\omega), E(\omega)\right): \omega \in \Omega \right)$, where $\iic(\omega)$ is the set of vertices of the (unique, infinite) connected component of $0$ in $\omega$, and $E(\omega)$ is the associated edge set. 
Let  $\Xcal=\big((X_n)_{n\ge0} , P_{\omega}^{x}, x \in \iic (\omega)\big)$ denote simple random walk on $\Gamma_{\mathsf{ \sss IIC}}$ started at $x$. While $(\Omega, \Fcal, \Piic)$ denotes the probability space of the random environment $\Gamma_{\mathsf{ \sss IIC}}$, we denote by $(\overline{\Omega}, \overline{\Fcal})$ a second space for the law of the random walk $\Xcal$ on $\Gamma_{\mathsf{ \sss IIC}}(\omega)$, so that the random walk on a random environment $\Xcal$ is defined on the product $\Omega \times \overline{\Omega}$.

\subsubsection*{Two important assumptions}
The proof of Theorem \ref{th:Eucexittimes} below uses the asymptotics of the \emph{extrinsic one-arm probability} of critical percolation, i.e.\ the probability that the origin of $\Zd$ is connected to a point at (at least) distance $r$:
\begin{customass}{O}\label{ass:onearm}
    The extrinsic one-arm probability satisfies
        \begin{equation}
            \Ppc(0 \conn \dqr) \le C r^{-2}
        \end{equation}
        for some constant $C>0$.
\end{customass}
Kozma and Nachmias proved this assumption for finite-range percolation models in high dimensions \cite{KozNac11}. In \cite{HeyHofHul12a} we prove that the one-arm probability of long-range percolation is bounded from below by $c r^{-(4 \wedge \alpha)/2}$. Hence, if this bound is sharp, Assumption \ref{ass:onearm} also holds for long-range percolation when $\alpha \ge 4$. But such an upper bound cannot hold when $\alpha<4$ (even if it turns out that the lower bound is not sharp). Theorem \ref{th:LRP-IIC} below illustrates how this phase transition at $\alpha=4$ affects the exit time of random walk.

Several proofs in this paper use the assumption that the backbone gives rise to a process on $\Zd$ whose scaling limit is either Brownian motion (for finite-range models and long-range models with $\alpha > 2$) or a symmetric $\alpha$-stable motion (for long-range models with $\alpha < 2$). To make this assumption precise we first have to define such a process on $\bbb(\omega)$.

We can give a unique ordering, say $\{e_i\}_{i=1}^{\infty}$ to the backbone-pivotal edges of $\bbb(\omega)$ by considering the order in which they are crossed by any and every infinite self-avoiding path on $\bbb(\omega)$ started at $0$. Since all $e_i$ are crossed exactly once by any infinite self-avoiding path, we can assign a \emph{top} and \emph{bottom} to these edges, e.g.\ $e_i = (\ole_i, \ule_i)$, such that any self-avoiding path started at $0$ crosses $\ule_i$ before it crosses $\ole_i$. Let $S_n = \ole_n$, then $(S_n)_{n=0}^{\infty}$ is the stochastic process of the position of the top of backbone pivotal edges (where we set $\ole_0 = 0$). Consider the rescaled process
\begin{equation}\label{e:Xn-def}
    X_n (t) \equiv (v_{\alpha} n)^{-1/\twa} S_{\lceil n t \rceil}, \qquad t \in [0,1].
\end{equation}
We assume the following behavior:
\begin{customass}{S}\label{ass:scaling}
    As $n \to \infty$, the process $X_n (t)$ converges in distribution to an $\alpha$-stable L\'evy motion when $\alpha < 2$, and to a Brownian motion when $\alpha > 2$.
\end{customass}
A proof of Assumption \ref{ass:scaling} is in preparation by the authors and Miermont \cite{HeyHofHulMie12a}. It should be remarked that this assumption is stronger than what we actually use, see Proposition \ref{prop-piv-bd} below.

\subsection{Main results.}
We now state our main results: first we present results about the Euclidean (extrinsic) metric, then we present results about the graph (intrinsic) distance metric. 
\subsubsection*{Results for extrinsic distances}
The first theorem gives upper and lower bounds for various quantities related to the exit time of random walk from extrinsic balls:
\bt[Extrinsic random walk geometry of the IIC]\label{th:Eucexittimes} Let $r \ge 1$. If the strong triangle condition is satisfied for some sufficiently small $\beta$, and if either the model is finite-range, or if Assumptions \ref{ass:onearm} and \ref{ass:scaling} hold, then the following holds:
    \begin{itemize}
        \item[(a)]  Uniformly in $r$,
                    \begin{equation}\label{e:RXI}
                        \Piic \left(\theta^{-1} r^{6} \le E_{\omega}^{0} \tau_{Q_r} \le \theta r^{6} \right) \to 1\qquad \text{ as }\quad\theta \to \infty.
                    \end{equation}
        \item[(b)]  There exists $r^\star \ge 1$ such that, for all $r \ge r^{\star}$,
                    \begin{equation}\label{e:RXVI}
                        c_1 r^{6} \le \Eiic[E_{\omega}^{0} \tau_{Q_r} ] \le c_2 r^{6}.
                    \end{equation}
        \item[(c)]  There exists $\gamma_1 < \infty$ and a subset $\Omega_0 \subset \Omega$ with $\Piic(\Omega_0)=1$ such that for all $\omega \in \Omega_0$ and $x \in \iic(\omega)$, there exists $R_x (\omega)< \infty$ such that
                    \begin{equation}\label{e:RXX}
                        (\log r)^{-\gamma_1} r^{6} \le E_{\omega}^{x} \tau_{Q_r} \le (\log r)^{\gamma_1} r^{6}, \qquad \forall r \ge R_x (\omega).
                    \end{equation}
        \item[(d)]  For all $(\omega,\bar\omega) \in \Omega_0\times\overline\Omega$ and $x \in \iic(\omega)$, there exists $\gamma_2 < \infty$ and $R_x (\omega, \bar \omega)$ such that $P_{\omega}^{x} (R_x < \infty) = 1$ and such that
                    \begin{equation}
                        (\log r)^{-\gamma_2} r^{6} \le \tau_{Q_r}(\omega, \bar \omega) \le (\log r)^{\gamma_2} r^{6}, \qquad \forall r \ge R_x (\omega, \bar \omega).
                    \end{equation}
    \end{itemize}
\et
As is explained in \cite{BarJarKumSla08}, it is unlikely that it will turn out that $\gamma_1, \gamma_2=0$, since Barlow and Kumagai in \cite{BarKum06} have shown that the tree analogue of the IIC exhibits $(\log \log r)^c$ fluctuations, for some $c>0$. It would be interesting to see whether random walk on the high-dimensional IIC also exhibits these $\log \log$ fluctuations.

We were not able to establish similar estimates for long-range percolation, but we found that  the behavior $\tau_{Q_r}$ changes radically for long-range percolation when $\alpha < 4$, as the following theorem shows:
\bt[Extrinsic random walk geometry of the LRP-IIC]\label{th:LRP-IIC}
    Consider long-range spread-out percolation with parameter $\alpha$ that satisfies the strong triangle condition for some sufficiently small $\beta$. There exist $c, \vep > 0$ and $r^\star (\lambda)$ such that, for any $r \ge r^\star (\lambda)$,
   \begin{equation}
        \mathsf{P}^\star \left(\tau_{Q_r} \le \lambda r^{3 (4 \wedge \alpha)/2}\right) \ge 1 -c/\lambda^\vep.
    \end{equation}
\et
Theorem \ref{th:LRP-IIC} is much weaker than Theorem \ref{th:Eucexittimes} but it does show that the exit time for long-range spread-out percolation when $\alpha<4$ typically comes \emph{much} sooner than it does for finite-range models. From Theorem \ref{ExpectationBounds} below, it can be seen that the effect of the presence of long-range edges is that clusters become more `smeared out' in space. If this was the only effect that the presence of long-range edges has, then we would expect the exit time to be of the order of $r^{3 \twa}$.
Heuristically, the reason that $r^{3\twa}$ is not the correct order for long-range percolation when $\alpha<4$ is that (1) typically, there are relatively many open edges with length at least $2r$ and with one end in $Q_r$ and (2) once the random walker in $\iic \cap Q_r$ crosses such an edge, it immediately enters $Q_r^c$ and the exit time is reached.

The following theorem involves the \emph{annealed law}
\begin{equation}\label{e:Pstar}
    \mathsf{P}^{\star}(\;\cdot\;) := \int P_\omega^0 (\;\cdot\;) \Piic(\d \omega).
\end{equation}
\bt[Extrinsic distance of random walk from the origin]\label{th:Eucwalkdist}
Let $n \ge 1$. If the strong triangle condition is satisfied for some sufficiently small $\beta$, and if either the model is finite-range, or if Assumption \ref{ass:scaling} holds, then the following holds:
    \begin{itemize}
        \item[(a)]  uniformly in $n$,
                    \begin{equation}
                        \mathsf{P}^{\star} \left(\lvert X_n\rvert < \theta n^{1/(3\twa)}\right) \to 1 \qquad \text{ as }\quad\theta \to \infty
                    \end{equation}
                    and
                    \begin{equation}
                        \mathsf{P}^{\star} \left( \theta^{-1} n^{1/(3\twa)}< 1 +\lvert X_n\rvert \right) \to 1 \qquad \text{ as }\quad\theta \to \infty;
                    \end{equation}
        \item[(b)]  letting $Z_n = \max_{0 \le k \le n} \lvert X_k \rvert$, there exists a subset $\Omega_0\subset \Omega$ with $\Piic(\Omega_0)=1$ such that for all $(\omega,\bar\omega) \in \Omega_0\times\overline\Omega$ and $x \in \iic(\omega)$ there exists $\zeta < \infty$ and $N_x (\omega, \bar \omega)$ such that $P_{\omega}^{x} (N_x < \infty) = 1$ and such that
                    \begin{equation}
                        (\log n)^{-\zeta} n^{1/(3\twa)} \le Z_n (\omega, \bar \omega) \le (\log n)^{\zeta} n^{1/(3\twa)}, \qquad n \ge N_x (\omega, \bar \omega).
                    \end{equation}
    \end{itemize}
\et

For random walk on the IIC backbone $\bbb(\omega)$ we have the following result:
\bt[Random walk on the IIC-backbone, extrinsic distance results]\label{th:IIC-bb}
	If the strong triangle condition is satisfied for some sufficiently small $\beta$, and if either the model is finite-range or Assumption \ref{ass:scaling} holds, then the conclusions of Theorems \ref{th:Eucexittimes} and \ref{th:Eucwalkdist} hold for random walk restricted to $\bbb(\omega)$ when the exponent $6$ is changed to $2\twa$ in Theorem \ref{th:Eucexittimes} and the exponent $1/(3\twa)$ is changed to $1/(2\twa)$ in Theorem \ref{th:Eucwalkdist}.
\et


In \cite[Section 5]{HeyHofHul12a} it is shown that the cluster at the other end of a long edge is small with high probability, so the random walk will with high probability not spend much time outside of $Q_r$ if it exits through a long edge. Thus, the only way a random walk can escape $Q_r$ for more than an instant is if it exits $Q_r$ through the backbone. We conjecture that the time the random walker spends on the other side of long edges is so short that these short excursions will not affect the exit time of the scaling limit (in the standard topologies), so that the scaling limit of the exit time will be proportional to $r^{3\twa}$. We propose a quantity that we believe is interesting to look at as a possible preliminary to studying the scaling limit of random walk on the IIC. We call this quantity the \emph{modified exit time} $\tau_{Q_r}^{\mathsf{ \sss mod}}$ and we define it as the exit time of a random walk that walks on the configuration of the graph that contains all edges \emph{touching} $\iic \cap Q_r$, and where the clock is only stopped if the random walk reaches $Q_r^c$ through the backbone. That is, we do not stop the clock when the random walk exits $Q_r$ through a long edge, but we do force it to return to $\iic \cap Q_r$ in the next step. For this model the exit time is typically much larger than the unmodified version when $\alpha<4$, as the following theorem demonstrates:
\bt[Extrinsic random walk geometry with a modified exit time]\label{th:IIC-mod}
    If Assumption \ref{ass:scaling} holds and the strong triangle condition is satisfied for some sufficiently small $\beta$, then the conclusion of Theorem \ref{th:Eucexittimes} holds for the modified exit time $\tau_{Q_r}^{\mathsf{ \sss mod}}$ when the exponent $6$ is changed to $3\twa$.
\et

\subsubsection*{Results for intrinsic distances}
The following theorem is the analogue of Theorem \ref{th:Eucexittimes} for the intrinsic metric:
\bt[Intrinsic random walk geometry of the IIC]\label{th:Intexittimes}
    Let $r \ge 1$.
    If the strong triangle condition is satisfied for some sufficiently small $\beta$, then the following holds:
    \begin{itemize}
        \item[(a)]  Uniformly in $r$,
                    \begin{equation}
                        \Piic \left(\theta^{-1} r^3 \le E_{\omega}^{0} \tau_{B_r} \le \theta r^3 \right) \to 1 \qquad \text{ as }\quad\theta \to \infty.
                    \end{equation}
        \item[(b)]  \begin{equation}\label{e:EoBrbd}
                        c_3 r^3 \le \Eiic[E_{\omega}^{0} \tau_{B_r} ] \le c_4 r^3.
                    \end{equation}
        \item[(c)]  There exists $\delta_1 < \infty$ and a subset $\Omega_0 \subset \Omega$ with $\Piic(\Omega_0)=1$ such that for all $\omega \in \Omega_0$ and $x \in \iic(\omega)$, there exists $R'_x (\omega)< \infty$ such that
                    \begin{equation}
                        (\log r)^{-\delta_1} r^{3} \le E_{\omega}^{x} \tau_{B_r} \le (\log r)^{\delta_1} r^{3}, \qquad \forall r \ge R'_x (\omega).
                    \end{equation}
        \item[(d)]  For all $\omega \in \Omega_0$ and $x \in \iic(\omega)$, there exists $\delta_2<\infty$ and $R'_x (\omega, \bar \omega)$ such that $P_{\omega}^{x} (R'_x < \infty) = 1$ and such that
                    \begin{equation}
                        (\log r)^{-\delta_2} r^{3} \le \tau_{B_r}(\omega, \bar \omega) \le (\log r)^{\delta_2} r^{3}, \qquad \forall r \ge R'_x (\omega, \bar \omega).
                    \end{equation}
    \end{itemize}
\et

An important geometric property of a graph $\Gamma=(G,E)$ is its \emph{spectral dimension} $d_s (\Gamma)$. It is defined in terms of the asymptotics of the return probability $p_{2n}^{\Gamma} (x,x)$: for any vertex $x \in G$, we set
\begin{equation}
    d_s (\Gamma) = -2 \lim_{n \to \infty} \frac{\log p_{2n}^{\Gamma} (x,x)}{\log n}
\end{equation}
if the limit exists. It is a classical result that $p_n^{\Zd} (x,x)\approx n^{-d/2}$ for random walk on $\Zd$, so that $d_s (\Zd)= d$. Furthermore, it has been proved that for \emph{supercritical percolation} on $\Zd$ the unique infinite cluster $\Ccal_{\infty}$ also  has $d_s (\Ccal_{\infty})=d$ \cite{Barl04}. But the situation is quite different when we consider the IIC.

Alexander and Orbach conjectured that $d_s (\iic)= 4/3$ whenever $d \ge 2$ \cite{AleOrb82}. Although this conjecture is not believed to be true for small $d$, Kozma and Nachmias did prove in \cite{KozNac09} that it holds for percolation models that have $\taupc(x)=c |x|^{d-2} (1+o(1))$. As mentioned before, this asymptotic relation holds for finite-range models, but it does not hold for long-range models with $\alpha <2$. The following theorem improves on their result because it implies that the Alexander-Orbach conjecture is true for any percolation model on $\Zd$ that satisfies the strong triangle condition. In particular, it holds for long-range percolation on $\Zd$ when $d > 3 \twa$.

Again, using the framework of Kumagai and Misumi \cite{KumMis08}, we can establish bounds on the return probability of random walk:
\bt[Asymptotics for the return probability of random walk on the IIC]\label{th:returnprob}
     If the strong triangle condition is satisfied for some sufficiently small $\beta$, then
     \begin{itemize}
        \item[(a)] for $n \ge 1$, uniformly in $n$,
            \begin{equation}
                \Piic \left(\theta^{-1}  \le n^{2/3} p_{2n}^{\omega} (0,0) \le \theta \right)\to 1 \qquad \text{ as }\quad\theta \to \infty.
            \end{equation}
        \item[(b)] there exists $\epsilon < \infty$ and a subset $\Omega_0\subset \Omega$ with $\Piic(\Omega_0)=1$ such that for all $\omega \in \Omega_0$ and $x \in \iic(\omega)$, there exists $N_x (\omega)< \infty$ such that
            \begin{equation}\label{e:pnprecise}
                \frac{(\log n )^{-\epsilon}}{n^{2/3}} \le p_{2n}^{\omega} (x,x) \le \frac{(\log n)^{\epsilon}}{n^{2/3}}, \qquad \forall n \ge N_x(\omega).
            \end{equation}
     \end{itemize}
    \et
Note that \refeq{pnprecise} implies the first limit of \refeq{AlexanderOrbach} in Theorem \ref{th:loglimits}, so the Alexander-Orbach conjecture holds under the strong triangle condition.

For the next theorem we recall the definition of the \emph{annealed law} $\mathsf{P}^{\star}$ in \eqref{e:Pstar}. 

\bt[Intrinsic distance of random walk from the origin]\label{th:Intwalkdist}
    Let $n \ge 1$. If the strong triangle condition is satisfied for some sufficiently small $\beta$, then
    \begin{itemize}
        \item[(a)]  uniformly in $n$,
                    \begin{equation}
                        \mathsf{P}^{\star} \left( d_{\omega} (0, X_n)< \theta n^{1/3}\right) \to 1 \qquad \text{as} \quad \theta \to \infty
                    \end{equation}
                    and
                    \begin{equation}
                        \mathsf{P}^{\star} \left( \theta^{-1}n^{1/3} < 1 + d_{\omega} (0,X_n)\right) \to 1 \qquad \text{as} \quad \theta \to \infty;
                    \end{equation}
        \item[(b)]  letting $Y_n = \max_{0 \le k \le n} d_{\omega} (0, X_k)$, there exists a subset $\Omega_0\subset \Omega$ with $\Piic(\Omega_0)=1$ such that for all $(\omega,\bar\omega) \in \Omega_0\times\overline\Omega$ and $x \in \iic(\omega)$ there exists $\eta <\infty$ and $N'_x (\omega, \bar \omega)$ such that $P_{\omega}^{x} (N'_x < \infty) = 1$ and such that
                    \begin{equation}
                        (\log n)^{-\eta} n^{1/3} \le Y_n (\omega, \bar \omega) \le (\log n)^{\eta} n^{1/3}, \qquad n \ge N'_x (\omega, \bar \omega).
                    \end{equation}
    \end{itemize}
\et

Similar to Theorem \ref{th:IIC-bb}, we have the following result for random walk on the IIC \emph{backbone} for intrinsic distances:
\bt[Random walk on the IIC-backbone, intrinsic distance results]\label{th:IIC-bb2}
    If the strong triangle condition is satisfied for some sufficiently small $\beta$, and if either the model is finite-range or Assumption \ref{ass:scaling} holds, then the conclusions of Theorems \ref{th:Intexittimes} -- \ref{th:Intwalkdist} hold for random walk restricted to $\bbb(\omega)$ 
   when the exponent $3$ is changed to $2$ in Theorem \ref{th:Intexittimes}, 
   the exponent $2/3$ is changed to $1/2$ in Theorem \ref{th:returnprob}, 
   and the exponent $1/3$ is changed to $1/2$ in Theorem \ref{th:Intwalkdist}.
\et
A direct result of the above theorem is that $d_s (\bbb)=1$.

\subsection{Discussion}
\subsubsection*{Literature}
The study of random walk on incipient infinite cluster was initiated by Kesten \cite{Kest86b}. 
Kesten, who initiated the mathematical study of random walk on the IIC proved that random walk on two-dimensional IIC is subdiffusive \cite{Kest86b}. Kozma and Nachmias \cite{KozNac09} established the Alexander-Orbach conjecture for random walk on (finite range) high-dimensional percolation. 
Their work was based on a paper by Barlow, Jarai, Kumagai and Slade \cite{BarJarKumSla08}, who proved the Alexander-Orbach conjecture for random walk on high-dimensional \emph{oriented} percolation. Related results (that predate \cite{BarJarKumSla08}) appear in a study of random walk on the IIC analogue on trees by Barlow and Kumagai \cite{BarKum06}.

In recent work, Jarai and Nachmias \cite{JarNac13} prove bounds on the effective resistance of \emph{branching random walk} in low dimension, thereby settling an open problem formulated in \cite{BarJarKumSla08} . 

Random walk on \emph{finite} critical clusters have been studied in spatial and non-spatial regime by Nachmias and Peres \cite{NacPer08}, Heydenreich and van der Hofstad \cite{HeyHof11} and Croydon, Hambly and Kumagai (in preparation). 

Croydon \cite{Croy08} studies random walk on a random walk trace, and proves that the scaling limit is Brownian motion. We expect a similar behaviour for random walk on the IIC backbone (cf.\ Theorems \ref{th:IIC-bb} and \ref{th:IIC-bb2}).

\subsubsection*{Overview of the proofs}
As mentioned, thanks to the framework of Kumagai and Misumi \cite{KumMis08}, the above theorems all follow once we prove the appropriate volume and effective resistance estimates of the intersection between the IIC (or backbone) and balls. In the next section we state our main technical result, Theorem \ref{th:Goodballs}, which establishes precisely these bounds.
But that Theorem \ref{th:Goodballs} implies Theorems \ref{th:loglimits}, \ref{th:Eucexittimes}, and \ref{th:Eucwalkdist} -- \ref{th:IIC-bb2} is not obvious. 
The main idea of the proofs of Kumagai and Misumi is to first show that simple bounds on the volume and effective resistance of balls imply bounds on, for instance, the exit time, if the configuration is `nice'. They then show that `most' configurations are nice, and thus get bounds that apply with high probability. 

The criteria of  \cite{KumMis08} are not directly applicable to the extrinsic metric case, so we have to make a few modifications (see Remark \ref{rk:KM} below).
As an example of the arguments that are involved, we give our modification of Kumagai and Misumi's proof to show how Theorem \ref{th:Eucexittimes} follows from Theorem \ref{th:Goodballs} in Appendix \ref{app:Pfth14}.

Our proofs add to the existing literature in three ways: 
we prove ball growth and effective resistance bounds in extrinsic geometry; 
we generalize exiting results for intrinsic metric to infinite range models; 
we establish bounds for the IIC backbone.

We prove Theorem \ref{th:Goodballs} in Sections \ref{sec:pftha} and \ref{sec-LB-vol-eff} for the extrinsic case, in Section \ref{AOproof} for the intrinsic case, in Section \ref{sec:Ext-bb} for both intrinsic and extrinsic cases on the backbone, and in Section \ref{sec:Ext-mod} for the modified case.
The proof of Theorem \ref{th:LRP-IIC} is different, so we present it in full in Section \ref{sec:LRP-IIC}.

\section{Definitions, important results and the main theorem}\label{sec:defandimpres}
In this section we state the theorem that implies all the theorems of the introduction (except Theorem \ref{th:LRP-IIC}). We also restate some important results on which our analysis is based and we introduce most of the definitions that we need.

A standard piece of notation are the generic constants $C \ge c > 0$. The value of these constants may change from line to line, or even within the same equation. We will make no attempt to determine their numerical value.

The following definitions apply to general graphs:
\bd[Random walk terminology]\label{def:properties}{\white .}\\
    \vspace{-.5cm}
   \begin{enumerate}
    \item Let $\Gamma =(G,E)$ be a graph with vertex set $G$ and edge set $E$. Define $\mu_x$ as the degree of vertex $x \in G$, and let $\mu(\;\cdot\;)$ be the extension of $\mu_x$ to a measure on $\Gamma$. The discrete-time simple random walk on $\Gamma$ is the Markov chain with transition probabilities
        \begin{equation}
            P^{x}_{\Gamma} (X_1 = y) = \frac{1}{\mu_x} \qquad \text{ if } \{x,y\} \in E.
        \end{equation}
        The \emph{transition density} (or discrete-time heat kernel) is defined as
        \begin{equation}
            p_n^{\Gamma} (x,y) = \frac{P^{x}_{\Gamma} (X_n = y)}{\mu_y}.
        \end{equation}
        Note that $p_n^{\Gamma} (x,y) = p_n^{\Gamma} (y,x)$.
    \item Define the \emph{edge volume of $\;\Gamma$} as
        \begin{equation}
            V(\Gamma) = \sum_{x \in G} \mu_x.
        \end{equation}
   \item  Define the \emph{quadratic form} $\Ecal$ as
        \begin{equation}
            \Ecal (f,g) = \frac12 \sum_{\{x,y\} \in E} (f(x) - f(y))(g(x)-g(y)).
        \end{equation}
        If we think of $\Gamma$ as an electrical network where the edges represent unit resistors, then we can think of $\Ecal (f,f)$ as the \emph{energy dissipation} of the electrical network when the potential at the vertices is given by $f$. We define the set $H^2 = \{f \in \R^{G} \colon \Ecal(f,f) < \infty\}$, so that $H^2$ can be viewed as the set of all `physical' potential functions. For two disjoint subsets $A$ and $B$ of $G$, we define the \emph{effective resistance} between $A$ and $B$ as
        \begin{equation}
            \Reff (A,B)^{-1} = \inf \{\Ecal(f,f) \colon f \in H^2, f \vert_{A} = 1, f \vert_{B}=0\}.
        \end{equation}
        For vertices $x,y \in V$, let $\Reff(x,y) = \Reff(\{x\},\{y\})$ and $\Reff(x,x) =0$. The effective resistance is a metric on subsets of $G$ that is dominated by the intrinsic metric, i.e.,
        \begin{equation}\label{e:reffled}
            \Reff(A,B) \le d_{\Gamma}(A,B).
        \end{equation}
        Many other useful properties of $\Reff(\;\cdot\;,\;\cdot\;)$ are known, cf. \cite{DoySne84}.
   \item Let $\sigma_A$ denote the \emph{random walk hitting time} of the set $A \subseteq G$, i.e.,
        \begin{equation}
            \sigma_A = \inf \{n \ge 0 \colon X_n \in A\}.
        \end{equation}
        The \emph{random walk exit time} $\tau_A$ from the set $A$ is given by $\tau_A = \sigma_{A^c}$.
     \item Let $G_A(x,y)$ denote the \emph{Green's function} of the set $A \subset G$:
     	\begin{equation}
		G_A (x,y) := \frac{1}{\mu_y} E_\Gamma^x [\# \text{ of visits of $X_n$ to $y$ before }\sigma_A]
	\end{equation}
	It is easy to establish that $G_A (x,y) \le G_A (x,x)$ and that $\Reff(x,A) = G_A(x,x).$
   \end{enumerate}
\ed

For (percolation) subgraphs of $\Zd$ we state the following definitions:
\bd[Percolation terminology]\label{def:def}
\label{def-onin2}{\white .}\\
    \vspace{-.5cm}
\begin{enumerate}
\item   Given an edge configuration and a set $A \subseteq \Z^d$, we write $\Ccal(A)$ for the set of vertices that are connected to $A$, i.e., $\Ccal(A)=\{y \in \Z^d: A \conn y\}$. Given an edge configuration and an edge set $B \subseteq E(\Zd)$, we define the restricted cluster $\tilde{\Ccal}^{B}(A)$ as the set of vertices $y \in \Ccal(A)$ that are connected to $A$ in the (possibly modified) configuration in which all edges in $B$ are closed. When $A = \{x\}$ for some $x \in \Zd$, as will often occur, we write $\Ccal(\{x\})=\Ccal(x)$.
\item For any pair $x,y \in \Zd$, we write $\{x,y\}$ for the \emph{undirected edge} between $x$ and $y$, and we write $(x,y)$ for the \emph{directed edge} from $x$ to $y$. When dealing with directed edges $(x,y)$, we call $x$ the `bottom' vertex, and $y$ the `top' vertex. We define $\Be = \{(x,y) : x \in Q_r, y \in \Zd\}$, the set of directed edges with the bottom vertex inside $Q_r$ and the top vertex in $\Zd$.
\item Let $\omega$ be an edge configuration and $b$ an (open or closed) edge. Let $\omega^{b}$ be the edge configuration $\omega$ with the status of the edge $b$ changed. We say an edge $b$ is a \emph{pivotal edge} for the configuration $\omega$ and the event $E$, if $\omega \in E$ and $\omega^{b} \,\nin\, E$, or if $\omega \,\nin\, E$ and $\omega^{b} \in E$. An edge $b$ that is pivotal for a configuration $\omega$ and a connection event $\{A \conn B\}$ will always be assumed to be directed, i.e., $b = (\ulb, \olb)$, in such a way that $\omega, \omega^{b} \in \{A \conn \ulb\}\cap\{\olb \conn B\}$. When we say that an edge is pivotal for an event, what we mean is that it is pivotal for that event in some fixed but unspecified configuration.
\item Given a (deterministic or random) set of vertices $A$ and an event $E$, we say that \emph{$E$ occurs on $A$}, and write $\{E$ on $A\}$, if the (possibly modified) configuration $\omega_A$ in which every edge that does not have both endpoints in $A$ is made vacant satisfies $\omega_A \in E$.  We adopt the convention that $\{x \conn x$ on $A\}$ occurs if and only if $x\in A$.
    Similarly, we say that \emph{$E$ occurs off $A$}, and write $\{E$ off $A\}$, if $\{E$ on $A^c\}$, where $A^c$ is the complement of $A$, and we say that \emph{$E$ occurs through $A$}, and write $\{E $ through $A\}$ for the event $E \setminus \{E $ off $A\}$.
\item Given an edge $e=\{e_1, e_2\}$ (open or closed) and a set of vertices $A \subseteq \Zd$, we say that $e$ \emph{touches} $A$ if either $e_1 \in A$, or $e_2 \in A$ or both.
\item We often abbreviate `pivotal edge' by either `pivotal' or simply by `piv', and similarly, we will often abbreviate `backbone-pivotal edge' by `backbone-pivotal' or `bb-piv'.
\end{enumerate}
\ed

We use the following pieces of notation repeatedly:
\bd[Specific definitions]\label{def:spec}{\white .}\\
    \vspace{-.5cm}
    \begin{itemize}
        \item[(i)] Define the \emph{$r$-restricted cluster} $\bqr(x)$ by
        \begin{equation}\label{e:rrestQ}
            \bqr(x) =\{y \in \Zd \colon x \conn y \text{ on }Q_r\}.
        \end{equation}
        We write $\bqr(0) = \bqr$. We similarly write $\bqr^{\sss \bbb}$ to denote the $r$-restricted cluster of $\bbb$, i.e.,
        \begin{equation}
            \bqr^{\sss \bbb}(x) =\{y \in \Zd \colon x \conn y \text{ on }Q_r \cap \Bb\}.
        \end{equation}
        \item[(ii)] We write $B_r^{\sss \bbb}(x;\omega)$ for the intrinsic ball of radius $r$ centered at $x$ on $\bbb(\omega)$. We write $\Reffbb(A,B)$ for the effective resistance between $A$ and $B$ in $\bbb$. Similarly, we write $\Reff^{\mathsf{\sss mod}}(A,B)$ for the effective resistance in the modified configuration described above Theorem \ref{th:IIC-mod}.
    \end{itemize}
\ed

To state the main theorem of this paper we define the following sets:
\bd[Sets of `good' radii]\label{def:radiussets}
    For $\lambda >1$ define
    \begin{enumerate}
        \item the \emph{extrinsic IIC radius set:}
            \begin{equation}
                \begin{split}
                    J_{{\rm E}}(\lambda) = \Big\{r \in [1,\infty]\,:&\,\lambda^{-1}r^{4} \le V(\bqr) \le \lambda r^{4},\,\lambda^{-1} r^{2} \le \ReffB \le \lambda r^{2}\\
                    & \text{\emph{ and }}\Reff(0,x) \le \lambda r^{2} \text{\emph{ for all }}x \in \bqr\Big\};
                \end{split}
            \end{equation}
        \item the \emph{intrinsic IIC radius set:}
            \begin{equation}
                J_{{\rm I}}(\lambda) = \Big\{r \in [1,\infty]\,:\,\lambda^{-1}r^{2} \le \EVBr \le \lambda r^{2} \text{\emph{ and }}\lambda^{-1} r \le \Reff(0,B_r^c)\Big\};
            \end{equation}
        \item the \emph{extrinsic backbone radius set:}
            \begin{equation}
                \begin{split}
                    J_{{\rm E}}^{\mathsf{ \sss Bb}}(\lambda) = \Big\{r \in [1,\infty]\,:&\,\lambda^{-1} r^{\twa} \le V(\bbr) \le \lambda r^{\twa},\,\lambda^{-1} r^{\twa} \le \Reffbb (0, Q_r^c) \le \lambda r^{\twa}\\
                    & \text{\emph{ and }}\Reffbb(0,x) \le \lambda r^{\twa}\text{\emph{ for all }}x \in \bbr\Big\};
                \end{split}
            \end{equation}
        \item the \emph{intrinsic backbone radius set:}
            \begin{equation}
                J_{{\rm I}}^{\mathsf{ \sss Bb}}(\lambda) = \Big\{r \in [1,\infty]\,:\,\lambda^{-1}r \le \EVBr \le \lambda r \text{\emph{ and }}\lambda^{-1} r \le \Reff(0,B_r^c)\Big\};
            \end{equation}
        \item the \emph{modified extrinsic IIC radius set:}
            \begin{equation}
                \begin{split}
                    J_{{\rm E}}^{\mathsf{ \sss mod}}(\lambda) = \Big\{r \in [1,\infty]\,:&\,\lambda^{-1} r^{2\twa} \le V(\bqr) \le \lambda r^{2\twa},\,\lambda^{-1} r^{\twa} \le \Reffmod (0, Q_r^c) \le \lambda r^{\twa}\\
                    & \text{\emph{ and }}\Reffmod(0,x) \le \lambda r^{\twa}\text{\emph{ for all }}x \in \bqr\Big\};
                \end{split}
            \end{equation}
    \end{enumerate}
\ed
\begin{rk}[About the modifications we make to the definitions of Kumagai and Misumi]\label{rk:KM} (i) Our definition of $\jel$ differs from the one proposed in \cite{KumMis08} in two ways. First, we require an upper bound on $\ReffB$. This bound is not required in \cite{KumMis08} because there the choice of metric is restricted to cases where the upper bound is trivially $r$. Second, we use the $r$-restricted cluster $\bqr$, which is not the extrinsic metric ball. One of the main insights of this paper is that $\bqr$ is much easier to work with than the `correct' ball $Q_r \cap \iic$ and it makes no difference in the proofs, since $\tau_{Q_r} = \tau_{\bqr}$ for any random walk started at $x \in \bqr$ and the edge volume of $\bqr$ is comparable to that of $Q_r \cap \iic$.

(ii) Our definition of $J_{{\rm I}}(\lambda)$ also differs from the one proposed in \cite{KumMis08}, as we make no restriction on $\Reff(0,x)$ for points $x \in B_r(0)$. These resistances are trivially bounded from above by $r$, and this turns out to be sufficient.
\end{rk}

We are now ready to state our main technical theorem:
\bt[Most balls are good]\label{th:Goodballs}
    If the strong triangle condition is satisfied for some sufficiently small $\beta$, then,\\
    (a) if Assumption \ref{ass:onearm} holds, and either the model is finite-range, or Assumption \ref{ass:scaling} holds, then there exist $r^\star (\lambda) = r^{\star} \ge 1$, $\lambda_1 > 1$ and $c_1,c'_1, q_1 > 0$ such that
    \begin{equation}
        \Piic(r \in J_{{\rm E}}(\lambda)) \ge 1 - \frac{c_1}{\lambda^{q_{1}}} \qquad \text{ for all }r \ge r^\star, \lambda \ge \lambda_1
    \end{equation}
    and
    \begin{equation}
        \Eiic[ \Reff(0, Q_r^c) \Vr] \le c'_1 r^{6};
    \end{equation}
    (b) there exist $\lambda_2 > 1$ and $c_2,c'_2, q_2 > 0$ such that
    \begin{equation}
        \Piic(r \in J_{{\rm I}}(\lambda)) \ge 1 - \frac{c_2}{\lambda^{q_{2}}} \qquad \text{ for all }r \ge 1, \lambda \ge \lambda_2
    \end{equation}
    and
    \begin{equation}
        \Eiic[ \Reff(0, B_r^c) \EVBr] \le c'_2 r^{3};
    \end{equation}
    (c) if the model is either finite-range or Assumption \ref{ass:scaling} holds, then there exist $r^\star (\lambda) = r^{\star} \ge 1$, $\lambda_3 > 1$ and $c_3, \, c'_3, q_3 > 0$ such that
    \begin{equation}
        \Piic(r \in J_{{\rm E}}^{\mathsf{ \sss Bb}}(\lambda)) \ge 1 - \frac{c_3}{\lambda^{q_3}} \qquad \text{ for all }r \ge r^\star, \lambda \ge \lambda_3
    \end{equation}
    and
    \begin{equation}
        \Eiic[ \Reffbb(0, Q_r^c) V(\bbr)] \le c'_3 r^{2 \twa};
    \end{equation}
    (d) there exist $\lambda_4 > 1$ and $c_4,c'_4, q_4 > 0$ such that
    \begin{equation}
        \Piic(r \in J_{{\rm I}}^{\mathsf{ \sss Bb}}(\lambda)) \ge 1 - \frac{c_4}{\lambda^{q_{4}}} \qquad \text{ for all }r \ge 1, \lambda \ge \lambda_4
    \end{equation}
    and
    \begin{equation}
        \Eiic[ \Reffbb(0, B_r^c) V(B_r^{\mathsf{ \sss Bb}})] \le c'_4 r^{2};
    \end{equation}
    (e) if Assumption \ref{ass:scaling} holds, then there exist $r^\star (\lambda) = r^{\star} \ge 1$, $\lambda_5 > 1$ and $c_5, \, c'_5, q_5 > 0$ such that
    \begin{equation}
        \Piic(r \in J_{{\rm E}}^{\mathsf{ \sss mod}}(\lambda)) \ge 1 - \frac{c_5}{\lambda^{q_5}} \qquad \text{ for all }r \ge r^\star, \lambda \ge \lambda_5
    \end{equation}
    and
    \begin{equation}
        \Eiic[ \Reffmod(0, Q_r^c) V(\bqr)] \le c'_5 r^{3 \twa}.
    \end{equation}
\et
\begin{rk}
    (i) By \cite[Propositions 1.3 and 1.4, and Theorem 1.5]{KumMis08}, all theorems stated in the introduction except Theorem \ref{th:LRP-IIC} are corollaries to the above theorem. In particular, Theorem \ref{th:Eucexittimes} is a corollary to (a), Theorems \ref{th:Intexittimes}, \ref{th:returnprob} and \ref{th:Intwalkdist} are corollaries to (b), Theorem \ref{th:IIC-bb} is a corollary to (c), Theorem \ref{th:IIC-bb2} is a corollary to (d), and Theorems \ref{th:Eucwalkdist} and \ref{th:IIC-mod} are corollaries to (e).

    (ii) We believe that the method of \cite{KumMis08} cannot be used to prove the analogue of Theorem \ref{th:Eucexittimes} for long-range spread-out percolation with $\alpha <4$ because the method requires that $\Reff(0,Q_r^c)$ and $\Reff(0,x)$ scale as the same power of $r$ for all $x \in \bqr$. But this does not appear to be the case for long-range percolation when $\alpha < 4$. In particular, Lemma \ref{lem:reffLRP} shows that with high probability, $\Reff(0,Q_r^c)$ is at most of order $r^{(4 \wedge \alpha)/2}$, whereas Lemma \ref{lem:Goodpoints1} can easily be modified to show that $\Reff(0,x) = O(r^{\twa})$ for all $x \in \bqr$ with high probability. We do not believe that the latter can be improved to give a bound of the same order as $\Reff(0,Q_r^c)$, since the relatively small value for $\Reff(0,Q_r^c)$ comes from the presence of edges of length $>2 r$ with one end in $Q_r$, whereas, by the definition of $\bqr$, these long edges cannot occur on the path from $0$ to any $x \in \bqr$ (because they are too long to be included in $\bqr$).
\end{rk}

We now describe some established results that we use in the course of the proof of Theorem \ref{th:Goodballs}.

The next theorem states bounds on the expected volume of certain balls. The proof of this theorem relies heavily on Fourier analysis. The techniques presented in \cite{HeyHofHul12a} are straightforward and can be applied to a wide range of similar quantities, but unfortunately the calculations are typically quite long. In the course of the proofs of Lemmas \ref{lem:UpperVolBoundBall} and \ref{lem:boundReffV} below we will use similar bounds, but we will omit their proofs. There, we refer the reader to \cite{HeyHofHul12a} and leave it at that.

Let $\Nbb(r)$ be the number of edges in the backbone of the IIC with the bottom vertex at extrinsic distance at most $r$ from 0, that is, $\Nbb(r)$ is the number of (directed) edges $e=(\ule,\ole)$ with $\ule \in Q_r \cap \iic$ such that $\{0 \conn \ule\}$, $\{e$ open$\}$ and $\{\ole \conn \infty\}$ occur disjointly.
\bt[Cluster and backbone volume bounds, \cite{HeyHofHul12a}]\label{ExpectationBounds}
     If the strong triangle condition is satisfied for some sufficiently small $\beta$, then the following bounds hold:
        \begin{eqnarray}
            \label{e:pcVolBound} c r^{\twa} \le& \Epc [|Q_r  \cap \Ccal(0)|] &\le C r^{\twa};\\
            \label{e:BBVolBound} c'r^{\twa} \le& \Eiic [\Nbb(r)] &\le C' r^{\twa};\\
            \label{e:IICVolBound}c'' r^{2 \twa} \le& \Eiic [|Q_r  \cap \iic|] &\le C'' r^{2\twa}.
        \end{eqnarray}
\et

Note that the expected cluster size can be written in two other, useful ways:
\begin{equation}\label{e:expident}
    \E_p [\lvert Q_r \cap \Ccal(0) \rvert] = \sum_{x \in Q_r} \taup (x) = (\tau_p * \indi_{Q_r})(0),
\end{equation}
where $\indi_{Q_r}$ is the indicator function on $\Zd$ of the set $Q_r$ and ``$*$'' denotes a convolution.

We also restate an important theorem by Kozma and Nachmias \cite{KozNac09}. It gives bounds on the expected size of an intrinsic ball of radius $r$ with respect to the critical percolation measure, and it gives a strong upper bound on the intrinsic one-arm probability of critical percolation. Given a graph $\Gamma = (G, E)$ and a subset $A \subseteq G$, define the random set
\begin{equation}
    \partial B_r (x; A) = \{y \in A \colon d_{\Gamma}(x,y)=r\}
\end{equation}
and the event
\begin{equation}
    H(r;A) = \{\partial B_r (x;A) \neq \emptyset\}.
\end{equation}
Note that $H(r,\omega_p)$ is not an increasing event in $p$. That is, if we have two percolation configurations, $\omega_p$ and $\omega_{q}$ with parameters $p$ and $q$ such that $p < q$, that have been coupled in the standard way (cf.\ \cite{Grim99}), then $H(r;\omega_p) \nsubseteq H(r;\omega_q)$. This is not hard to see, as adding edges to a configuration could actually reduce the length of the shortest path between two points. This makes it difficult to bound the intrinsic one-arm probability $\Pp(H(r;\omega_p))$. Kozma and Nachmias get around this problem by considering instead
\begin{equation}\label{e:koznacintonearm}
    \Gamma_{p} = \sup_{A \subseteq \Zd \times \Zd} \Pp(H(r;A)).
\end{equation}
Clearly, an upper bound on $\Gamma_p$ implies an upper bound on $\Pp(H(r;A))$ for any subgraph $A$.
\bt[Properties of critical percolation clusters \cite{KozNac09}]\label{th:KozNac}
   If the strong triangle condition is satisfied for some sufficiently small $\beta$, then the following bounds hold:
   \begin{eqnarray}
        \label{e:Bound1} \Epc[\lvert B_r(0; \omega) \rvert] &\le& C r;\\
        \label{e:Bound2} \Gamma_{p_c} (r) &\le& C/r.
   \end{eqnarray}
\et
Kozma and Nachmias give the proof in \cite{KozNac09}, where they also prove the accompanying lower bounds (but we won't use the lower bounds here). Sapozhnikov recently presented a short and easy alternative proof of \refeq{Bound1} \cite{Sapo10}.

 An important tool in the upcoming analysis is the \emph{van den Berg-Kesten inequality} (or BK-inequality) \cite{BerKes85},\cite{Grim99}. We say an event $A$ is \emph{increasing} if for any two configurations $\omega$ and $\omega'$ such that $\omega \preceq \omega'$ (that is, any edge that is open in $\omega$ is also open in $\omega'$), $\omega \in A$ implies $\omega' \in A$. Hence, by a standard coupling argument, if $A$ is increasing, then $\Pp(A) \le \P_{p'}(A)$ whenever $p < p'$. 
 
One version of the BK inequality that is valid for infinite lattices is the following. 
Let $F^{\sss A}_i$ $(i\in I)$ and $F^{\sss B}_j$ ($j\in J$) be families of \emph{finite} subsets of the edges $E(\Zd)$. We consider increasing events $A$ and $B$ of the form 
\begin{equation}
	A=\{\exists i\in I:\text{ all bonds in $F^{\sss A}_i$ are open}\} \quad \text{ and } \quad B=\{\exists j\in J:\text{ all bonds in $F^{\sss B}_j$ are open}\}.
\end{equation} 
The \emph{disjoint occurence} of such $A$ and $B$ is given by 
\begin{equation}
	A\circ B=\{\exists i\in I,j\in J:F^{\sss A}_i\cap F^{\sss B}_j=\varnothing, \text{ all bonds in $F^{\sss A}_i$ and $F^{\sss B}_j$ are open}\}.
\end{equation}
The BK-inequality then states
\begin{equation}
    \P_p(A\circ B)\le \P_p(A)\,\P_p(B).
\end{equation}

\section{The backbone limit reversal lemma}
Backbone events are by definition \emph{not} cylinder events, and hence it is a priori not clear whether the limiting scheme that gives $\Piic$ can be reversed. (By `reversing the limit' of $\Piic(E)$ for a given event $E$, we mean that we can express $\Piic(E)$ as the limiting scheme \eqref{e:IICdef} applied to some family of sets $E_x$, $x \in \Zd$.) The aim of this section is to show that in many cases we can.

We call an open edge $b=\{x,y\}\in \Zd$ \emph{backbone-pivotal} when every infinite self-avoiding walk in the IIC starting at the origin uses this edge.

It is not difficult to show that there is an infinite number of backbone-pivotal edges $\Piic$-a.s. Indeed, having a finite number of backbone-pivotal edges implies that there exist at least two disjoint infinite paths from the top of the last backbone pivotal. In Theorem 1.4(ii) of \cite{HofJar04} it is proved that in the finite-range setting this does not happen $\Piic$-a.s. This proof is easily modified to the infinite-range setting.

The backbone-pivotal edges can be ordered as $(b_i)_{i=1}^{\infty}$ so that
every infinite self-avoiding walk starting at $0$ passes through $b_i$ before passing through $b_{i+1}$.
Also, we can think of the backbone-pivotal edges as being \emph{directed} edges $b=(x,y)$, where the direction is such that $\{0\conn x\}$ uses different edges than $\{y\conn \infty\}$.
For a directed edge $b=(x,y)$, we let $\ulb=x$ denote its \emph{bottom}, and $\olb=y$ its \emph{top}.
Writing $b_m$ for the $m$th backbone-pivotal edge, we define
\begin{equation}
     \smi\equiv \tilde\Ccal^{b_m}(0) \setminus \tCcal^{b_{m-1}}(0)
\end{equation}
to be the subgraph of the $m$th ``sausage'' (where, by convention, $\tCcal^{b_{0}}(0) = \varnothing$).

If $0$ is connected to $\dqr$ and there are precisely $n$ open pivotal edges for this connection, then we can again impose an ordering on the open pivotal edges $(b_i)_{i=1}^n$ in such a way that any self-avoiding path from 0 to $\dqr$ passes through $b_i$ before passing through $b_{i+1}$.
If $m \le n$, we let $\smr\equiv \tilde\Ccal^{b_m}(0) \setminus \tCcal^{b_{m-1}}(0)$ and we let $\smr=\varnothing$ whenever $0\nleftrightarrow \dqr$ or $m>n$.

In the same way, we let $\smx\equiv \tilde\Ccal^{b_m}(0)\setminus \tCcal^{b_{m-1}}(0)$ where $b_m$ now is the $m$th open pivotal edge for $\{0\conn x\}$, and $\smx=\varnothing$ if no $m$th pivotal bond exists for the connection $\{0\conn x\}$.

We are interested in events that take place on the first $m$ sausages. To this end, we define
\begin{equation}
    \sni \equiv \bigcup_{i=1}^m \mathsf{S}_{i}^{\sss \infty}, \qquad \snr \equiv \bigcup_{i=1}^m \mathsf{S}_{i}^{\sss (r)}, \qquad \text{ and } \qquad \snx \equiv \bigcup_{i=1}^m \mathsf{S}_{i}^x.
\end{equation}
Note that $\snr$ and $\snx$ may contain fewer than $m$ sausages.

We can similarly define the set of the first $m$ sausages of the backbone graph $\Bb$ as $\Bb_m$, the subgraph of $\Bb$ that contains $0$ when the $m$th pivotal edge is removed.
The prelimit analogue of $\Bb_m$ is constructed as follows: consider first the backbone of $\{0 \conn x\}$, that is, let $\Bb^x$ be the graph that is induced by the set $\{z : \{0 \conn z\} \circ \{z \conn x\}\},$ and define $\Bb_m^x$ as the subgraph of $\Bb^x$ that contains $0$ when the $m$th pivotal edge is removed.
 
Even though events occurring on $\sni$ are not necessarily cylinder events, it is still possible to reverse the IIC-limit for such events, as the following lemma demonstrates.

\bl[Backbone limit reversal lemma] \label{lem:BackboneReversal}\label{lem:BPLR}\label{BPLR}\label{BPRL}
Consider a model such that for all cylinder events $F$, $\Piic(F) = \limp \Q_p (F)$. Then, for any event $E$ of the form
\begin{equation}
    E = \left\{\{b_i\}_{i=1}^n \text{ is occupied}\right\}
\end{equation}
for a set of edges $\{b_i\}_{i=1}^n$, and any $m \in \N$,
\begin{equation}\label{e:LimitReversal}
    \Piic\left(E \text{\emph{ on }}\sni\right) = \limp \frac{1}{\chi(p)} \sum_{x \in \Zd} \Pp\left(\{E \text{\emph{ on }} \snx\}\cap\{0\conn x\}\right),
\end{equation}
and
\begin{equation}\label{e:LimitReversalBackbone}
	\Piic(E  \text{\emph{ on }} \Bb_m)	 = \limp \frac{1}{\chi(p)} \sum_{x \in \Zd} \Pp \left(\{E \on \Bb_m^{x}\} \cap \{0 \conn x\}\right).
\end{equation}
\el
\proof
We only prove \eqref{e:LimitReversal}; the proof of \eqref{e:LimitReversalBackbone} follows almost the exact same strategy.
Fix $m$ throughout the proof. We prove the lemma via a comparison of $\sni$, $\snr$ and $\snx$.
To this end, we define the events
\begin{eqnarray}
    \liR &\equiv& \left\{\omega \colon \sni = \snR\right\};\\
    \lrR &\equiv& \left\{\omega \colon \snr = \snR \text{ and }\snR \text{ contains at least $m$ pivotals for $0 \conn \dqr$}\right\};\\
    \lrx &\equiv& \left\{\omega \colon \snr = \snx \text{ and }\snx \text{ contains at least $m$ pivotals for $0 \conn x$}\right\}.
\end{eqnarray}
We show that it is improbable that these sets are different when we compare them on the same configuration and near the origin. Therefore, we may replace one with the other once we take a suitable limit.

We start by observing that $\liR \subseteq \Lambda^{\sss \infty}_{\sss (R+1)}$ for any $R$, and that $\lrR \subseteq \Lambda^{\sss (r+1)}_{\sss (R)}$ and $\lrx \subseteq \Lambda^{\sss (r+1)}_{\sss x}$ for all $r<R$ and $x \in Q_R^c$.

For any $R$ we can write
\begin{equation}\label{e:splitF}
    \left\{E \text{ on } \sni\right\} = \left(\left\{E \text{ on } \sni \right\} \cap \liR \right)\dot \cup\left(\left\{E \text{ on } \sni\right\} \cap \bliR \right) \equiv F_m^1 (R) \dot \cup F_m^2 (R)
\end{equation}
(we write $\dot \cup$ to indicate that this is the union of two mutually exclusive events).
At the end of the proof we take the limit $R \to \infty$.
In this limit, the event $\bliR$ has probability $0$ under $\Piic$ for the following reasons: The occurrence of $\bliR$ implies that there exists a path from one of the first $m$ sausages to $Q_R^c$ that is disjoint of the backbone. In the limit $R \to \infty$ this implies that there exist \emph{two} disjoint connections to $\infty$ and this event does not occur $\Piic$-almost surely. Indeed, since $\liR \subseteq \Lambda^{\sss \infty}_{\sss (R+1)}$, by monotone convergence,
\begin{equation}
    \limsup_{R \to \infty} \Piic(F_m^2(R)) \le  \lim_{R \to \infty} \Piic\left(\bliR\right) = \Piic\left(\lim_{R \to \infty} \bliR\right) = 0.
\end{equation}

For $F_m^1 (R)$, the occurrence of $\liR$ implies $\{E \text{ on }\sni\} = \{E \text{ on }\snR\}$.
Furthermore, for any $r$ such that $0 < r < R$ we can write
\begin{equation}\label{e:splitG}
    F_m^1 (R) = \left(\left\{E \text{ on } \snR \right\} \cap \liR \cap \lrR \right) \dot \cup\left(\left\{E \text{ on } \snR\right\} \cap \liR \cap \blrR \right) \equiv G_m^1 (R,r) \dot \cup G_m^2 (R,r).
\end{equation}
In the double limit where first $R\to\infty$ and then $r\to\infty$, the probability of $G_m^2 (R,r)$ vanishes as
\begin{equation}\label{e:G2bound}
	\lim_{r \to \infty} \lim_{R \to \infty} \Piic(G_m^2 (R,r))
	\le\lim_{r \to \infty} \lim_{R \to \infty} \Piic(\blrR) = \lim_{r \to \infty} \Piic(\blir) =0.
\end{equation}
Here we again used the argument that in the limit there must exist two disjoint paths to $\infty$.

We can rewrite $G_m^1 (R,r)$ as follows:
\begin{equation}\label{e:splitH}
    G_m^1 (R,r) = \left(\left\{E \text{ on } \snR\right\} \cap \lrR \right) \setminus \left(\left\{E \text{ on } \snR\right\} \cap \lrR \cap \bliR\right) \equiv H_m^1 (R,r) \setminus H_m^2 (R,r).
\end{equation}
Since $H_m^2 (R,r) \subseteq \bliR$ we again have that $\Piic(H_m^2 (R,r)) \to 0$ as first $R \to \infty$ followed by $r\to\infty$.

Now, $H_m^1 (R,r)$ is a cylinder event, so that \eqref{e:IICdef} applies,
\begin{equation}
    \Piic(H_m^1 (R,r)) = \limp \frac{1}{\chi(p)} \sum_{x \in Q_R^c} \Pp\big(H_m^1 (R,r)\cap \{0 \conn x\}\big)
\end{equation}
(where the sum over $x\in Q_R$ vanishes in the $p\nearrow p_c$ limit).

The crucial observation is that  $\lrR \cap \{0 \conn x\} =\lrx \cap \{0 \conn x\}$ for $r<R$ and $x\in Q_R^c$, so that
\begin{equation}
	\left\{E \text{ on } \snR\right\} \cap \lrR \cap \{0 \conn x\}
	=\left\{E \text{ on } \snx\right\} \cap \lrR \cap \{0 \conn x\}.
\end{equation}
It follows that
\begin{equation}\label{e:splitM}
    \begin{split}
        H_m^1 (R,r) \cap \{0 \conn x\} &= \left\{E \text{ on } \snR\right\} \cap \lrR \cap \{0 \conn x\}
        	 = \left\{E \text{ on }\snx\right\}\cap \lrR \cap \{0 \conn x\} \\
                &= \left(\left\{E \text{ on } \snx\right\} \cap \{0 \conn x\}\right) \setminus \left(\left\{E \text{ on } \snx\right\} \cap \blrR \cap \{0 \conn x\}\right)\\
                &\equiv M_m^1 (x) \setminus M_m^2 (R,r,x).
    \end{split}
\end{equation}
For $M_m^2 (R,r,x)$ we note that $\blrR$ is a cylinder event, so that \eqref{e:IICdef} implies
\begin{equation}\label{e:M2bound}
    \begin{split}
        \lim_{r \to \infty}\lim_{R \to \infty}  \limp \frac{1}{\chi(p)} \sum_{x \in Q_R^c} \Pp(M_m^2 (R,r,x))
        &\le \lim_{r \to \infty} \lim_{R \to \infty} \limp \frac{1}{\chi(p)} \sum_{x \in Q_R^c} \Pp(\blrR, 0 \conn x)\\
        &\le \lim_{r \to \infty} \lim_{R \to \infty} \Piic(\blrR) = \lim_{r \to \infty} \Piic(\blir) =0.
    \end{split}
\end{equation}

Combining \refeq{splitF}--\refeq{splitM},
\begin{eqnarray}
    \Piic\left(E  \text{ on }  \sni\right) &= &
    \Piic(F_m^2 (R)) + \Piic(G_m^2 (R,r)) - \Piic(H_m^2 (R,r))\nnb
     &&{}+ \limp \frac{1}{\chi(p)} \sum_{x \in Q_R^c}
     \big(\Pp(M_m^1 (x))  - \Pp(M_m^2 (R,r,x))\big).
\end{eqnarray}
Now we add $0=\limp {\chi(p)^{-1}} \sum_{x \in Q_R} \Pp(M_m^1(x))$
to the right-hand side, so that the term involving $M^1_m(x)$ is independent of $r$ and $R$.
Then we let $R\to\infty$, so that $\Piic(F_m^2 (R))$ and $\Piic(H_m^2 (R,r))$ vanish.
After this we let $r\to\infty$, so that the terms involving $G_m^2 (R,r)$ and $M_m^2 (R,r,x)$ also disappear, by \eqref{e:G2bound} and \eqref{e:M2bound}.
The result is
\begin{equation}
    \Piic\left(E \text{ on } \sni\right) = \limp \frac{1}{\chi(p)} \sum_{x\in \Zd} \Pp(M_m^1 (x)) = \limp \frac{1}{\chi(p)} \sum_{x \in \Zd} \Pp\left(\left\{E \text{ on } \snx\right\}\cap\{0\conn x\}\right),
\end{equation}
completing the proof.
\qed

\section{Upper bounds for the extrinsic case}\label{sec:pftha}
In this section we prove all the upper bounds that are needed to establish Theorem \ref{th:Goodballs}(a). Note that all bounds that we prove in this section also hold for long-range percolation with $\alpha \le 2$ when all occurrences of $r^2$ are replaced with $r^\alpha$, all occurrences of $r^4$ are replaced with $r^{2 \alpha}$ and all occurrences of $r^6$ are replaced with $r^{3 \alpha}$.

\subsubsection*{Effective resistance between $0$ and $Q_r^c$.}
We start by proving that the effective resistance between the origin and the boundary of $Q_r$ is with high probability bounded by $\lambda r^{2}$ under $\Piic$.
 \bl[IIC effective resistance upper bound]\label{lem:Reffub}
    Suppose $\alpha > 2$. Whenever the strong triangle condition is satisfied with sufficiently small $\beta$, there exists a constant $C >0$ such that $\Eiic[\ReffB] \le C r^{2}$. As a result, for all $\lambda>0$
    \begin{equation}
        \Piic \left(\ReffB \ge \lambda r^{2}\right) \le C/ \lambda.
    \end{equation}
\el
\medskip
\proof We would like to use the limiting scheme \refeq{IICdef} to evaluate this bound. This limit is established for cylinder events, but the random variable $\Reff(0, Q_r^c)$ is not necessarily defined in terms of cylinder events when the edge lengths are unbounded. Indeed, when we are dealing with finite-range percolation models, we can use that $\Reff(0,Q_r^c)$ can be determined by inspecting the status of edges with both ends inside $Q_{r+L}$ only, but this is not the case when dealing with long-range percolation, and so it is not immediately clear that the limit in \refeq{IICdef} can be reversed. Of course this is only a technical issue: fixing a large $R$ and closing all edges that are longer than $R$ gives a usable upper bound on $\Reff(0,Q_r^c)$ whose value is measurable with respect to $Q_{r+R}$. Thus we can use the IIC limiting scheme \refeq{IICdef}. We will write down this dependence on $R$ explicitly below.

But before we use the limiting scheme we observe that $\Piic$-a.s. $\ReffB = \ReffB \indi_{\{0 \conn Q_r^c\}}$, so that:
\begin{equation}\label{e:Reffsplit1}
    \begin{split}
        \Eiic [\Reff(0, Q_{r}^c)]       =& \Eiic[\ReffB \indi_{\{0 \conn Q_r^c\}}]\\
                                        =& \limp \frac{1}{\chi(p)} \sum_{x \in \Zd} \E_p [\ReffB \indi_{\{0 \conn x\}}\indi_{\{0 \conn Q_r^c\}}]\\
                                        =& \limp \frac{1}{\chi(p)} \sum_{x \in Q_r} \E_p [\ReffB \indi_{\{0 \conn x\}}\indi_{\{0 \conn Q_r^c\}}]\\
                                         & + \limp \frac{1}{\chi(p)} \sum_{x \in \dqr} \E_p [\ReffB \indi_{\{0 \conn x\}}\indi_{\{0 \conn Q_r^c\}}].
    \end{split}
\end{equation}
We have split up the contributions from $x \in Q_r$ and from $x \notin Q_r$ and we treat them separately.

To bound the first term on the right-hand side we use that $\ReffB\indi_{\{0 \conn Q_r^c\}}$ can be bounded from above by $C_d r^d$, for some constant $C_d$ that depends only on $d$: indeed, the random variable is only non-zero when there is a path from $0$ to $Q_r^c$ and the maximum effective resistance of any such configuration is achieved by a configuration where every vertex lies on the same unique path of open edges from the origin to a single vertex on the boundary of the ball. By the series law for resistances, the effective resistance of such a configuration is proportional to the number of vertices in the ball, and so it can be bounded by $C_d r^d$. Hence,
\begin{equation}\label{e:rhs1}
    \begin{split}
       \limp \frac{1}{\chi(p)} \sum_{x \in Q_r} \E_p [\ReffB \indi_{\{0 \conn x\}}\indi_{\{0 \conn Q_r^c\}}] &\le \limp \frac{1}{\chi(p)} C_d r^d \sum_{x \in Q_r} \tau_p (x)\\
       & \le \limp \frac{1}{\chi(p)} C_d r^d C r^{2}  = 0,
    \end{split}
\end{equation}
where in the second inequality we use the bound from Theorem \ref{ExpectationBounds}, and in the final equality we use that $\chi(p)$ diverges in the limit, while $r$ is fixed.

To bound the second term on the right-hand side of \refeq{Reffsplit1}, first note that for $x \notin Q_r^c$, $\{0 \conn x\} \subseteq \{0 \conn Q_r^c\}$, so we can drop the second indicator function. Recall \refeq{reffled} and note that the intrinsic distance between two connected sets $A$ and $B$ is bounded from above by the number of open edges between them, that is,
\begin{equation}\label{e:dlenoedges}
    d(A,B)\indi_{\{A \conn B\}} \le \sum_{e \in E} \indi_{\{A \conn \ule\}\circ\{e \text{ open}\}\circ\{\ole \conn B\}}.
\end{equation}

A configuration in $\{0 \conn Q_r^c, 0 \conn x\}$ with $x \in Q_r^c$ may contain two or more (partially overlapping) paths from 0 to $Q_{r}^c$. The graph distance from 0 to $Q_{r}^c$ is equal to the length of the shortest such path. Therefore, the effective resistance $\ReffB$ is at most the number of edges from 0 to $Q_{r}^c$ along a path that ends at $x$. Hence we can apply \refeq{reffled} and \refeq{dlenoedges} with $A=\{0\}$ and $B= Q_{r}^c \cap \Ccal(x)$ to bound $\ReffB$:
\begin{equation}\label{e:reffbd}
    \E_p[\ReffB \indi_{\{0 \conn x\}}] \le \sum_{e \in \Be} \E_p [\indi_{\{0 \conn \ule\}\circ\{e \text{ open}\}\circ\{\ole \conn x\}}]
\end{equation}
(recall the definition of $\Be$, Definition \ref{def:def}(ii)).
We get
\begin{equation}
    \limp \frac{1}{\chi(p)} \sum_{x \in \dqr} \E_p [\ReffB \indi_{\{0 \conn x\}}] \le \limp \frac{1}{\chi(p)} \sum_{x \in \dqr} \sum_{e \in \Be} \E_p [\indi_{\{0 \conn \ule\}\circ\{e \text{ open}\}\circ\{\ole \conn x\}}].
\end{equation}
Now we apply the BK inequality to the right-hand side to get the upper bound
\begin{equation}
    \limp \frac{1}{\chi(p)} \sum_{x \in \dqr} \sum_{e \in \Be} \tau_p (\ule) pD(e) \tau(x - \ole).
\end{equation}
If we extend the summation of $x$ to $\Zd$, then we get a factor $\chi(p)$. The summation over $p D(e)$ then leaves us with a factor $p$ and the summation over edges in $\Be$ as a result simplifies to the summation over one end of the edges (we choose the bottom end). We also take the limit $p \ua p_c$ to get the upper bound
\begin{equation}\label{e:rhs2}
   \limp \frac{1}{\chi(p)} \sum_{x \in \dqr} \E_p [\ReffB \indi_{\{0 \conn x\}}] \le p_c \sum_{\ule \in Q_r} \tau_{p_c} (\ule)= p_c \E_{p_c} [\lvert Q_r \cap \Ccal(0) \rvert] \le Cr^{2},
\end{equation}
where the last inequality is due to our volume bound in Theorem \ref{ExpectationBounds}. Combining \refeq{rhs1} and \refeq{rhs2} we conclude
\begin{equation}
    \Eiic[\ReffB] \le C r^{2}.
\end{equation}
Using Markov's inequality with this bound we get
\begin{equation}
    \Piic \left(\ReffB \ge \lambda r^{2}\right) \le \frac{1}{\lambda r^{2}} \Eiic [\ReffB] \le C / \lambda.
\end{equation}
\qed

\subsubsection*{The edge volume of $\bqr$.}
Recall Definitions \ref{def:properties}(ii) and \ref{def:spec}(i). Next we turn to the upper bound on the probability that the edge volume of $\bqr$ is larger than $\lambda r^{4}$.
\bl[IIC edge volume upper bound] \label{lem:UpperVolBoundBall}
    Suppose $\alpha > 2$. If the strong triangle condition is satisfied for some sufficiently small $\beta$, then there exists a $C>0$ such that for all $\lambda>0$,
    \begin{equation}
        \Piic \left(\Vr \ge \lambda r^{4}\right) \le C / \lambda.
    \end{equation}
\el
\proof By Markov's inequality,
\begin{equation}\label{e:MarkovVol}
    \Piic\left(\Vr \ge \lambda r^{4}\right) \le \frac{\Eiic[\Vr]}{\lambda r^{4}}.
\end{equation}
Note that
\begin{equation}\label{e:VcqrIndi}
    \Vr \le V(Q_r \cap \iic) = \sum_{b \in \Be} \indi_{\{0 \conn \ulb, b \text{ open}\}}.
\end{equation}
Therefore,
\begin{equation}
        \Eiic[\Vr] \le \sum_{b \in \Be}  \Eiic[\indi_{\{0 \conn \ulb, b \text{ open}\}}]\\
                    = \limp \frac{1}{\chi(p)} \sum_{b \in \Be} \sum_{x \in \Zd} \P_p(0 \conn \ulb, b \text{ open},0 \conn x)
\end{equation}
The event in $\P_p$ can be contained as follows:
\begin{equation}\label{e:TGvol}
    \{0 \conn \ulb, b \text{ open},0 \conn x\} \subseteq \bigcup_{z \in \Zd} \{0 \conn z\}\circ\{z \conn \ulb, b \text{ open}\}\circ\{z \conn x\}.
\end{equation}
Furthermore,
\begin{equation}\label{e:splitvol}
    \{z \conn \ulb, b \text{ open}\} = \{z \conn \ulb\} \circ \{b \text{ open}\} \cup \{z \conn \olb\} \circ \{b \text{ open}\}.
\end{equation}
Applying \refeq{TGvol}, \refeq{splitvol} and the BK-inequality, we get
\begin{equation}
    \Eiic[\Vr] \le \limp \frac{p}{\chi(p)} \sum_{b \in \Be} \sum_{z,x \in \Zd} \taup(z)[\taup(\ulb-z)+ \taup(\olb-z)]D(b)\taup(x-z).
\end{equation}
Summing over $x$ gives a factor $\chi(p)$. Taking the limit $p \ua p_c$ yields
\begin{equation}\label{e:mod1}
    \begin{split}
        \Eiic[\Vr] &\le p_c \sum_{b \in \Be} \sum_{z \in \Zd} \taupc(z)[\taupc(\ulb-z)+ \taupc(\olb-z)]D(b)\\
                    &= p_c(D* \indi_{Q_r} *\taupc *\taupc)(0) + p_c(\indi_{Q_r}*D*\taupc *\taupc)(0) \le C r^{4}.
    \end{split}
\end{equation}
The second inequality can be established in the same way as was done in the course of the proof of Theorem \ref{ExpectationBounds} (cf.\ \cite[Theorem 1.5]{HeyHofHul12a}). Applying this bound to \refeq{MarkovVol} completes the proof.
\qed

\subsubsection*{An upper bound on $\ReffB \cdot \Vr$.}
The second part of Theorem \ref{th:Goodballs}(a) states an upper bound on the product of the edge volume and the effective resistance of $\bqr$. We prove this bound here.
\bl[Bound on the expectation of the product of the volume and resistance of a ball]\label{lem:boundReffV}
    Let $r \ge 1$, and suppose $\alpha >2$. If the strong triangle condition is satisfied for some sufficiently small $\beta$, then there exists a constant $C > 0$ such that
    \begin{equation}
        \Eiic[\ReffB \Vr] \le C r^{6}.
    \end{equation}
\el
\proof Both $\ReffB$ and $\Vr$ are random variables that can be bounded in terms of indicator functions of two-point function events and of cylinder events (we use the same arguments as we used at the beginning of the proof of Lemma \ref{lem:Reffub}). After we have done so, we may reverse the limit:
\begin{equation}
    \begin{split}
    \Eiic[\ReffB \Vr]   & = \Eiic[\ReffB \Vr \indi_{\{0 \conn Q_r^c\}}] \\
                        & \le \sum_{b \in \Zd \times \Zd} \Eiic[\ReffB \indi_{\{\ulb \in U_r, b \text{ open}\}}\indi_{\{0 \conn Q_r^c\}}]\\
                        & \le \sum_{b \in \Be} \limp \frac{1}{\chi(p)} \sum_{x \in \Zd} \E_p [\ReffB \indi_{\{0 \conn \ulb, b \text{ open}\}} \indi_{\{0 \conn x\}}],
    \end{split}
\end{equation}
where in the second inequality, we used the same arguments as in \refeq{rhs1} to argue that the contributions from $x \in Q_r$ vanish in the limit $p \ua p_c$. And because the summation is now restricted to $x \in Q_r^c$ the indicator $\indi_{\{0 \conn Q_r^c\}}$ was dropped.
Since $x \in Q_r^c$, we can apply the same bound on $\ReffB \indi_{\{0 \conn x\}}$ as was used in \refeq{reffbd}:
\begin{equation}
    \ReffB \indi_{\{0 \conn x\}} \le \sum_{e \in \Be} \indi_{\{0 \conn \ule\}\circ\{e \text{ open}\}\circ\{\ole \conn x\}}.
\end{equation}
Hence,
\begin{equation}
    \Eiic[\ReffB \Vr]    \le \sum_{b \in \Be} \limp \frac{1}{\chi(p)} \sum_{x \in Q_r^c} \sum_{e \in \Be} \E_p[\indi_{\{0 \conn \ule\}\circ\{e \text{ open}\}\circ\{\ole \conn x\} \cap \{0 \conn \ulb, b \text{ open}\}}]
\end{equation}
The event in the indicator function implies that there exists some vertex $z \in \Zd$ such that the path $0 \conn x$ and the path $0 \conn \ulb$ split at $z$ and that $z$ lies either before or after $e$ on the path $0 \conn x$, i.e.,
\begin{multline}
    \bigl\{\{0 \conn \ule\}\circ\{e \text{ open}\}  \circ\{\ole \conn x\}\bigr\} \cap \bigl\{\{0 \conn \ulb, b \text{ open}\}\bigr\}\\
         \subseteq \bigcup_{z \in \Zd} \Bigl\{\bigl\{\{0 \conn z\}\circ\{z \conn \ule\}\circ\{e \text{ open}\}\circ\{\ole \conn x\}\circ \{z \conn \ulb, b \text{ open}\}\bigr\}\\
         \cup \bigl\{\{0 \conn \ule\}\circ\{e \text{ open}\}\circ\{\ole \conn z\}\circ\{\ole \conn x\}\circ \{z \conn \ulb, b \text{ open}\}\bigr\}\Bigr\}.
\end{multline}
Making this replacement, applying \refeq{splitvol} and using the BK-inequality, we get:
\begin{multline}
        \Eiic[\ReffB \Vr] \le \limp \frac{p^2}{\chi(p)} \sum_{x \in \dqr}\sum_{z,\in \Zd} \sum_{e,b \in \Be} \bigl[ \taup(z)\taup(\ule-z) \taup(x-\ule) \\
        + \taup(\ule)\taup(z-\ole)\taup(x-z)\bigr] [\taup(\ulb -z) + \taup(\olb -z)] D(b)D(e).
\end{multline}
Now we can sum over $x$ and take the limit to get
\begin{multline}
        \Eiic[\ReffB \Vr] \le p_c^2 \sum_{z\in \Zd} \sum_{e \in \Be} \bigl[\taupc(z)\taupc(\ule-z) + \taupc(\ule)\taupc(z-\ole)\bigr]\\
       \times \big[(\taupc * \indi_{Q_r})(z) + (\taupc *D* \indi_{Q_r})(z)\big] D(e).
\end{multline}
where we rewrote the terms involving $b$ as a convolution. The factor $D$ dropped out in the first term because, by definition, $\sum_{x \in \Zd}D(x-y)=1$ for any $y$.
Since both convolutions attain their maximum at $z=0$ and since both contribute at most a factor $C r^2$, we have, by Theorem \ref{ExpectationBounds},
\begin{equation}\label{e:mod2}
    \begin{split}
        \Eiic[\ReffB \Vr]    &\le C r^2 \sum_{z \in \Zd} \sum_{e \in \Be} \bigl[\taupc(z)\taupc(\ule-z)D(e) + \taupc(\ule)D(e)\taupc(z-\ole)\bigr]\\
                                &\le C r^2\bigl[(\taupc * \taupc * \indi_{Q_{r}})(0) + (\taupc * D * \taupc * \indi_{Q_{r}})(0)\bigr].
    \end{split}
\end{equation}
We can bound the convolutions using the Fourier-space techniques introduced in \cite{HeyHofHul12a}. The result is that both convolutions can be bounded by $C' r^{4}$, so it follows that
\begin{equation}
    \Eiic[\ReffB \Vr] \le C r^{6},
\end{equation}
as required.
\qed

\subsubsection*{The effective resistance between $0$ and $x \in \bqr$.}
To show that the final upper bound in Definition \ref{def:radiussets}(i) is satisfied we need to show that with a probability proportional to $1-1/\lambda$ there are no vertices $x$ in $\bqr$ such that $\Reff(0,x)$ exceeds $\lambda r^{2}$.
\bl[All vertices in $\bqr$ are well behaved]\label{lem:Goodpoints1}
    Let $\alpha > 2$. If the strong triangle condition is satisfied for some sufficiently small $\beta$, then there exists a $C>0$ such that for all $\lambda>0$
    \begin{equation}\label{e:lemstat}
        \Piic\left(\exists x  \in \bqr \text{\emph{ such that }}\Reff(0,x) \ge \lambda r^{2}\right) \le C/\lambda.
    \end{equation}
\el

\proof We write $m = \lambda r^{2}/2$. Recall the definition of $\bqr^{\sss \bbb}$ in Definition \ref{def:spec}(i). We start by splitting up the event in $\Piic$ in \refeq{lemstat} according to whether $\Sbbr$ is greater than $m$ or not, i.e.,
\begin{multline}\label{e:eventsplit}
       \Piic\left(\exists x  \in \bqr \text{ such that }\Reff(0,x) \ge \lambda r^{2}\right) = \Piic\left(\exists x  \in \bqr \text{ s.t. }\Reff(0,x) \ge \lambda r^{2}, \Sbbr \le m \right)\\
         + \Piic\left(\exists x  \in \bqr \text{ s.t. }\Reff(0,x) \ge \lambda r^{2},\Sbbr  > m \right).
\end{multline}

We can bound the second term on the right-hand side using $\Sbbr \le \Nbb(r)$ and Markov's inequality,
\begin{equation}\label{e:markov}
    \begin{split}
        \Piic\left(\exists x  \in \bqr \text{ s.t. }\Reff(0,x) \ge \lambda r^{2}, \Sbbr > m \right) & \le \Piic\left( \Sbbr  > m \right) \le \Piic(\Nbb(r) > m)\\
        & \le \frac{\Eiic[\Nbb (r) ]}{m} \le \frac{2 C r^{2}}{\lambda r^{2}} \le \frac{C}{\lambda},
    \end{split}
\end{equation}
where the second-to-last inequality follows from Theorem \ref{ExpectationBounds}.

Bounding the first term on the right-hand side of \refeq{eventsplit} is more involved. Since $\Reff(0,x) \le d(0,x)$ the probability of this event is bounded from above by the probability of the event
\begin{equation}\label{e:simplereffevent}
    \{\exists x \in \bqr \text{ s.t. }d(0,x) \ge 2m,\, \Sbbr \le m\}.
\end{equation}
This event implies that $x \in \bqr \setminus \bbr$, since, if $x \in \bbr$, then $d(0,x)$ would be bounded from above by $m$. This means that there exists a (directed) open edge $(u,v)$ such that it is the only open edge on the path from $0$ to $x$ with $u \in \bbr$ and $v \in \bqr \setminus \bbr$. This edge is pivotal for the connection from $0$ to $x$. Furthermore, since $d(0,u) \le m$ and $d(0,x) \ge 2m$, it follows that $d(v,x) \ge m-1$. Finally, since $\{0 \conn u\}$ occurs on $Q_r$ (by the definition of $\bqr$), the vertex $u$ must lie on $\mathsf{Z}_R^{\sss \infty}$, for some $R$ that depends only on $r$ and $d$. Hence, we can contain \refeq{simplereffevent} by the following event:
\begin{equation}
    \bigcup_{u,v \in Q_r} \left\{0 \conn u  \text{ on }\tCcal^{(u,v)}(0) \cap \mathsf{Z}_R^{\sss \infty}, \,(u,v) \text{ open},\, \partial B_{m-1} (v) \neq \emptyset \text{ off }\tCcal^{(u,v)}(0) \cup Q_r^c \right\},
\end{equation}
so that
\begin{multline}\label{e:Piicsplit}
    \Piic(\exists x \in \bqr \text{ s.t. }d(0,x) \ge 2m, \Sbbr \le m) \\
    \le \sum_{u,v \in Q_r} \Piic\bigl(0 \conn u  \text{ on }\tCcal^{(u,v)}(0) \cap \mathsf{Z}_R^{\sss \infty}, (u,v) \text{ open},\\ \partial B_{m-1} (v) \neq \emptyset \text{ off } \tCcal^{(u,v)}(0) \cup Q_r^c\bigr).
\end{multline}
The event inside $\Piic$ on the right-hand side is the intersection between one backbone event and two cylinder events. Thus, by Lemma \ref{lem:BPLR}, it follows that we may reverse the limit for $\Piic$ on the right-hand side of \refeq{Piicsplit} to get
\begin{multline}\label{e:Piiclimrevbbcyl}
    \sum_{u,v \in Q_r}  \Piic\Big(0 \conn u  \text{ on }\tCcal^{(u,v)}(0) \cap \mathsf{Z}_R^{\sss \infty}, (u,v) \text{ open},\, \partial B_{m-1} (v) \neq \emptyset \text{ off } \tCcal^{(u,v)}(0) \cup Q_r^c\Big) \\
    = \limp \frac{1}{\chi(p)} \sum_{y \in \Zd} \sum_{u,v \in Q_r} \Pp \Big(\{0 \conn u\}\circ\{u \conn y\} \text{ on }\tCcal^{(u,v)}(0), (u,v) \text{ open},\\
     \partial B_{m-1} (v) \neq \emptyset \text{ off } \tCcal^{(u,v)}(0) \cup Q_r^c\Big).
\end{multline}

We next use the \emph{Factorization Lemma} (see \cite[Lemma 2.2]{HofHolSla07b}). The variant of the Factorization Lemma that we need states that for two events $E$ and $F$, a vertex $y$ and a directed edge $(u,v)$ with $E \subseteq \{u \in \tCcal^{(u,v)}(y), v \notin \tCcal^{(u,v)}(y)\}$, the following equality holds:
\begin{multline}\label{e:faclem}
    \Pp( E \text{ on }\tCcal^{(u,v)}(y),\, (u,v) \text{ open},\, F \text{ off } \tCcal^{(u,v)} (y))\\
     = p D(u-v) \E_p \left[\indi_{\{E \text{ on }\tCcal^{(u,v)}(y)\}} \P_\Ccal  \left( F \text{ off }\tCcal^{(u,v)}(y)\right)\right],
\end{multline}
where $\P_\Ccal$ denotes that the cluster $\tCcal^{(u,v)}(y)$ is fixed with respect to $\P_{\Ccal}$ (but is random with respect to $\Epc$).

By dropping the restriction that $\{\partial B_{m-1} (v)\neq \emptyset\}$ occurs off $Q_r^c$ on the right-hand side of \refeq{Piiclimrevbbcyl}, we get an upper bound that fits precisely with the formulation of the Factorization Lemma as stated above. Applying \refeq{faclem}, we thus get the upper bound
\begin{multline}\label{e:afterfaclem}
    \limp \frac{1}{\chi(p)} \sum_{y \in \Zd} \sum_{u,v \in Q_r} p D(u -v) \E_p \Bigl[\indi_{\{\{0 \conn u \}\circ\{u \conn y\} \text{ on }\tCcal^{(u,v)}(0)\}}\\
    \times \P_{p,\Ccal} \bigl(\partial B_{m-1} (v) \neq \emptyset \text{ off } \tCcal^{(u,v)}(0)\bigr)\Bigr].
\end{multline}

From \refeq{koznacintonearm} and Theorem \ref{th:KozNac} it follows that, uniformly in $p \in [0,p_c]$ and for any vertex $x \in \Zd$, any set $A \subset \Zd$ and any $n \ge 1$,
\begin{equation}
    \Pp (\partial B_{n}(x) \neq \emptyset \text{ off } A) \le C/n\quad\quad \text{ uniformly in } A.
\end{equation}
Hence, using the above bound, $\Pp (a \conn b $ on $A) \le \Pp (a \conn b)$ and the BK-inequality, we can bound \refeq{afterfaclem} by
\begin{equation}\label{e:BKRplusKozNac}
    \limp \frac{C}{m \chi(p)}\sum_{y \in \Zd} \sum_{u,v \in Q_r} \taup(u) \taup(y-u) p D(u-v).
\end{equation}
We can sum over $y$ to get a factor $\chi(p)$. This cancels the factor $1/\chi(p)$, so we can take the limit $p \ua p_c$ to get the upper bound on \refeq{BKRplusKozNac},
\begin{equation}\label{e:BKRpluslim}
    \frac{C}{m} \sum_{u,v \in Q_r} \taupc(u) p_c D(u-v) \le \frac{C r^{2}}{m} = \frac{C}{\lambda},
\end{equation}
where the final bound follows from Theorem \ref{ExpectationBounds}, $m=\lambda r^2/2$, and the fact that $\sum_x D(x) =1$.

The result is that
\begin{equation}
    \Piic\left(\exists x  \in \bqr \text{ s.t. }\Reff(0,x) \ge \lambda r^{2}, \Sbbr \le m \right) \le C / \lambda,
\end{equation}
and this, combined with \refeq{eventsplit} and \refeq{markov}, completes the proof of Lemma \ref{lem:Goodpoints1}. \qed

\section{Lower bounds for the extrinsic case}\label{sec-LB-vol-eff}
In this section, we prove all the lower bounds that are needed for the proof of Theorem \ref{th:Goodballs}(a). In this analysis we will use the intrinsic-metric ball $B_r$ and of the related ball of pivotal edges, $\Bpiv{r}$: for $x\in \Z^d$ and $A\subset \Z^d$, we let $\Npiv{x,A}$ denote the number of open and pivotal edges for $x\conn A$, where by convention, $\Npiv{x,A}=\infty$ when $x\nconn A$.
We define
    \eqn{
    \lbeq{Bpiv-def}
    \Bpiv{r}=\{x\colon \Npiv{0,x}\leq r\},
    \qquad
    \Vpiv{r}=|\Bpiv{r}|.
    }

The following proposition, the main result in this section, gives lower bounds on $\Reff(0,\BER^c)$ and $\Vr$.
\begin{prop}[Lower bounds on the extrinsic effective resistance and volume]
\label{prop-piv-eff-vol}
Assume that all of the following are satisfied:
\begin{enumerate}
	\item the strong triangle condition is satisfied for some sufficiently small $\beta$;
	\item Assumption \ref{ass:onearm} holds;
	\item the model is finite-range or satisfies Assumption \ref{ass:scaling};
	\item there exist $C', \vep, \eta >0$ such that
    \eqn{
    \lbeq{VIR-bd}
    \Piic\left(\VIR\leq \vep r^{2}\right)\leq C' \vep^{\eta}.
    }
    \end{enumerate}
Then there exists $r^{\star}(\vep) = r^{\star} \ge 1$, $C>0$ and $\kappa>0$ such that, for $r \ge r^{\star}$,
    \eqan{
    \Piic\left(\Reff(0,\BER^c)\leq \vep r^{2}\right) &\leq C\vep^{\kappa};
    \lbeq{eff-res-lbd}\\
    \Piic\left(V(U_r) \leq \vep r^{4}\right) &\leq C\vep^{\kappa}
    \lbeq{vol-lbd}.
    }
\end{prop}
The assumption \refeq{VIR-bd} is proved in Lemma \ref{lem:Cballbounds} below (it is equivalent to \refeq{Bound3}).

\proof[Proof of Theorem \ref{th:Goodballs}(a) subject to Proposition \ref{prop-piv-eff-vol}]
Combining Lemmas \ref{lem:Reffub} -- \ref{lem:Goodpoints1} and Proposition \ref{prop-piv-eff-vol} with $\lambda=1/\vep$ establishes Theorem \ref{th:Goodballs}(a). \qed
\medskip

The proof of Proposition \ref{prop-piv-eff-vol} is organized as follows.
In Section \ref{sec-reduct} we use Assumptions \ref{ass:onearm} and \ref{ass:scaling} to reduce Proposition \ref{prop-piv-eff-vol} to a bound on the number of backbone pivotal edges in $\BER$. This is formulated in Proposition \ref{prop-piv-bd} below. Then, in Section \ref{sec-pf-prop-prop-piv-bd} we use Assumption \ref{ass:scaling} to prove Proposition \ref{prop-piv-bd}.

\subsection{Reduction to the number of backbone-pivotals in an extrinsic ball}
\label{sec-reduct}
We start by bounding the probability that $\Reff(0,\BER^c)$ and $V(U_r)$ are small in terms of the number of pivotals for $\{0\conn \BER^c\}$:
\begin{lemma}[Bounds in terms of the number of pivotals]
\label{lem-piv-eff-vol}
For each $r\geq 1$ and $a\in (0,1)$, the following bounds hold:
    \eqan{
    \Piic\left(\Reff(0,\BER^c)\leq \vep r^{2}\right)
    &\leq \Piic\left(\Npiv{0,\BER^c} \leq \vep r^{2}\right);
    \lbeq{eff-res-piv}\\
    \Piic\left(V(U_r) \leq \vep r^{4}\right)
    &\leq \Piic\left(\Npiv{0,\BER^c} \leq \vep^a r^{2}\right) +\Piic\left(\VI{\vep^a r^{2}}\leq \vep r^{4}\right)
    \lbeq{vol-piv}.
    }
\end{lemma}

\proof The proof of \refeq{eff-res-piv} follows from the bound
    \eqn{
    \lbeq{Reff-piv-bd}
    \Reff(0,\BER^c)\geq \Npiv{0,\BER^c}.
    }
Indeed, the effective resistance of series of elements is the sum of the effective resistances of the elements.
When an edge is pivotal for $\{0\conn \BER^c\}$, then all paths from $0$ to $\BER^c$ must pass through the edge. Thus, we can think of the pivotals, and the intermediate sausages, as lying in series. Since the effective resistance of an edge equals 1, we get \refeq{Reff-piv-bd} by the series law of resistances.

For \refeq{vol-piv}, we note that $V(U_r) \ge | U_r | -1$, so it suffices to prove this bound for $|U_r|$. Also note that
\begin{multline}\lbeq{Qiic-ver-split}
    \Piic\left(\VER\leq \vep r^{4}\right)
    =\Piic\left(\VER\leq \vep r^{4}, \Npiv{0,\BER^c}\leq \vep^a r^{2}\right)\\
    +\Piic\left(\VER\leq \vep r^{4}, \Npiv{0,\BER^c}>\vep^a r^{2}\right).
\end{multline}
The first term can be bounded by the right-hand side of \refeq{eff-res-piv}.
For the second term of \refeq{Qiic-ver-split} we note that if $\Npiv{0,\BER^c}>\vep^a r^{2}$,
then $\VER\geq \VI{\vep^a r^{2}},$ so that
$\VI{\vep^a r^{2}}\leq \VER\leq \vep r^{4}$ holds a.s.
\qed

Now we state the main technical result of this section.
We start by introducing some notation. Recall the definition of $S_n$ above \refeq{Xn-def} and define the \emph{pivotal exit time} of the ball $\BER$ by
    \eqn{
    \lbeq{HR-def}
    H(r)=\inf \left\{n\colon S_{n+1}\in \BER^c \right\}.
    }
Thus, $H(r)$ is the number of backbone-pivotals up to the time at which the process $(S_n)_{n=0}^{\infty}$ leaves the extrinsic ball for the first time.
\begin{prop}[A bound on the number of pivotals]
\label{prop-piv-bd}
If Assumption \ref{ass:onearm} holds and if the strong triangle condition is satisfied for some sufficiently small $\beta$, and furthermore there exists $r^\star (\vep) = r^\star \ge 1$ and $C',q>0$ such that
    \eqn{
    \lbeq{HR-bd}
    \Piic\big(H(r) \leq \vep r^{2}\big)
    \leq C' \vep^q
    }
then there exists $C>0$ and $\gamma>0$ such that, for all $r\geq r^{\star}$,
    \eqn{
    \lbeq{Num-piv-bd-prf}
    \Piic\left(\Npiv{0,\BER^c} \leq \vep r^{2}\right)
    \leq C \vep^\gamma.
    }
\end{prop}
\medskip
The assumption \refeq{HR-bd} is proved in Lemma \ref{prop:piv-exit-time} below. It is a consequence of Assumption \ref{ass:scaling}, but in fact it holds under the considerably weaker assumption that $\limsup_{r \to \infty} \Piic(H(r) \le \vep r^2 ) \le 1-\delta$ for some $\delta>0$. This assumption can be proved for without any knowledge of the scaling limit (but the proof does appear to require a suitable upper bound on the one-arm probability). 

We defer the proof of Proposition \ref{prop-piv-bd} to Subsection \ref{sec-pf-prop-prop-piv-bd}, and now focus on its consequences.
We are ready to prove Proposition \ref{prop-piv-eff-vol}.

\proof[Proof of Proposition \ref{prop-piv-eff-vol} subject to Proposition \ref{prop-piv-bd}]
The bound in \refeq{eff-res-lbd} follows directly from \refeq{eff-res-piv} and Proposition \ref{prop-piv-bd} with $C'=C$ and $\kappa=\gamma$.
For the proof of \refeq{vol-lbd}, we use \refeq{vol-piv}.
The first term in \refeq{vol-piv} can be bounded using Proposition \ref{prop-piv-bd}: for each $a\in (0,1)$,
    \eqn{
    \Piic\left(\Npiv{0,\BER^c} \leq \vep^a r^{2}\right)
    \leq C\vep^{a \gamma}.
    }
To bound the second term in \refeq{vol-piv} we need to show that there exists $a\in(0,1)$ for which we can find a $\kappa>0$ such that
    \eqn{
    \Piic\left(\VI{\vep^a r^{2}}\leq \vep r^{4}\right)
    \leq C\vep^{\kappa}.
    }
Rewriting with $r_\vep=\vep^{a/2}r$ yields
    \eqn{
    \lbeq{VI-bd}
    \Piic\left(\VI{\vep^a r^{2}}\leq \vep r^{4}\right)
    =\Piic\left(\VI{r_\vep^{2}}\leq \vep^{1-2a} r_\vep^{4}\right)\leq C\vep^{\eta(1-2a)},
    }
where the bound follows from the assumption \refeq{VIR-bd}.

By \refeq{VI-bd}, the second term on the right-hand side is bounded by $C\vep^{\eta(1-2a)}$ when $a\in (0,1/2)$.
Thus, for any $a\in (0,1/2)$, \refeq{vol-lbd} follows with
$\kappa=\min\{a \gamma, \eta(1-2a)\}>0$.
\qed

\subsection{Proof of Proposition \ref{prop-piv-bd}}
\label{sec-pf-prop-prop-piv-bd}

We start by relating $\Npiv{0,\BER^c}$ to the number of edges that are pivotal for $0 \conn \dqr$ that are inside a smaller ball.
We fix $a\in (0,1)$ and bound
    \eqn{
    \lbeq{Ng-def}
    \Npiv{0,\BER^c}\geq \Ngood(\vep^a r)
    \equiv \# \left\{e\in Q_{\vep^a r} \colon e \text{ occupied and pivotal for } 0\conn \BER^c \right\}.
    }

Recall the definition of $H(r)$ in \refeq{HR-def} above.
Since any pivotal for $\{0\conn \BER^c\}$ is also backbone-pivotal, we can split, for any $a\in (0,1)$,
    \eqn{
    H(\vep^a r)=\Ngood(\vep^a r)+\left(H(\vep^a r) - \Ngood(\vep^a r)\right) \equiv \Ngood(\vep^a r) + \Nbad(\vep^a r),
    }
By definition, $\Nbad(\vep^a r)$ is the number of edges $e$ with $\bbe \in Q_{\vep^a r}$ that are backbone-pivotals such that all earlier backbone-pivotals are also in $Q_{\vep^a r}$, but that are not pivotal for $\{0\conn Q_{r}^c\}$. Clearly, $\Nbad(\vep^a r)\geq 0$, but the idea is to show that not many pivotal edges are `bad'.\footnote{This argument fails for long-range spread-out percolation with $\alpha <4$, as the number of `bad' pivotals is \emph{not} small. In this regime the probability of having long open edges (say longer than $2r$) in the IIC near the origin is not small, and the presence of such an edge will render all of the backbone-pivotals beyond this point `bad'.}

We can bound
    \eqan{
    \lbeq{Npivsplit}
    \Piic\left(\Npiv{0,\BER^c} \leq \vep r^{2}\right)
    &\leq
    \Piic\left(\Ngood(\vep^a r) \leq \vep r^{2}\right)\\
    &\leq
    \Piic\left(H(\vep^a r) \leq 2\vep r^{2}\right)
    +\Piic\left(\Nbad(\vep^a r)\geq \vep r^{2}\right).\nn
    }
We start by bounding $\Piic\left(\Nbad(\vep^a r)\geq \vep r^{2}\right)$.
Using Markov's inequality gives
    \eqn{
    \lbeq{fmm-Nbad}
    \Piic\left(\Nbad(\vep^a r)\geq \vep r^{2}\right)
    \leq \frac{\Eiic[\Nbad(\vep^a r)]}{\vep r^{2}}.
    }
It is not hard to see that for any edge $e$ that counts toward $\Nbad(\vep^a r)$ there exists a vertex $z\in Q_{\vep^a r}$ such that the event
\begin{equation}
    \{0\conn z\}\circ \{z\conn \BER^c\}\circ \{z\conn \bbe\}\circ \{e\text{ open}\}\circ \{\te \conn \infty\}
\end{equation}
occurs. Thus,
    \eqn{
    \lbeq{expec-N-H}
    \Eiic[\Nbad(\vep^a r)]
    =\sum_{e \in \Bear} \sum_{z\in Q_{\vep^a r}}
    \Piic\left(\{0\conn z\} \circ \{z\conn \BER^c\}\circ \{z\conn \bbe\}\circ \{e\text{ open}\}\circ \{\te \conn \infty\}\right).
    }
By Lemma \ref{lem:BPLR} and the BK-inequality,
\begin{equation} \lbeq{expec-N-H-2}
    \begin{split}
    \Eiic[\Nbad(\vep^a r)]
    &\leq \sum_{e \in \Bear} \sum_{z\in Q_{\vep^a r}} \limp \frac{1}{\chi(p)} \sum_{x \in \Zd} \Ppc\left(\{0\conn z\} \circ \{z\conn \BER^c\}\circ \{z\conn \bbe\}\circ \{e\text{ open}\} \circ \{\ole \conn x\}\right)\\
    & \le \sum_{e \in \Bear} \sum_{z\in Q_{\vep^a r}} \limp \frac{1}{\chi(p)} \sum_{x \in \Zd} \taupc(z) \Ppc(z \conn \dqr) \taupc(\ule -z) p D(e) \taupc(x-\ole)\\
    &\leq \sum_{e \in \Bear} p_c D(e)\sum_{z \in Q_{\vep^a r}} \tau_{p_c}(z)\tau_{p_c}(\bbe-z)
    \P_{p_c}(z\conn \BER^c).
    \end{split}
\end{equation}
Since $z\in Q_{\vep^a r}$,  we have by Assumption \ref{ass:onearm} that uniformly in $z$,
    \eqn{
    \P_{p_c}(z\conn \BER^c)\leq \P_{p_c}(0\conn Q_{r/2}^c)\leq C r^{-2}.
    }
Therefore, we end up with
    \eqan{
    \lbeq{expec-N-H-3}
    \Eiic[\Nbad(\vep^a r)]
    &\leq C p_c r^{-2} \sum_{\ule \in Q_{\vep^a r}} (\tau_{p_c}*\tau_{p_c})(\ule) \le C r^{-2} (\vep^a r)^{4}.
    }
where the last inequality follows from Theorem \ref{ExpectationBounds}.
Thus, we arrive at
    \eqan{
    \lbeq{expec-N-H-4}
    \Eiic[\Nbad(\vep^a r)]
    &\leq C \vep^{4a} r^{2} ,
    }
so that by \refeq{fmm-Nbad}
    \eqn{
    \lbeq{Nbadbd}
    \Piic\left(\Nbad(\vep^a r)\geq \vep r^{2}\right)
    \leq C \vep^{4a-1},
    }
which satisfies \refeq{Num-piv-bd-prf} with $\gamma=4a-1>0$ when $a > 1/4$. Applying \refeq{Nbadbd} to \refeq{Npivsplit} with $a \in (1/4, 1/2)$ completes the proof of Proposition \ref{prop-piv-bd}. \qed

\subsection{A bound on the pivotal exit time}
Next, we investigate the assumption in Proposition \ref{prop-piv-bd} that $\Piic(H(\vep^a r) \leq 2\vep r^{2}) \le C' \vep^q$. The following lemma establishes this bound. The lemma is stated in two ways: in (a) it is stated for high-dimensional finite-range models, where the claim can be proved without assumptions, and in (b) it is stated under Assumption \ref{ass:scaling} for the more general setting that includes long-range spread out percolation models. With Assumption \ref{ass:scaling} we can use the lemma in the proofs of Theorem \ref{th:Goodballs}(c) and (e) below as well. We thus include the exponent $\alpha$ in the statement and the proof, but keep in mind that the result also holds for finite-range models if we set $\alpha$ to $\infty$. Also note that in the statement of Proposition \ref{prop-piv-bd} it suffices to take $r\geq r^\star$, where $r^\star$ may depend on $\vep$.

\bl[A lower bound on the pivotal exit time]\label{prop:piv-exit-time}\color{white}.\color{black}

(a) Assume that
    \begin{equation}\label{eq-twoptass}
        \taupc(x) \le \frac{C}{|x|^{d-\twa}}
    \end{equation}
    and
    \begin{equation}\label{eq-boundaryass}
         \sum_{\ule \in Q_r}\sum_{\ole \in \dqr} \Ppc(0 \conn \ule \text{ in }Q_r)D(\ole-\ule) \le C'.
    \end{equation}
    Then
    \begin{equation}
        \limsup_{r \to \infty} \Piic (H(r) \le \vep r^{\twa}) \le 1-\delta
    \end{equation}
    for $\vep > 0$ and some $0 < \delta < 1.$
    
(b) If Assumption \ref{ass:scaling} holds and if the strong triangle condition is satisfied for some sufficiently small $\beta$, then there exists $r^\star = r^\star (\vep) \ge 1$ and $C, q>0$ such that for $a < 1/\twa$,
\begin{equation}
    \Piic(H(\vep^a r) \leq \vep r^{\twa}) \le C \vep^q.
\end{equation}
\el
The assumption \eqref{eq-twoptass} has been discussed already in the introduction. It is known to hold for high-dimensional finite-range models (cf.\ \cite{HarHofSla03} and \cite{Hara08}), and for a certain class of long-range models \cite{CheSak12}, but it is not known to hold under the strong triangle condition.

The assumption \eqref{eq-boundaryass} has been proved for finite-range models by van der Hofstad and Sapozhnikov \cite[Theorem 1.6]{HofSap11}. It is not known to hold for long-range models.

\proof[Proof of (a)] For a non-negative integer-valued random variable $X$ and probability measure $\P$ we have the elementary inequality
\begin{equation}
    \P(X \ge 1) \ge \frac{\E[X]^2}{\E[X^2]}.
\end{equation}
We apply this inequality as follows:
\begin{equation}\label{e:secmommeth}
    		\Piic (H(r) \ge \vep r^{\twa}) = \Piic (H(r) \indi_{\{H(r) \ge \vep r^{\twa} \}} \ge 1)
		 \ge \frac{\Eiic[H(r) \indi_{\{H(r) \ge \vep r^{\twa} \}}]^2}{\Eiic[H(r)^2 \indi_{\{H(r) \ge \vep r^{\twa} \}}]}.
\end{equation}
Proving the proposition is equivalent to proving that the right-hand side of \refeq{secmommeth} has a uniform, positive lower bound. To achieve this, we bound the expectations on the right-hand side of \refeq{secmom} separately.

We start with an upper bound on the denominator. Trivially,
\begin{equation}
    \Eiic[H(r)^2 \indi_{\{H(r) \ge \vep r^{\twa} \}}] \le \Eiic[H(r)^2].
\end{equation}
We note that
\begin{equation}
    H(r) \le \#\{b \in \Be : b \text{ is bb-piv}\} = \sum_{b \in \Be} \indi_{\{b \text{ is bb-piv}\}}.
\end{equation}
Hence, we have
\begin{equation} \label{e:secmom}
        		\Eiic[H(r)^2]    \le \sum_{b_1, b_2 \in \Be} \Eiic[\indi_{\{b_1 \text{ is bb-piv}\}} \indi_{\{b_2 \text{ is bb-piv}\}}]
		 = \sum_{b_1, b_2 \in \Be} \Piic (b_1, b_2 \text{ are bb-piv}).
\end{equation}
We can apply Lemma \ref{BPLR} to the right-hand side of \refeq{secmom}:
\begin{equation}
    \Eiic[H(r)^2] \le \limp \frac{1}{\chi(p)} \sum_{x \in \Zd} \sum_{b_1, b_2 \in \Be} \Pp( b_1, b_2 \text{ are piv for }0 \conn x).
\end{equation}
The event on the right-hand side can be contained in a disjoint union of events:
\begin{equation}
        \{b_1, b_2 \text{ are piv for }0 \conn x\} 
        \subset \left\{\{0 \conn \ulb_1\}\circ\{\olb_1 \conn \ulb_2\}\circ\{\olb_2 \conn x\}\right\}
         \cup\left\{\{0 \conn \ulb_2\}\circ\{\olb_2 \conn \ulb_1\}\circ\{\olb_1 \conn x\}\right\}.
\end{equation}
Making this replacement and applying the BK-inequality, we obtain an upper bound:
\begin{equation}
    \Eiic[H(r)^2] \le \limp \frac{1}{\chi(p)} \sum_{x \in \Zd} \sum_{b_1, b_2 \in \Be} [\taup(\ulb_1)\taup(\ulb_2 -\olb_1)\taup(x -\olb_2)
     + \taup(\ulb_2) \taup(\ulb_1 - \olb_2)\taup(x - \olb_1)].
\end{equation}
Summing over $x$ and summing the two terms over $\olb_2$ and $\olb_1$, respectively, and then taking the limit, we obtain
\begin{equation}
    \Eiic[H(r)^2] \le p_c \sum_{b_1 \in \Be} \sum_{\ulb_2 \in Q_r} \taupc(\ulb_1)\taupc(\ulb_2 - \olb_1) 
    + p_c \sum_{b_2 \in \Be} \sum_{\ulb_1 \in Q_r} \taupc(\ulb_2)\taupc(\ulb_1 - \olb_2).
\end{equation}
Both sums can be bounded using the Fourier-space techniques described in the proof of Theorem \ref{ExpectationBounds}. For a constant $c_a > 0$ we obtain
\begin{equation}\label{e:secmombd}
    \Eiic [H(r)^2] \le c_a r^{2\twa}.
\end{equation}

We are left to find a lower bound on the numerator of \refeq{secmommeth}. We start by noting
\begin{equation}
    \Eiic[H(r) \indi_{\{H(r) \ge \vep r^{\twa}\}}] = \Eiic[H(r)] - \Eiic[H(r) \indi_{\{H(r) < \vep r^{\twa}\}}].
\end{equation}
We use a trivial upper bound for the second term:
\begin{equation}
    \Eiic[H(r) \indi_{\{H(r) < \vep r^{\twa}\}}] \le \vep r^{\twa}.
\end{equation}
The lower bound on $\Eiic[H(r)]$ is the most involved part of the proof. Let $\Fe = \{(x,y) \colon x \in \dqr, y \in \Zd\}$. Recall that we can order the backbone pivotals from the origin outward. We say that two backbone pivotal edges $e$ and $b$ are \emph{ordered} if $e$ comes before $b$ in this ordering. We start by observing that for $\eta \in (0,1)$:
\begin{equation}
    \begin{split}
        H(r)   &=  \#\{b \in \Be : b \text{ is bb-piv and } \nexists e \in \Fe \text{ s.t. }e,b \text{ are ordered bb-piv}\} \\
                & \ge \#\{ b \in \Beh : b \text{ is bb-piv and } \nexists e \in \Fe \text{ s.t. } e,b \text{ are ordered bb-piv}\}\\
                & = \#\{ b \in \Beh: b \text{ is bb-piv}\} - \#\{b \in \Beh: \exists e \in \Fe \text{ s.t. } e,b \text{ are ordered bb-piv}\}\\
                & \equiv \Nbb (\eta r) - \Zer.
    \end{split}
\end{equation}
Hence,
\begin{equation}
    \Eiic[H(r)] \ge \Eiic[\Nbb (\eta r)] - \Eiic[\Zer].
\end{equation}

Using Theorem \ref{ExpectationBounds}, we can bound
\begin{equation}
    \Eiic[\Nbb (\eta r) ]\ge c (\eta r)^{\twa}.
\end{equation}

Finally, then, we need an upper bound on $\Eiic[\Zer]$. We write
\begin{equation}
    \Zer = \sum_{e \in \Fe} \sum_{b \in \Beh} \indi_{\{e,b \text{ are ord. bb-piv}\}}.
\end{equation}
It follows that
\begin{equation}
    \Eiic [\Zer] = \sum_{e \in \Fe} \sum_{b \in \Beh} \Piic(e,b \text{ are ord. bb-piv}).
\end{equation}
We can apply Lemma \ref{BPLR}:
\begin{equation}\label{e:QiicZer}
    \Eiic[\Zer] = \limp \frac{1}{\chi(p)} \sum_{x \in \Zd} \sum_{e \in \Fe} \sum_{b \in \Beh} \Pp(e,b \text{ are ord. piv for }0 \conn x).
\end{equation}
The events on the right-hand side can be contained as follows: the fact that $e$ is pivotal and comes before $b$ along the path from $0$ to $x$ means that the path from $0$ to $x$ leaves $Q_r$ (using a lexicographical ordering to break ties). Let $y \in \Be$ be the first edge along the path with $\uly \in Q_r$ and $\oly \in Q_r^c$. There has to be a path from $\oly$ back to $\ulb$ as well (this path then containing $e$), and there has to be a path from $\olb$ to $x$. These three paths are disjoint.
So, ignoring the position of the edge $e$, we can contain the event as follows:
\begin{equation}
    \bigcup_{e \in \Fe} \{e,b \text{ are ord. piv for }0\conn x\} 
    \subseteq \bigcup_{\uly \in Q_r}\bigcup_{\oly \in Q_r^c} \{0 \conn y\}\circ\{y \text{ open}\}\circ\{\oly \conn \ulb\}\circ\{\olb \conn x\}.
\end{equation}
Making this replacement and applying the BK-inequality, we obtain
\begin{equation}
    \Eiic[\Zer] 
    \le \limp \frac{1}{\chi(p)} \sum_{x \in \Zd}\sum_{\uly \in Q_r} \sum_{\oly \in Q_r^c} \sum_{b \in \Beh} \Pp(0 \conn y \text{ in }Q_r)
     p D(\oly-\uly) \taup(\ulb-\oly)\taup(x-\olb).
\end{equation}
Summing over $x$ and $\olb$ and taking the limit, we obtain
\begin{equation}
    \Eiic[\Zer] \le p_c^2 \sum_{\uly \in Q_r} \sum_{\oly \in Q_r^c} \sum_{\ulb \in Q_{\eta r}} \Ppc( 0 \conn \uly \text{ in }Q_r)D(\oly - \uly) \taupc(\ulb - \oly).
\end{equation}

All pairs $\ulb$ and $\oly$ are at least at distance $r-\eta r$, and there are $C_d (\eta r)^d$ vertices in $Q_{\eta r}$, so we can apply  \eqref{eq-twoptass} followed by \eqref{eq-boundaryass} to bound
\begin{equation}
	\begin{split}
    		\Eiic[\Zer] \le & p_c^2 C_d (\eta r)^d \frac{1}{(r-\eta r)^{d-\twa}}
		 \sum_{\uly \in Q_r}\sum_{\oly \in Q_r^c} \Ppc(0 \conn \uly \text{ in }Q_r)D(\oly - \uly)\\
				 \le & \tilde{c} \frac{\eta^d}{(1-\eta)^{d-\twa}} r^{\twa}.
	\end{split}
\end{equation}
We end up with
\begin{equation}\label{e:firmombd}
    \begin{split}
        \Eiic[H(r) \indi_{\{H(r) \ge \vep r^\twa\}}]^2    &\ge\bigl(\Eiic[\Nbb  (\eta r)]-\Eiic[\Zer]
        -\Eiic[H(r) \indi_{\{H(r) < \vep r^{\twa}\}}]\bigr)^2 \\
                        &\ge \left(c \eta^{\twa} -\tilde{c} \frac{\eta^d}{(1-\eta)^{d-\twa}} - \vep\right)^2 r^{2 \twa} 
                        = c_b r^{2 \twa}
    \end{split}
\end{equation}
when we choose $\eta>0$ sufficiently small.

Plugging \refeq{secmombd} and \refeq{firmombd} into \refeq{secmommeth}, we obtain
\begin{equation*}
    \Piic(H(r) \ge \vep r^{\twa}) \ge \frac{c_b}{c_a} >0.
\end{equation*}
This completes the proof of (a).
\qed
\color{black}

\proof[Proof of (b)]
Recall \refeq{Xn-def} and note that if $X_n$ converges to Brownian motion or stable motion, then
    \eqn{
    \label{e:weakconv}
    r^{-\twa} H(r) \convd E_1,
    }
where $E_1$ is the exit time of Brownian motion or stable motion from the unit ball in ${\mathbb R}^d$. Let $r_\vep = \vep^a r$. By Assumption \ref{ass:scaling}, \refeq{weakconv} and the Portmanteau Theorem (cf.\ \cite[Theorem 2.1]{Bill99})
\begin{equation}
    \limsup_{r \to \infty}\,  \Piic\left(H(r_\vep) \le \vep^{1-a \twa} r_\vep^{\twa}\right) \le \P \left(E_1 \le \vep^{1-a \twa}\right).
\end{equation}

Hence, for any $\delta > 0$ there exists $r_\vep^\star(\delta)$ such that, for all $r_\vep \ge r_{\vep}^\star(\delta)$,
\begin{equation}
    \Piic\left(H(r_\vep) \le \vep^{1-a \twa} r_{\vep}^{\twa}\right) \le \delta + \P\left(E_1 \le \vep^{1-a \twa}\right).
\end{equation}
Let $X_t \cdot e_i $ be the projection of $X_t$ onto its $i$th coordinate. Observe that by the reflection principle
\begin{equation}
    \P(E_1 \le t) \le 2\, \P(\lvert X_t \rvert \ge 1) \le 2\, \P(\max_i \lvert X_t \cdot e_i \rvert \ge 1/d) \le 2 d\, \P (\lvert X_t \cdot e_1 \rvert \ge 1/d).
\end{equation}
Finally, since $X_t \cdot e_1$ is a one-dimensional Brownian motion or stable motion, we have that $\P(\lvert X_t \cdot e_1 \rvert \ge 1/d) \le C_{X} t^{\twa}$, where $C_{X}$ is a constant, so that
\begin{equation}
    \Piic\left(H(r_\vep) \le \vep^{1-a \twa} r_{\vep}^{\twa}\right) \le \delta + C' \vep^{\twa(1-a\twa)}.
\end{equation}
Let $q = \twa(1-a \twa)$. Then, setting $\delta = \vep^q$ and $r^\star(\vep) = \vep^{-a} r_\vep^\star(\vep^q)$ proves the claim.  \qed

\section{Intrinsic distances for the IIC: proof of Theorem \ref{th:Goodballs}(b)}\label{AOproof}
Our proof of Theorem \ref{th:Goodballs}(b) is a slight adaptation of the proofs in \cite[Section 2]{KozNac09} so we only discuss the changes needed and refer the reader to \cite{KozNac09} for the full details of the proof. Define $\partial B_r (x ; \omega)$ as the boundary of $B_r(x ; \omega)$, that is,
\begin{equation}
    \partial B_r (x; \omega) = \{y \in \Zd \colon d_{\omega}(x,y)=r\}.
\end{equation}
We write $\{x \stackrel{r}{\longleftrightarrow} y\}$  for the event that there exists a path of at most $r$ open edges connecting the vertices $x$ and $y$.

In this section we write $B_r$ for $B_r (0;\omega)$ and $\partial B_r$ for $\partial B_r (0; \omega)$.

\bl[Typical volume and effective resistance of an intrinsic-metric ball]\label{lem:Cballbounds}
    If the strong triangle condition is satisfied for some sufficiently small $\beta$, then Theorem \ref{th:KozNac} implies for $\lambda > 1$,
    \begin{eqnarray}
        \label{e:Bound3} \Piic\big( \lambda^{-1} r^2 \le V(B_r) \le \lambda r^2 \big) &\ge& 1 - c / \lambda ;\\
        \label{e:Bound4} \Piic\big( \Reff(0, \partial B_r) \ge \lambda^{-1} r \big) &\ge& 1 - c/ \sqrt{\lambda}.
    \end{eqnarray}
\el

\proof[Proof of \refeq{Bound3}]
Unless stated otherwise, all sums below are taken over $\Zd$.

We can write $\Eiic[V(B_r)]$ as a sum over edges:
\begin{equation}
    \Eiic[V(B_r)] = \sum_{e \in \Zd \times \Zd} \Piic\left(0 \stackrel{r}{\longleftrightarrow} \ule, e \text{ open}\right).
\end{equation}
Recall the definition of $\zri$ and $\zrx$ in Definition \ref{def:spec}(iii). Since $\{0 \stackrel{r}{\longleftrightarrow} \ule, e \text{ open}\}$ is measurable with respect to $ \mathsf{Z}_{r+1}^{\sss \infty}$ we may reverse the IIC limit by Lemma \ref{lem:BPLR}, so that
\begin{equation}
    \begin{split}
        \Eiic[V(B_r)] &=  \sum_{e \in \Zd \times \Zd} \Piic\left(0 \stackrel{r}{\longleftrightarrow} \ule, e \text{ open}\right)\\
         & = \sum_{e \in \Zd \times \Zd} \Piic\left(\left\{0 \stackrel{r}{\longleftrightarrow} \ule, e \text{ open}\right\} \text{ on }\mathsf{Z}_{r+1}^{\sss \infty}\right)\\
         & = \sum_{e \in \Zd \times \Zd} \limp \frac{1}{\chi(p)} \sum_{x} \P_p \left(\left\{0 \stackrel{r}{\longleftrightarrow} \ule, e \text{ open}\right\} \text{ on }\mathsf{Z}_{r+1}^{x}, 0 \conn x\right)\\
         &= \limp \frac{1}{\chi(p)} \sum_x \E_p [V(B_r)\indi_{\{0 \conn x\}}].
    \end{split}
\end{equation}
For any integer $r\ge1$, by the BK-inequality,
\begin{equation}\label{e:intdistdecomp}
    \begin{split}
        \E_p \!\left[V(B_r) \indi_{\{0\conn x\}}\right]
            =& \sum_{z,z'} \P_p\big(0\stackrel{r}{\longleftrightarrow} z, \{z,z'\} \text{ open},  0\conn x\big)\\
            \le&  \sum_{y,z,z'} \P_p\big(\{0\stackrel{r}{\longleftrightarrow} y\} \circ \{y\stackrel{r}{\longleftrightarrow} z\} \circ \{\{z,z'\} \text{ open}\}\circ\{y\conn x\}\big)\\
            \le&  \sum_{y,z,z'} \P_p\big(0\stackrel{r}{\longleftrightarrow} y\big)\, \P_p\big(y\stackrel{r}{\longleftrightarrow} z\big)\, p D(z'-z)\, \P_p\big(y\conn x\big)\\
            =& p \sum_{y,z}\P_p \big(0 \stackrel{r}{\longleftrightarrow} y\big) \, \P_p \big(y \stackrel{r}{\longleftrightarrow} z \big) \,\P_p \big(y \conn x \big).
    \end{split}
\end{equation}
Hence,
\begin{equation}\label{e:expiicvolint}
        \Eiic[V(B_r)] \le \limp \frac{p}{\chi(p)} \sum_{x,y,z}\P_p \big(0 \stackrel{r}{\longleftrightarrow} y\big) \, \P_p \big(y \stackrel{r}{\longleftrightarrow} z \big) \,\P_p \big(y \conn x \big) =  p_c\; \Epc[\lvert B_r \rvert]^2.
\end{equation}
Finally, by Theorem \ref{th:KozNac} and Markov's inequality,
\begin{equation}
    \Piic(V(B_r) \ge \lambda r^2) \le \frac{\Eiic[V(B_r)]}{\lambda r^2} \le \frac{C}{\lambda}.
\end{equation}

Now we derive the bound for the other end of the interval, i.e., the bound on $\Piic(V(B_r) \le \lambda^{-1} r^2)$.
Since $V(B_r) \ge \lvert B_r \rvert-1$ for any configuration, it is sufficient to prove the statement for the vertex volume $|B_r|$ instead of the edge volume $V(B_r)$.

Observe that since $B_r \subseteq \zri$, by Lemma \ref{lem:BPLR},
\begin{equation}
	\begin{split}
       		\Piic\left(|B_r| \le\lambda^{-1} r^2\right) &=  \Piic\left(\big\{|B_r| \le\lambda^{-1} r^2\big\}\text{ on } \zri\right)\\
			& = \limp \frac1{\chi(p)}\sum_x \P_p\left(\big\{|B_{r}|\le \lambda^{-1}\,r^2\big\}\text{ on }\zrx , 0\conn x\right)\\
			&\le \limp \frac1{\chi(p)}\sum_x \P_p\big(0 \conn x,\, |B_{r}|\le \lambda^{-1}\,r^2\big).
	\end{split}
\end{equation}

Define $\Br = B_{j}$ with $j$ being the smallest integer in $[r/2,r]$ satisfying $|\partial B_{j}|\le2\lambda^{-1}r$. Such a $j$ always exists when $|B_{r}|\le\lambda^{-1}\,r^2$.
Then
\begin{equation}\label{e:proof3}
    \big\{0\conn x,\,|B_{r}|\le \lambda^{-1}\,r^2\big\}
    \subseteq
    \bigcup^{\centerdot}_{A \text{ adm.}} \{0 \conn x,\Br=A\}
\end{equation}
where the disjoint union over ``$A$ adm.'' is over all sets $A\subset\Zd$ that are \emph{admissible}. Here, \emph{admissible} means that
$\P_p(\Br=A)>0$ and $|\partial A|\le 2\lambda^{-1}\,r$.
It follows that, for $x\nin A$,
\begin{eqnarray}
    \P_p(|B_{r}|\le \lambda^{-1}\,r^2 , 0\conn x)
    &\le& \sum_{A \text{ adm.}} \P_p(0\conn x\mid \Br=A)\,\P_p(\Br=A)\nnb
    &\le& \sum_{A \text{ adm.}} \sum_{y\in\partial A} \P_p(y\conn x)\,\P_p(\Br=A)
\end{eqnarray}
as $\Pp(y \conn x $ off $A \vert B^{(r)}=A) \le \Pp(y \conn x)$.
For $x\in A$,
\begin{equation}
    \P_p(|B_{r}|\le \lambda^{-1}\,r^2 , 0\conn x)
    =\sum_{A \text{ adm.}} \P_p(\Br=A).
\end{equation}

Using translation invariance,
\begin{equation}
    \begin{split}
    \chi(p)^{-1}\sum_x \P_p\big(|B_{r}|\le \lambda^{-1}\,r^2 , 0\conn x\big) \le& \sum_{A \text{ adm.}} \P_p(\Br=A)\sum_{y\in\partial A} \chi(p)^{-1}\sum_{x\nin A} \P_p(y\conn x)\\
            \le& \sum_{A \text{ adm.}} \P_p(\Br=A)\cdot 2\lambda^{-1}\,r\cdot 1,
    \end{split}
\end{equation}
since $|\partial A|\le2\lambda^{-1}\,r$.
Finally,
\begin{equation}
    \bigcup^{\centerdot}_{A \text{ adm.}} \{\Br=A\}\subseteq \big\{\partial B_{r/2}\neq\emptyset\big\},
\end{equation}
and by Theorem \ref{th:KozNac} the probability of the event on the right-hand side is bounded above by $C/ r$. Thus,
\begin{equation}
	\Piic(|B_r| \le \lambda^{-1} r^2) \le 2 \lambda^{-1}\, r \Ppc\left(\partial B_{r/2} \neq \emptyset \right) \le C/\lambda,
\end{equation}
completing the proof. \qed
\medskip

We now prove the bound on the effective resistance $\Reff(0,B_r^c)$ in \refeq{Bound4}.
We need the following two lemmas:
\bl[A lower bound on the effective resistance for critical percolation]
    Consider a percolation model such that the triangle condition holds.
    For any $p\le p_c$, $\lambda>1$, $r\ge1$,
    \begin{equation}\label{e:proof1}
        \P_p(\Reff(0, \partial B_r)\le \lambda^{-1} r)\le C / (\lambda r).
    \end{equation}
\el
\proof
The statement is proved along the lines of \cite[Proof of Lemma 2.6]{KozNac09}. In particular, see the last displayed inequality in that proof. \qed
\bl[An upper bound on $\Piic$]
    If $E$ is an event measurable with respect to $B_r$, $r\ge1$, then
    \begin{equation}\label{e:proof2}
        \Piic(E)\le C\,\sqrt{r\,\P_{p_c}(E)}.
    \end{equation}
\end{lemma}
\proof
Since $E$ is measurable with respect to $B_r$ and therefore measurable with respect to $\zri$, by Lemma \ref{lem:BPLR},
\begin{equation}\label{e:proof00}
	\Piic(E) = \Piic\left( E \text{ on } \zri\right) = \limp  \frac1{\chi(p)}\sum_x\P_p\left(\{E\text{ on } \zrx\}\cap\{0\conn x\}\right).
\end{equation}
Fix $M>0$ and $r\ge1$, and let $p<p_c$. (The constant $M$ will be optimized below.) We can bound
\begin{multline}\label{e:proof-split}
    P_p\left(\{E\text{ on }\zrx\}\cap\{0\conn x\}\right) \le \P_p\left(E\cap\{0\conn x\}\right)\\
    \le \P_p \left(\{|\dBr|>M\}\cap\{0\conn x\}\right)
    + \P_p \left(E\cap\{|\dBr|\le M,0\conn x\}\right).
\end{multline}
For the first term on the right hand side we use the BK-inequality to estimate
\begin{equation}
    \E_p\big[|\dBr|\,\indi_{\{0\conn x\}}\big]
    \le \sum_y \P_p\big(\{0\stackrel{r}{\longleftrightarrow} y\}\circ\{y\conn x\}\big)
    \le \sum_y \P_p\big(0\stackrel{r}{\longleftrightarrow} y\big)\,\P_p\big(y\conn x\big).
\end{equation}
Hence Markov's inequality implies
\begin{equation}\label{e:proof8}
    \sum_x\P_p\big(|\dBr|>M,{0\conn x}\big)
    \le \frac{1}{M}\sum_y \P_p\big(0\stackrel{r}{\longleftrightarrow} y\big)\,\sum_x\P_p\big(y\conn x\big)
    \le \frac{1}{M}\,\E_{p_c}[|B_{r}|]\,\chi(p),
\end{equation}
and this is bounded above by $C r\chi(p)/M$ by \refeq{Bound1}.

For the last term in \refeq{proof-split} we proceed like \refeq{proof3} by writing
\begin{equation}\label{e:proof4}
    \P_p \left(E\cap\{|\dBr|\le M,0\conn x\}\right)
    =\sum_{A \text{ adm.}} \P_p(B_{r}=A)\,\P_p(0\conn x\mid B_{r}=A)
\end{equation}
where the sum over ``$A$ adm.'' now is the sum over all pairs of sets $A, \partial A$ satisfying $\{B_{r}=A\}\subset E$, $|\partial A|\le M$, and $\P_p(\dBr=\partial A)>0$.
For each such admissible $A$ (in particular using $|\partial A|\le M$),
\begin{equation}
    \sum_x\P_p(0\conn x\mid\dBr=A)
    \le \sum_x\sum_{y\in\partial A}\P_p(y\conn x)
    \le M\chi(p),
\end{equation}
where we used translation invariance to get the last inequality.
Since $E$ is measurable with respect to $B_{r}$,
\begin{equation}\label{e:proof5}
    \sum_{A \text{ adm.}} \P_p(B_{r}=A)\le \P_p(E).
\end{equation}
Combining \refeq{proof4}--\refeq{proof5} yields
\begin{equation}\label{e:proof6}
    \sum_x\P_p \left(E\cap\{|\dBr|\le M,0\conn x\}\right)
    \le \P_p(E)\,M\,\chi(p).
\end{equation}

Now, using \refeq{proof00}  and \refeq{proof-split} together with \refeq{proof8} and \refeq{proof6}, we get
\begin{equation}\label{e:proof7}
   \Piic(E)= \limp \frac1{\chi(p)}\sum_x\P_p\left(\{E\text{ on }\zrx\}\cap\{0\conn x\}\right)
    \le \frac{C r}{M}+\Ppc (E)\,M.
\end{equation}
Letting $M=\sqrt{r/\Ppc(E)}$ proves the lemma.
\qed

\proof[Proof of \refeq{Bound4}]
The statement \refeq{Bound4} follows from \refeq{proof1} and \refeq{proof2}.
\qed
\medskip

The final bound we need to establish for the proof of Theorem \ref{th:Goodballs}(b) is
\bl[An expectation bound]\label{lem:Boringbound}
    If the strong triangle condition is satisfied for some sufficiently small $\beta$, then
    \begin{equation}\label{e:expreffvolint}
        \Eiic[\Reff(0, \partial B_r) V(B_r)] \le C r^3.
    \end{equation}
\el
\proof
Note that $\Reff(0, \partial B_r) \le r$ since the intrinsic distance metric dominates the effective resistance metric. So \refeq{expreffvolint} follows immediately from \refeq{expiicvolint}.
\qed
\proof[Proof of Theorem \ref{th:Goodballs}(b)]
Combining Lemmas \ref{lem:Cballbounds} and \ref{lem:Boringbound} completes the proof. \qed

\section{Random walk on the backbone: proof of Theorem \ref{th:Goodballs}(c) and (d)}\label{sec:Ext-bb}
\subsubsection*{The extrinsic distance metric}
Using the bounds that we have established in the previous two sections and of Theorem \ref{ExpectationBounds}, we can easily establish most of the bounds that we need to prove Theorem \ref{th:Goodballs}(c). Recall Definitions \ref{def:spec}(i) and (ii).
\bl[Bounds on the extrinsic volume and effective resistance of the backbone]\label{lem:bbballbounds}
 If the strong triangle condition is satisfied for some sufficiently small $\beta$ and if either the model is finite-range, or Assumption \ref{ass:scaling} holds, then there exists a $C>0$ and a $\xi > 0$ such that for all $\lambda>0$ and all $r \ge r^\star = r^\star(\lambda)$,
    \begin{eqnarray}
        \label{e:bbvolprob} \Piic(\lambda^{-1} r^{\twa} \le V(\bbr) \le \lambda r^{\twa}) & \ge & 1- C /\lambda^{\xi};\\
        \label{e:bbreffprob} \Piic(\lambda^{-1} r^{\twa} \le \Reffbb(0, Q_r^c) \le \lambda r^{\twa}) &\ge &1- C / \lambda^{\xi}.
    \end{eqnarray}
\el
\proof We start with \refeq{bbvolprob}. Note that $V(\bbr) \le  V(\bbb \cap Q_r)$. Furthermore, from the definition of $\Nbb(r)$ it follows that $V(\bbb \cap Q_r) \le 2 \Nbb(r)$. So by Markov's inequality, uniformly in $r \ge 1$,
\begin{equation}
    \Piic(V(\bbr) > \lambda r^{\twa}) \le \frac{2 \Eiic[\Nbb(r)]}{\lambda r^{\twa}} \le \frac{C r^{\twa}}{\lambda r^{\twa}} = \frac{C}{\lambda},
\end{equation}
where the final inequality follows from Theorem \ref{ExpectationBounds}.

For the other bound in \refeq{bbvolprob}, we note that the number of edges in $\bbr$ exceeds the number of pivotal edges for the connection between $0$ and $Q_r^c$, so that $V(\bbr) \ge H(r)$ follows by the definition of $V(\bbr)$ and $H(r)$. Hence,
\begin{equation}
    \Piic( V(\bbr) < \lambda^{-1} r^{\twa}) \le \Piic(H(r) < \lambda^{-1} r^{\twa}) \le C/ \lambda^{\xi},
\end{equation}
where the final inequality follows from Lemma \ref{prop:piv-exit-time}.

To bound \refeq{bbreffprob} we first note that $\Reffbb(0,Q_r^c) \le \Nbb(r)$ by the series law of resistances. Furthermore, also by the series law of resistances, $\Reffbb(0,Q_r^c) \ge H(r)$, so we can apply the same arguments as we did for the proof of \refeq{bbvolprob} to establish \refeq{bbreffprob}. This completes the proof. \qed

The next lemma is similar to Lemma \ref{lem:Goodpoints1}, though much easier to prove:
\bl[All vertices in $\bbr$ are well behaved]\label{lem:Goodpointsbb}
    If the strong triangle condition is satisfied for some sufficiently small $\beta$, then there exists a $C>0$ such that for all $\lambda>0$,
    \begin{equation}\label{e:lemstat1}
        \Piic\left(\exists x  \in \bbr \text{\emph{ such that }}\Reff(0,x) \ge \lambda r^{\twa}\right) \le C/\lambda.
    \end{equation}
\el
\proof Having $\Reff(0,x) \ge \lambda r^{\twa}$ for some $x \in \bbr$ implies that $d_{\omega}(0,x) \ge \lambda r^{\twa}$ also. This means that $\bbr$ contains at least $\lambda r^{\twa}$ vertices, i.e.,
\begin{equation}
    \Piic\left(\exists x \in \bbr \text{ such that }\Reff(0,x) \ge \lambda r^{\twa}\right) \le \Piic (\Sbbr \ge \lambda r^{\twa}) \le C /\lambda,
\end{equation}
where the final inequality is due to \refeq{markov}. \qed

\medskip
The last lemma we need concerns an upper bound on the expectation of $\Reffbb(0,Q_r^c) V(\bbr)$.
\bl[Upper bound on the expectation of $\Reffbb(0,Q_r^c) V(\bbr)$]
    Let $r \ge 1$. If the strong triangle condition is satisfied for some sufficiently small $\beta$, then there exists a constant $C > 0$ such that
    \begin{equation}\label{e:bbprodbd}
        \Eiic[\Reffbb(0,Q_r^c) V(\bbr)] \le C r^{2 \twa}.
    \end{equation}
\el
\proof Since $V(\bbr) \le 2 \Nbb(r)$ and also $\Reffbb(0,Q_r^c) \le \Nbb(r)$ we can bound
\begin{equation}\label{e:bbprodbd1}
    \Eiic[\Reffbb(0,Q_r^c) V(\bbr)] \le 2 \Eiic[\Nbb(r)^2].
\end{equation}
We can express $\Nbb(r)$ as a sum of indicator functions:
\begin{equation}\label{e:Nbbindi}
    \Nbb(r) = \sum_{b \in \Be} \indi_{\{0 \conn \ulb\}\circ \{b \text{ open}\} \circ \{\olb \conn \infty\}}.
\end{equation}

After substituting \refeq{Nbbindi} into \refeq{bbprodbd1}, a corollary to Lemma \ref{lem:BPLR} shows that we may then reverse the limit \refeq{IICdef} (i.e., \cite[Corollary 4.3]{HeyHofHul12a}),
\begin{multline}\label{e:bbprodbd2}
    2 \sum_{b, e \in \Be} \Eiic[\indi_{\{0 \conn \ulb\}\circ \{b \text{ open}\} \circ \{\olb \conn \infty\}}\indi_{\{0 \conn \ule\}\circ \{e \text{ open}\} \circ \{\ole \conn \infty\}}]\\
        =\limp \frac{2}{\chi(p)} \sum_{b, e \in \Be} \sum_{x \in \Zd} \Pp\Bigl(\big\{\{0 \conn \ulb\}\circ \{b \text{ open}\} \circ \{\olb \conn x\} \big\}\\
        \cap \big\{\{0 \conn \ule\}\circ \{e \text{ open}\} \circ \{\ole \conn x\}\big\}\Bigr).
\end{multline}

We want to write the event inside $\Pp$ in terms of disjointly occurring events, so that we may apply the BK-inequality. To this end, we define the following four events that together cover the event inside $\Pp$ (leaving dependence on $b, e$ and $x$ implicit on the left-hand sides):
\begin{eqnarray}
    E_1 &\equiv& \{e = b\} \cap \big(\{0 \conn \ule\} \circ \{e \text{ open}\} \circ \{\ole \conn x\}\big);\\
    E_2 &\equiv& \{e \neq b\} \cap \big(\{0 \conn \ule\}\circ \{e \text{ open}\} \circ \{\ole \conn \ulb\} \circ \{b \text{ open}\}\circ \{\olb \conn x\}\big);\\
    E_3 &\equiv& \{e \neq b\} \cap \big(\{0 \conn \ulb\}\circ \{b \text{ open}\} \circ \{\olb \conn \ule\} \circ \{e \text{ open}\}\circ \{\ole \conn x\}\big);\\
    \nonumber E_4 &\equiv& \{e \neq b\} \cap \Big( \bigcup_{z,w \in \Zd} \{0 \conn z\}\circ\{z \conn \ule\}\circ\{e \text{ open}\} \\
    && \qquad\circ\{\ole \conn w\}\circ \{z \conn \olb\} \circ \{b \text{ open}\} \circ \{\olb \conn w \} \circ \{w \conn x\}\Big).
\end{eqnarray}
It follows that
\begin{equation}\label{e:5eventsplit}
    \left(\{0 \conn \ule\}\circ\{e \text{ open}\}\circ\{\ole \conn x\}\right) \cap \big(\{0 \conn \ulb\}\circ\{b \text{ open}\}\circ\{\olb \conn x\}\big) \subseteq E_1 \cup E_2 \cup E_3 \cup E_4.
\end{equation}

Substituting the right-hand side of \refeq{5eventsplit} into \refeq{bbprodbd2}, we get the upper bound
\begin{equation}\label{e:ABCDE}
    \limp \frac{2}{\chi(p)} \sum_{x \in \Zd} \sum_{e,b \in \Be} \left[\Pp(E_1) + \Pp(E_2) + \Pp(E_3) + \Pp(E_4)\right].
\end{equation}

We bound the four sums separately. For each bound we start by applying the BK-inequality, sum over $x$ and take the limit to get a factor $\chi(p)$ and then take the limit $p \ua p_c$. For the sums involving $E_2$, $E_3$ and $E_4$ we drop the requirement that $\{e \neq b\}$ for an upper bound. We start with $\Pp(E_1)$:
\begin{equation}
     \limp \frac{2}{\chi(p)} \sum_{x \in \Zd} \sum_{e,b \in \Be} \Pp(E_1) \le 2 p_c \sum_{e \in \Be} \taupc (\ule) D(e) \le C r^{\twa},
\end{equation}
where the bound follows from the fact that $\sum_x D(x) =1$, Theorem \ref{ExpectationBounds} and \refeq{expident}.

To bound $\Pp(E_2)$ we rewrite the it as a convolution,
\begin{equation}
    \limp \frac{2}{\chi(p)} \sum_{x \in \Zd} \sum_{e,b \in \Be} \Pp(E_2) \le 2 p_c^2 (\indi_{Q_r} * \taupc * D * \taupc * D)(0) \le C r^{2 \twa}.
\end{equation}
The second bound follows after applying methods similar to those used in the proof of Theorem \ref{ExpectationBounds} as given in \cite{HeyHofHul12a}.
Interchanging the labels $e$ and $b$ shows that the same bound holds for $\Pp(E_3)$.

Finally, $\Pp(E_4)$ is upper bounded by
\begin{multline}
        \limp \frac{2}{\chi(p)} \sum_{x \in \Zd} \sum_{e,b \in \Be} \Pp(E_4) \\
        \le 2 p_c^2 \sum_{z,w \in \Zd} \sum_{e,b \in \Be} \taupc(z) \taupc(\ule -z) D(e) \taupc(w-\ole) \taupc(\ulb - z) D(b) \taupc(w - \olb).
\end{multline}
We extend the summation over $\ulb$ to $\Zd$, shift each term in the summation by $-\ule$ and relabel to get
\begin{equation}\label{e:pathplustriangle}
    2 p_c^2 \sum_{\ule' \in Q_r} \sum_{\ole',z',w',\ulb',\olb' \in \Zd}  \taupc(z') \taupc(\ule' -z') D(\ole') \taupc(w'-\ole') \taupc(\ulb'-w') D(b') \taupc(z' - \olb').
\end{equation}
From \cite[Theorem 1.3]{BorChaHofSlaSpe05b} we have the following bound:
\begin{equation}
    \sup_{z' \in \Zd}\sum_{\ole',w',\ulb',\olb'\in \Zd}D(\ole')\taupc(w'-\ole')\taupc(\ulb'-w')D(b')\taupc(z'-\olb') \le C \beta,
\end{equation}
(where $\beta$ is the same constant as given below \refeq{trianglecond}).
We can apply this bound to \refeq{pathplustriangle} to get the upper bound
\begin{equation}
    C \beta \sum_{\ule' \in Q_r} \sum_{z' \in \Zd} \taupc(z') \taupc(\ule'-z') = C \beta (\indi_{Q_r} * \taupc * \taupc)(0) \le C \beta r^{2 \twa},
\end{equation}
where the final bound again follows from applying methods similar to those used in the proof of Theorem \ref{ExpectationBounds} (cf. \cite{HeyHofHul12a}).

Adding the bounds for the sums over $\Pp(E_1)$, $\Pp(E_2)$, $\Pp(E_3)$ and $\Pp(E_4)$ establishes that $C r^{2 \twa}$ is an upper bound on \refeq{ABCDE}, and so we get the desired upper bound in \refeq{bbprodbd}, which in turn completes the proof. \qed

\subsubsection*{The intrinsic distance metric}
\bl[Bounds on the intrinsic volume and effective resistance of the backbone]
    If the strong triangle condition is satisfied for some sufficiently small $\beta$, then there exist $C, C'>0$ such that for all $\lambda>0$,
    \begin{eqnarray}
        \label{e:bbvolprobint} \Piic(r \le V(B_r^{\sss \bbb}) \le \lambda r) & \ge & 1- C /\lambda;\\
        \label{e:bbreffprobint} \Piic(\lambda^{-1} r \le \Reffbb(0, \partial B_r^{\sss \bbb})) &\ge &1- C / \sqrt{\lambda};\\
        \label{e:bbreffvolint}  \Eiic[V(B_r^{\sss \bbb})\Reffbb(0, \partial B_r^{\sss \bbb})] &\le& C' r^2.
    \end{eqnarray}
\el
\proof
We start by observing that the lower bound on $V(B_r^{\sss \bbb})$ holds trivially since it takes at least $r$ edges to reach distance $r$ in the intrinsic distance metric and $\partial B_r^{\sss \bbb} \neq \emptyset$ $\Piic$-a.s.

For the upper bound, start by applying Markov's inequality,
\begin{equation}\label{e:Vb1stmom}
    \Piic(V(B_r^{\sss \bbb}) \ge \lambda r) \le \frac{\Eiic[V(B_r^{\sss \bbb})]}{\lambda r}.
\end{equation}
The random variable $V(B_r^{\sss \bbb})$ is measurable with respect to $\zri$, so, by Lemma \ref{lem:BPLR}, we may reverse the IIC limiting scheme and apply arguments similar to \refeq{intdistdecomp} to yield
\begin{equation}\label{e:expbbvolint}
    \begin{split}
        \Eiic[V(B_r^{\sss \bbb})]    &=\limp \frac{1}{\chi(p)} \sum_{x} \E_p [V(B_r^{\sss \bbb})\indi_{\{0 \conn x\}}]\\
                                &\le \frac{1}{\chi(p)} \sum_{x,z,z'} \Pp(0 \stackrel{r}{\longleftrightarrow} z)p D(z'-z) \Pp(z' \conn x) \le p_c \, \Epc[\lvert B_r \rvert] \le C r,
    \end{split}
\end{equation}
where the final inequality follows from Theorem \ref{th:KozNac}. Substitution in \refeq{Vb1stmom} yields the required upper bound in \refeq{bbvolprobint}.

To prove \refeq{bbreffprobint} it suffices to observe that by the cutting law for resistances, $\Reffbb(0, \partial B_r^{\sss \bbb})\ge \Reff(0, \partial B_r)$, so that the required bound follows by Lemma \ref{lem:Cballbounds}.

Finally, since $\Reffbb(0, \partial B_r^{\sss \bbb}) \le r$ holds trivially, \refeq{bbreffvolint} follows from \refeq{expbbvolint}. \qed

\section{The modified exit time: proof of Theorem \ref{th:Goodballs}(e)}\label{sec:Ext-mod}
\bl[Bounds on the modified volume and effective resistance]\label{lem:modballs}
  If the strong triangle condition is satisfied for some sufficiently small $\beta$ and if Assumption \ref{ass:scaling} holds, then there exist $C, C', C''>0$ and a $\xi > 0$ such that for all $\lambda>0$ and $r \ge r^\star (\lambda)$,
    \begin{eqnarray}
        \label{e:modvolprob} \Piic(\lambda^{-1} r^{2\twa} \le V(\bqr) \le \lambda r^{2\twa}) & \ge & 1- C /\lambda^{\xi};\\
        \label{e:modreffprob} \Piic(\lambda^{-1} r^{\twa} \le \Reffmod(0, Q_r^c) \le \lambda r^{\twa}) &\ge &1- C / \lambda^{\xi};\\
        \label{e:modreffunif} \Piic\left(\exists x  \in \bqr \text{\emph{ such that }}\Reff(0,x) \ge \lambda r^{\twa}\right) &\le& C'/\lambda;\\
        \label{e:modvolreff} \Eiic[\Reffmod(0,Q_r^c) V(\bqr)] &\le& C'' r^{3 \twa}.
    \end{eqnarray}
\el
\proof By Theorem \ref{ExpectationBounds}, the upper bounds in \refeq{modvolprob}--\refeq{modvolreff} follow by the same proof as for Lemmas \ref{lem:Reffub}, \ref{lem:UpperVolBoundBall}, \ref{lem:boundReffV}, and \ref{lem:Goodpoints1}, respectively, when we replace the bound $C r^2$ by $C r^{\twa}$ wherever this bound is used (i.e., in \refeq{rhs1}, \refeq{rhs2}, \refeq{mod1}, \refeq{mod2}, \refeq{markov} and \refeq{BKRpluslim}).

The lower bounds in \refeq{modvolprob} and \refeq{modreffprob} follow by almost the same proof as Proposition \ref{prop-piv-eff-vol}. The main change due to the modification is that now $\Reffmod(0,Q_r^c) \ge H(r)$, instead of $\Npiv{0,Q_r^c}$. The result is that we do not need a bound on $\Nbad(0,Q_r^c)$ as in Proposition \ref{prop-piv-bd}, so we do not need Assumption \ref{ass:onearm}, and hence the proof works for long-range spread-out percolation as well. \qed

Lemma \ref{lem:modballs} implies Theorem \ref{th:Goodballs}(e). Theorem \ref{th:IIC-mod} then follows by the same proof as Theorem \ref{th:Eucexittimes} (with the appropriate changes made to the exponents).

\section{Exit time for the long-range IIC: proof of Theorem \ref{th:LRP-IIC}}\label{sec:LRP-IIC}

In this section we only consider distributions $D$ of the form \refeq{DLRP} and we write $\varrho =(4 \wedge \alpha)/2$.
Let $a \in (0,1/2)$ and $0<b< (1+a)/3$ and consider these fixed for the rest of this section. Write $n= \lfloor \lambda^{-(1+a)/3} r^{\varrho} \rfloor$. Throughout this section we will write $e_n$ for the $n$th backbone pivotal edge. It should be viewed as a $\Zd \times \Zd$ valued random variable under the measure $\Piic$. We will write $b_n$ when we mean that the edge is fixed.

Define
\begin{equation}\label{e:Jtilde}
    \tilde J(\lambda) = \left\{r \in [1,\infty]\,:\,  \Reff (0, Q_r^c) \le \lambda^b r^\varrho,\,\Reffpn \le \lambda n,\,V(\tCcal^{e_{n}}(0)) \le \lambda n^2 \right\}.
\end{equation}

\bl[Conditional bound on the exit time]\label{lem:condbd}
   There exists $\kappa>0$ such that for any $\lambda$, $r$, and $\omega$ such that $r \in \tilde J (\lambda)$,
    \begin{equation}
       P^0_{\omega} (\tqr > \lambda r^{3 \varrho}) \le 2/\lambda^\kappa.
    \end{equation}
\el
\proof Let $\sigma_{\ole_n}$ be the hitting time of $\ole_n$. Then we can bound
\begin{equation}\label{e:splitP0}
    P^0_{\omega} (\tqr > \lambda r^{3 \varrho}) \le P^0_{\omega} (\sigma_{\ole_n} > \lambda r^{3 \varrho}) + P^0_{\omega} (\sigma_{\ole_n} < \tqr).
\end{equation}
To bound the second term, we use the following standard bound (cf. \cite[(4)]{BerGanPer03})
\begin{equation}
    P^0_{\omega} (\sigma_{\ole_n} < \tqr)  \le \frac{\Reff(0, Q_r^c)}{\Reff(0, \ole_n)}.
\end{equation}
For $\Reff(0,  \ole_n)$ we have that $n$, the number of pivotals for $0 \conn \ole_n$, is a lower bound by the series law of resistances, so by the definition of $\tilde J (\lambda)$,
\begin{equation}
    P^0_{\omega} (\sigma_{\ole_n} < \tau_{Q_r}) \le \frac{\lambda^b r^\varrho}{n} \le \frac{1}{\lambda^{(1+a)/3-b}}.
\end{equation}

For the first term in \refeq{splitP0} we use Markov's inequality:
\begin{equation}
    P^0_{\omega} (\sigma_{\ole_n} > \lambda r^{3 \varrho}) \le \frac{E_{\omega}^0 \sigma_{\ole_n}}{\lambda r^{3 \varrho}}.
\end{equation}
By the Green's function interpretation of the hitting time (see Defintion \ref{def:properties}(v)):
\begin{equation}
	E_\omega^0 \sigma_{\ole_n} = \sum_{y \in \tCcal^{e_n} (0)} G_{\tCcal^{e_n} (0)} (0, y) \mu_y
			 \le \sum_{y \in \tCcal^{e_n} (0)} G_{\tCcal^{e_n} (0)} (0,0) \mu_y
			 = \Reff(0, \ole_n) (V(\tCcal^{e_n} (0))+1)
\end{equation}
where the last inequality again follows from our choice of $\tilde J (\lambda)$.
So it follows that
\begin{equation}
    P^0_{\omega} (\sigma_{\ole_n} > \gamma r^{3 \varrho}) \le \frac{\lambda^2 n^3}{\lambda r^{3 \varrho}} \le 1 /\lambda^a.
\end{equation}
Let $\kappa = \min \{a, \frac{1+a}{3} -b\} \in (0,1/2)$, then
\begin{equation}
    P^0_{\omega} (\tqr > \lambda r^{3 \varrho}) \le 2/\lambda^\kappa.
\end{equation}
\qed

The other ingredient needed for the proof of Theorem \ref{th:LRP-IIC} is the following proposition:
\begin{prop}\label{prop:GoodballLRP}
    For any sufficiently large $\lambda$ and for all $r \ge r^*(\lambda)$, let $\nu = \min\{b, 1/2\}$. There exists a constant $c>0$ such that
    \begin{equation}
        \Piic(r \in \tilde J(\lambda)) \ge 1-c/\lambda^\nu.
    \end{equation}
\end{prop}
The proof of this proposition is given in the next subsection.

\proof[Proof of Theorem \ref{th:LRP-IIC} subject to Proposition \ref{prop:GoodballLRP}]
By Lemma \ref{lem:condbd} and Proposition \ref{prop:GoodballLRP} and the definition of $\mathsf{P}^\star$, \refeq{Pstar}, we have (with $\vep =\min \{\kappa,\nu\}$ and $\kappa, \nu$ as defined above)
\begin{equation}
    \mathsf{P}^\star (\tau_{Q_r} \ge \lambda r^{3 \varrho}) \le  \int\limits_{\{r \in \tilde J (\lambda)\}} P^0_{\omega} (\tau_{Q_r} \ge \lambda r^{3 \varrho}) \Piic(\d \omega) + \Piic(r \nin \tilde J (\lambda)) \le C/\lambda^{\vep}.
\end{equation}
\qed

\subsection{The proof of Proposition \ref{prop:GoodballLRP}}
The proof of Proposition \ref{prop:GoodballLRP} is given in the three lemmas below; one lemma for each of the three restrictions in \refeq{Jtilde}.

\bl[An upper bound on the effective resistance for LRP]\label{lem:reffLRP}
    There exists a constant $c>0$ such that for any sufficiently large $\lambda$, $b>0$ and for all $r \ge r^\star(\lambda)$,
    \begin{equation}
        \Piic \left(\Reff (0, Q_r^c) \le \lambda^b r^{\varrho}\right) \ge 1 - c/\lambda^b.
    \end{equation}
\el
\proof
By Lemma \ref{lem:Reffub} the statement holds when $\alpha \ge 4$, so that we only need to prove it for the case $\alpha <4$, that is, when $\varrho=\alpha/2$.

Write $m = \lceil \lambda^b r^{\alpha/2} \rceil$. We start by noting that $\{\Reff (0, Q_r^c) \le \lambda^b r^{\varrho}\} \supseteq \{\lvert B_m \cap Q_r^c\rvert \neq 0\}$, so that
\begin{equation}
    \Piic(\Reff(0,Q_r^c) \le m) \ge 1- \Piic(\lvert B_m \cap Q_r^c \rvert =0).
\end{equation}
Define $\Ccal^r (x)$, the \emph{$r$-truncated cluster} of $x$, as the modified configuration of $\Ccal(x)$ where all edges of length at least $2r$ have been closed, and define the \emph{$r$-truncated intrinsic ball} of radius $m$ as
\begin{equation}
    B_m^{(r)} = \{x : 0 \mconn x \text{ on } \Ccal^r (0) \}.
\end{equation}

Then,
\begin{equation}\label{e:BmQrcsplit}
    \Piic(\lvert B_m \cap Q_r^c \rvert = 0) \le \Piic(\lvert B_m^{(r)} \rvert < m^2/\lambda^b) + \Piic(\lvert B_m^{(r)} \rvert \ge m^2 /\lambda^b, \lvert B_m \cap Q_r^c \rvert=0).
\end{equation}
We will bound both terms separately.

For the first term we observe that
\begin{equation}\label{e:Smrsplit}
    \Piic(\lvert B_m^{(r)} \rvert \ge m^2 / \lambda^b) \ge \Piic(\lvert B_m\rvert \ge 2 m^2 / \lambda^b) - \Piic(\lvert B_m \setminus B_m^{(r)} \rvert \ge m^2 / \lambda^b).
\end{equation}
By Lemma \ref{lem:Cballbounds} the first term is bounded from below by $1-c_1 /\lambda^b$. An upper bound on the second term follows by an application of Markov's inequality:
\begin{equation}\label{e:BmSmrfirstmom}
    \Piic(\lvert B_m \setminus B_m^{(r)} \rvert \ge m^2 /\lambda) \le \frac{\Eiic[\lvert B_m \setminus B_m^{(r)} \rvert]}{m^2 / \lambda^b}.
\end{equation}
Note that
\begin{equation}
    \Eiic[\lvert B_m \setminus B_m^{(r)} \rvert] = \sum_{x \in \Zd} \Piic(x \in B_m \setminus B_m^{(r)}),
\end{equation}
and furthermore, since $\{x \in B_m \setminus B_m^{(r)}\}$ is measurable with respect to $\mathsf{Z}_{m}^{\sss \infty} $, we can use Lemma \ref{lem:BPLR} to reverse the IIC-limit:
\begin{equation}\label{e:BmSmrlimrev}
    \begin{split}
        \sum_{x \in \Zd} \Piic\left(x \in B_m \setminus B_m^{(r)}\right) &= \limp \frac{1}{\chi(p)} \sum_{x,y \in \Zd} \Pp(x \in B_m \setminus B_m^{(r)}, 0 \conn y)\\
            & =  \limp \frac{1}{\chi(p)} \sum_{x,y \in \Zd}\bigl(\Pp(x \in B_m\setminus B_m^{(r)}, y \in B_{3m})\\
            &\qquad + \Pp(x \in B_m \setminus B_m^{(r)},y \notin B_{3m}, 0 \conn y)\bigr)\\
            & = \limp \frac{1}{\chi(p)} \sum_{x,y \in \Zd}\Pp(x \in B_m \setminus B_m^{(r)},y \notin B_{3m}, 0 \conn y).
    \end{split}
\end{equation}
The last equality follows since the sum over $y \in B_{3m}$ almost surely gives at most a finite contribution, whereas $\chi(p)$ diverges in the limit $p \ua p_c$.
It follows by the definition of $B_m^{(r)}$ that
\begin{multline}
        \{x \in B_m \setminus B_m^{(r)},y \notin B_{3m}, 0 \conn y\}\\
        \subseteq  \left(\bigcup_{\ulb,\olb,z \in \Zd} \{0 \mconn z\}\circ\{z \mconn \ulb\}\circ\{b \text{ open},\lvert b\rvert > 2r\} \circ\{\ole \mconn z\} \circ\{z \conn y\}\right)\\
         \cup \left(\bigcup_{\ulb,\olb,z \in \Zd} \{0 \mconn \ulb\}\circ\{b \text{ open}, \lvert b \rvert >2r\} \circ \{\olb \mconn z\} \circ \{z \mconn x\} \circ \{z \mconn y\} \right).
\end{multline}
Substituting the right-hand side in the last line of \refeq{BmSmrlimrev}, and applying the BK-inequality, we get the upper bound
\begin{multline}\label{e:subsplit}
        \limp \frac{1}{\chi(p)} \sum_{x,y,z,\ulb \in \Zd} \sum_{\olb \in \Zd : \lvert b \rvert > 2r} \Big(  \Pp(0 \mconn z) \Pp(z \mconn \ulb) p D(b) \Pp(\olb \mconn x) \Pp(z \conn y) \\
         + \Pp(0 \mconn \ulb) p D(b) \Pp(\olb \mconn z) \Pp(z \mconn x) \Pp(z \conn y) \Big).
\end{multline}
Taking the sum over $y$ gives a factor $\chi(p)$. After this we can take the limit $p \ua p_c$. Then, by translation invariance, the sum over $x$ gives a factor $\Epc[\lvert B_m \rvert]$. We get that \refeq{subsplit} is equal to
\begin{multline}\label{e:BmSmrbd1}
    \Epc[\lvert B_m \rvert] \sum_{z,\ulb \in \Zd} \sum_{v \in Q_{2r}^c}\bigl( \Ppc(0 \mconn z) \Ppc(z \mconn \ulb) p D(v)\\
     + \Ppc(0 \mconn \ulb) p D(v) \Ppc(\ulb+v \mconn z) \bigr).
\end{multline}
By the definition of $D(x)$ in \refeq{DLRP} there exist constants $\xi \ge \zeta >0$ such that for all $r \ge 1$,
\begin{equation}\label{e:Dbds}
    \frac{\zeta}{r^{\alpha}} \le \sum_{x \in Q_{2r}^c} D(x) \le \frac{\xi}{r^{\alpha}}.
\end{equation}
Summing over $z,v$ and $\ulb$ in \refeq{BmSmrbd1} and applying \refeq{Dbds}, we get the following upper bound:
\begin{equation}
    \Eiic[\lvert B_m \setminus B_m^{(r)} \rvert] \le \frac{2 \xi}{r^{\alpha}} \Epc[\lvert B_m \rvert]^3 \le \frac{C m^3}{r^{\alpha}},
\end{equation}
where the last inequality follows from Theorem \ref{th:KozNac}. Substituting this bound in \refeq{BmSmrfirstmom}, and using $m=\lceil \lambda^b r^{\alpha/2} \rceil$, we get
\begin{equation}
    \Piic(\lvert B_m \setminus B_m^{(r)} \rvert \ge m^2 /\lambda^b) \le \frac{C m^3 r^{-\alpha}}{m^2 /\lambda^b} = C_1 \lambda^{2b} r^{-\alpha/2}
\end{equation}
and then substituting this result into \refeq{Smrsplit}, we get,
\begin{equation}\label{e:Smrbd}
    \Piic(\lvert B_m^{(r)} \rvert < m^2 / \lambda^b) \le \frac{c_1}{\lambda^b} + C_1 \lambda^{2b} r^{-\alpha/2}.
\end{equation}

Now we prove that the second term in \refeq{BmQrcsplit} is also small. We start by splitting the probability once more:
\begin{multline}\label{e:BmSmrbd0}
    \Piic(\lvert B_m^{(r)} \rvert \ge m^2 /\lambda^b, \lvert B_m \cap Q_r^c \rvert=0) \\
    \le \Piic(\lvert B_{2m} \rvert > 4 \lambda^b m^2) + \Piic(\lvert B_m^{(r)} \rvert \ge m^2 /\lambda^b, \lvert B_m \cap Q_r^c \rvert =0, \lvert B_{2m} \rvert \le 4 \lambda^b m^2).
\end{multline}
By Lemma \ref{lem:Cballbounds}, the first term can be bounded from above by $c_2 /\lambda^b$. For the second term we observe that the event is measurable with respect to $\mathsf{Z}_{2m}^{\sss \infty} $, so we can reverse the IIC-limit to get
\begin{equation}\label{e:BmSmrbd2}
    \limp \frac{1}{\chi(p)} \sum_{x \in \Zd} \Pp(\lvert B_m^{(r)} \rvert \ge m^2 / \lambda^b, \lvert B_m \cap Q_r^c \rvert =0, \lvert B_{2m} \rvert \le 4 \lambda^b m^2, 0 \conn x).
\end{equation}
We will use the `admissibility method' of \cite{KozNac09} (also used in the proof of Lemma \ref{lem:Cballbounds} above). Whenever $\lvert B_{2m} \rvert \le 4 \lambda^b m^2$ occurs, there must exists at least one $j^* \in [m,2m]$ such that $\lvert \partial B_{j^*} \rvert \le 4 \lambda^b m$. Let $j$ be the first such $j^*$ and define $B^{(2m)}=B_j$. We call a set $A \subset \Zd$ `admissible' when $\Pp(B^{(2m)} = A)>0$ and $\lvert \partial A \rvert \le 4 \lambda^b m$. Then
\begin{equation}
    \{0 \conn x, \lvert B_{2m} \rvert \le 4 \lambda^b m^2 \} \subseteq \dot{\bigcup_{A \text{ adm.}}} \{0 \conn x, B^{(2m)} = A\}.
\end{equation}
For $x \in A$, we note that this contribution vanishes in the limit $p \ua p_c$. For $x \notin A$ the event $\{\partial A \conn x$ off $A\}$ is independent of the status of the edges in $A$, so
\begin{multline}
    \Pp(\lvert B_m^{(r)} \rvert \ge m^2/\lambda^b, \lvert B_m \cap Q_r^c \rvert = 0, \lvert B_{2m} \rvert \le 4 \lambda^b m^2, 0 \conn x) \\
    \le \sum_{A \text{ adm.}} \sum_{y \in \partial A} \Pp(y \conn x) \Pp(\lvert B_m^{(r)} \rvert \ge m^2/\lambda^b, \lvert B_m \cap Q_r^c \rvert = 0, B^{(2m)}=A).
\end{multline}
We substitute the right-hand side in \refeq{BmSmrbd2} and sum over $x$ and $y$, using that $\lvert \partial A \rvert \le 4 \lambda^b m$, to get
\begin{equation}\label{e:BmSmrbd3}
    \begin{split}
    \limp \frac{1}{\chi(p)} \sum_{x \in \Zd} \sum_{A \text{ adm.}} \sum_{y \in \partial A} & \Pp(y \conn x) \Pp(\lvert B_m^{(r)} \rvert \ge m^2/\lambda^b, \lvert B_m \cap Q_r^c \rvert = 0, B^{(2m)}=A)\\
        &\le 4\lambda^b m \sum_{A \text{ adm.}} \Ppc(\lvert B_m^{(r)} \rvert \ge m^2/\lambda^b, \lvert B_m \cap Q_r^c \rvert = 0, B^{(2m)}=A)\\
        &\le 4\lambda^b m \Ppc(\lvert B_m^{(r)} \rvert \ge m^2/\lambda^b, \lvert B_m \cap Q_r^c \rvert = 0)
    \end{split}
\end{equation}
where in the last step we used that all $A$ are mutually disjoint.

We say that an edge $e$ is \emph{long} when $\lvert e \rvert > 2r$. For the last steps of the proof we use that
\begin{equation}
    \{\lvert B_m^{(r)} \rvert \ge m^2 / \lambda^b, \lvert B_m \cap Q_r^c \rvert = 0 \} \subseteq \{\lvert B_m^{(r)} \rvert \ge m^2 /\lambda^b, \nexists \text{ long edge touching }B_m^{(r)}\}.
\end{equation}
The status of long edges is translation invariant and independent of the status of the edges in $B_m^{(r)}$, so that by the above and by \refeq{Dbds} we can bound the right-hand side of \refeq{BmSmrbd3} from above by
\begin{multline}
    4 \lambda^b m (1-\Ppc(\exists \text{ long edge touching }0))^{m^2/\lambda^b} \Ppc(\lvert B_m^{(r)} \rvert \ge m^2 /\lambda^b) \\
     \le 4 \lambda^b m \left(1-\frac{\zeta}{r^{\alpha}}\right)^{m^2/\lambda^b} \Ppc(\lvert B_m^{(r)} \rvert \ge m^2 /\lambda^b).
\end{multline}

For the final bound, we use that
\begin{equation}
    \Ppc(\lvert B_m^{(r)} \rvert \ge m^2 /\lambda^b) \le \Ppc(\lvert \Ccal^{r} (0) \rvert \ge m^2 /\lambda^b) \le \frac{C}{\sqrt{m^2 /\lambda^b}},
\end{equation}
where the final inequality follows from \cite{BarAiz91} and \cite{HeyHofSak08} (see also \refeq{Vvolsc} below).
Recall that $ m = \lceil \lambda^b r^{\alpha/2}\rceil $, so
\begin{equation}\label{e:BmSmrbd4}
    \Piic(\lvert B_m^{(r)} \rvert \ge m^2 /\lambda^b, \lvert B_m \cap Q_r^c \rvert=0) \le C \lambda^{3b/2} \left(1-\frac{\zeta}{r^{\alpha}}\right)^{m^2/\lambda^b} \le C_2 \lambda^{3b/2} e^{-\zeta \lambda^b}.
\end{equation}

Finally, combining the bounds \refeq{BmQrcsplit}, \refeq{Smrbd}, \refeq{BmSmrbd0} and \refeq{BmSmrbd4} we get for $r \ge r^\star (\lambda) = \lambda^{6b/\alpha}$,
\begin{equation}
    \Piic (\Reff (0, Q_r^c) \le \lambda^b r^{\alpha/2}) \ge 1 - \frac{(c_1 + c_2)}{\lambda^b} - \frac{C_1 \lambda^{2b}}{r^{\alpha/2}} -C_2  \lambda^{3b/2} e^{-\zeta \lambda^b} \ge 1-c/\lambda^b.
\end{equation}
\qed

The next lemma shows that it is unlikely that the effective resistance between $0$ and the $n$th backbone pivotal edge $e_n$ is large. Write $\Bb_n$ for the subgraph of $\Bb$ induced by the backbone up to $e_n$ . 
\begin{lemma}\label{lem:reffenlambda} If the strong triangle condition is satisfied for some sufficiently small $\beta$, then, for any sufficiently large $\lambda$ and for all $n \ge 1$,
\begin{equation}\label{e:reffenlambda1}
	\Eiic[|\Bb_n|] \le C n,
\end{equation}
and as a result,
\begin{equation}\label{e:reffenlambda2}
	\Piic \big(\Reff(0, \ole_n) \le \lambda n\big) \ge1- C /\lambda.
\end{equation}
\end{lemma}
\proof
We start by deriving  \eqref{e:reffenlambda1} from  \eqref{e:reffenlambda2}.
The effective resistance is dominated by the number of edges $\Bb_n$, and since each edge edge in $\Bb_n$ connects to precisely $2$ vertices in $\Bb_n$, 
\begin{equation}\label{e:reffenmarkov}
	\Piic(\Reff(0, \ole_n) \ge \lambda n) \le \Piic(|\Bb_n| \ge \lambda n / 2) \le \frac{2 \Eiic[|\Bb_n|]}{\lambda n}.
\end{equation}
It thus remains to show \eqref{e:reffenlambda2}. 
We can bound the graph distance from above by the number of edges in the backbone up to the $n$th pivotal, so by Lemma \ref{BPRL},
\begin{equation}
	\begin{split}
	\Eiic[ |\Bb_n|]  	
					&= \sum_{z} \Piic(\{0 \conn z\} \text{ on } \Bb_n)\\
					& = \sum_{z} \limp \frac{1}{\chi(p)} \sum_{x \in \Zd} \Pp (\{0 \conn z\} \on \Bb_n^{x}, 0 \conn x).
	\end{split}
\end{equation}
where $\Bb_n^x$ is the backbone graph for the connection from $0$ to $x$ up to the $n$th pivotal. 

If $z$ is a vertex in one of the $n$ first sausages between $0$ and $x$, then it follows that $z$ is connected to $0$ with fewer than $n$ pivotal edges, and there is a disjoint connection from $z$ to $x$, i.e.,
\begin{equation}
	\Ppc(\{0 \conn z\} \on \Bb_n^{x}, 0 \conn x) \le \Ppc (\{0 \conn z \text{ with }\le n \text{ pivotals}\}  \circ \{z \conn x\}).
\end{equation}
Applying the BK-inequality, summing over $x$ and taking the limit we get
\begin{equation}
	\Eiic[|\Bb_n|] \le \sum_{z} \P(0 \conn z \text{ with $\le n$ pivotals}).
\end{equation}
In Section \ref{sec:pivball} below we will prove
\begin{lemma}\label{prop-LinearBoundGn}\label{lem:pivball}
If $d>3(\alpha\wedge2)$ then there exists $C>0$ such that
    \begin{equation}
    \label{eqM1}
    \sum_{y \in \Zd} \Ppc ( 0 \conn y \text{ with $\le n$ pivotal edges}) \leq Cn.
    \end{equation}
\end{lemma}
Therefore, it follows that
\begin{equation}
	\Eiic[|\Bb_n|] \le Cn.
\end{equation}
Inserting this bound into \eqref{e:reffenmarkov} completes the proof. \qed

\bl \label{lem:Volcn}
If the strong triangle condition is satisfied for some sufficiently small $\beta$, then, for any sufficiently large $\lambda$ and for all $n \ge 1$,
\begin{equation}
            \label{e:Volcn0} \Piic \left(V(\zni) \le \lambda n^2\right) \ge 1-c/\sqrt{\lambda}.
\end{equation}
\el

\proof
Since $\Bb_n \subset \zni$ we have that $V(\zni) \le V(\Bb_n) + V(\zni \setminus \Bb_n)$, so that
\begin{equation}\label{e:Volcn3a}
    \Piic(V(\zni) \ge k) \le  \Piic(V(\zni \setminus \Bb_n) \ge k/2) + \Piic(V(\Bb_n) \ge k/2).
\end{equation}
It follows from Lemma \ref{lem:reffenlambda} that we can bound the second term on the right hand side by $2C/k$.

For the first term we start by applying Lemma \ref{BPRL}:
\begin{equation}\label{e:vollim}
	\Piic(V(\zni \setminus \Bb_n) \ge k/2) = \limp \frac{1}{\chi(p)} \sum_{x \in \Zd} \Pp(V(\znx \setminus \Bb^x_n) \ge k/2, 0 \conn x),
\end{equation}
where $\Bb_n^x$ is the backbone graph for the connection from $0$ to $x$ up to the $n$th pivotal. 

Observe that the set $\znx \setminus \Bb^x_n$ consists of all the vertices of $\Ccal(0)$ up to the $n$th pivotal edge that we can disconnect from $0$ by closing a single edge that touches $\Bb^x_n$, that is, let $\Bcal^x_n \equiv \{b: \,\ulb \in \Bb^x_n, b \text{ open }, \olb \notin \Bb^x_n\}$, then
\begin{equation}
    V(\znx \setminus\Bb^x_n) = \sum_{b \in E(\Zd)} V(\tCcal^b (\olb))\indi_{\{b \in \Bcal^x_n\}}.
\end{equation}
We use this identity to conclude that either of the following events must happen:
\begin{itemize}
    \item[(I)] for all $b \in \Bcal^x_n$ we have $V(\tCcal^b (\olb)) < k/4$;
    \item[(II)] there exists a $b \in \Bcal^x_n$ such that $V(\tCcal^b (\olb)) \ge k/4$.
\end{itemize}

If (I) happens, then, by the fact that  $\E[X \indi_{\{X \le R\}}] = \sum_{k=1}^{\lfloor R \rfloor} k\, \P (X =k) \le \sum_{k=1}^{\lfloor R \rfloor} \P (X \ge k)$, we get the upper bound
\begin{equation}
    \begin{split}
        \Pp \Biggl(\sum_b V(\tCcal^b (\olb)) \indi_{\{b \in \Bcal^x_n\}} &\indi_{\{V(\tCcal^b (\olb)) < k/4\}} \indi_{\{0 \conn x\}} \ge k/2\Biggr)\\
         &\le \frac{2}{k} \sum_{b} \E_p \left[V(\tCcal^b (\olb)) \indi_{\{b \in \Bcal^x_n\}} \indi_{\{V(\tCcal^b (\olb)) < k/4\}}  \indi_{\{0 \conn x\}} \right]\\
        & \le \frac{2}{k} \sum_{\ell = 1}^{\lceil k/4 \rceil} \sum_{b} \Pp \left( V(\tCcal^b (\olb)) \ge \ell,\, b \in \Bcal^x_n, 0 \conn x \right).
    \end{split}
\end{equation}
By construction the cluster $\tCcal^b(\olb)$ is ``off'' $\Bb^x$, so that we can apply the Factorization Lemma \cite[Lemma 2.2]{HofHolSla07b} again to get
\begin{equation}\label{e:Volcn4}
    \begin{split}
        \frac{2}{k} \sum_{\ell = 1}^{\lceil k/4 \rceil} \sum_{b} \Pp &\left( V(\tCcal^b (\olb)) \ge \ell,\, b \in \Bcal^x_n, 0 \conn x \right) \\
        & = \frac{2}{k}  \sum_{\ell = 1}^{\lceil k/4 \rceil} \sum_{b} p_c D(b) \E_p\left[ \indi_{\{\ulb \in \Bb_n^x, 0 \conn x \text{ on } \Bb^x\}} \P_p^{\Bb^x} \left(V(\tCcal^b (\olb)) \ge \ell \text{ off } \Bb^x \right)\right]\\
        & \le \frac{2 p}{k} \sum_{\ell=1}^{\lceil k/4 \rceil} \Ppc\left(V(\Ccal(0)) \ge \ell\right) \E_p[|\Bb_n^x|\indi_{\{0 \conn x\}}]
    \end{split}
\end{equation}
where the inequality follows from the fact that $\Pp(F $ off $A) \le \Pp(F)$ for any increasing event $F$ and any set $A$ and translation invariance of $\Pp$.
(As usual when applying the Factorization Lemma, we interpret the set $\Bb^x$ as fixed set w.r.t.\ the measure $\P_p^{\Bb^x}$, but as random set w.r.t.\ the expectation $\E_p$.)

To bound the probability on the right-hand side above we use the following lemma:
\bl[Edge volume scaling]\label{lem:Evolscaling} If the strong triangle condition is satisfied for some sufficiently small $\beta$, there exist $0 <c \le C$ such that
    \begin{equation}\label{e:Evolsc}
        \frac{c}{\sqrt{n}} \le \Ppc(V(\Ccal(0)) \ge n) \le \frac{C}{\sqrt{n}}.
    \end{equation}
\el
We will prove Lemma \ref{lem:Evolscaling} in Section \ref{sec:proofEvolscaling}. 

Inserting this bound into \eqref{e:vollim} we get the upper bound for case (I):
\begin{equation}
	\limp \frac{1}{\chi(p)}  \frac{2 p}{k} \sum_{x \in \Zd} \E_p[|\Bb_n^x|\indi_{\{0 \conn x\}}] \sum_{\ell=1}^{\lceil k/4 \rceil} \Ppc\left(V(\Ccal(0)) \ge \ell\right)   \le  \frac{C n}{k} \sum_{\ell=1}^{\lceil k/4 \rceil} \Ppc\left(V(\Ccal(0)) \ge \ell\right) \le \frac{C n}{\sqrt{k}},
\end{equation}
where the first bound can be proved similarly to Lemma \ref{lem:reffenlambda}, and the second bound follows from Lemma \ref{lem:Evolscaling}.
\medskip

If case (II) happens, then
\begin{multline}\label{e:Volcn5}
     \Pp (\exists b \in \Bcal^x_n \text{ s.t. } V(\tCcal^b (\olb)) \ge k /4, 0 \conn x) \\
     =  \sum_{(E, F) \text{adm.}} \Pp(\exists b \in F \text{ s.t. } V(\tCcal^b (\olb)) \ge k /4 \mid \Bb^x = E, \Bcal^x_n = F) \Pp(\Bb^x = E, \Bcal^x_n = F),
\end{multline}
where the sum over admissible $(E,F)$ is over all finite sets of edges $E$ and over all finite sets of directed edges $F$ such that $ \Pp(\Bb^x = E, \Bcal^x_n = F)>0$. Now, for any fixed sets $E$ and $F$, by  Boole's inequality,
\begin{equation}
    \begin{split}
        \Pp \left(\exists b \in F \text{ s.t. } V(\tCcal^b (\olb)) \ge k /4 \mid \Bb^x = E, \Bcal^x_n = F\right)& \le \sum_{b \in F} \Ppc(V(\tCcal^b(\olb)) \ge k /4 \text{ off }E \mid \Bb^x = E, \Bcal^x_n = F)\\
        & \le |F|\,\Ppc(V(\Ccal(0)) \ge k/4),
    \end{split}
\end{equation}
where the second inequality follows by the independence of $\{V(\tCcal^b(\olb)) \ge k /4 \text{ off }E\}$ and $ E = \Bb^x$, the fact that $\Ppc(A $ off $E) \le \Pp (A)$ for any increasing event $A$, and translation invariance of $\Ppc$. Applying this bound and Lemma \ref{lem:Evolscaling} to the right-hand side of \refeq{Volcn5} we get an upper bound on \eqref{e:vollim} for case (II):
\begin{equation}\label{e:Volcn6}
	\begin{split}
    		\limp \frac{1}{\chi(p)} \sum_x \frac{C}{\sqrt{k}} \sum_{(E, F) \text{adm.}} |F|\, \Ppc(\Bb^x = E, \Bcal^x_n = F)& = \limp \frac{1}{\chi(p)} \sum_{x} \frac{C}{\sqrt{k}} \E_p [|\Bcal^x_n| \indi_{\{0 \conn x\}} ] \\
		& \le \limp \frac{1}{\chi(p)}  \frac{C}{\sqrt{k}} \sum_x \E_p[V(\Bb_n^x)\indi_{\{0 \conn x\}}] \\
		& \le \frac{C n}{\sqrt{k}},
	\end{split}
\end{equation}
where the final inequality follows from a proof similar to that of Lemma \ref{lem:reffenlambda}.

Combining the bounds for cases (I) and (II) with $k = \lambda n^2$ we thus get
\begin{equation}
	 \Piic(V(\zni) \ge \lambda n^2) \le \frac{C n}{\sqrt{\lambda n^2}} + \frac{C' n}{\sqrt{\lambda n^2}} \le \frac{C+C'}{\sqrt{\lambda}}.
\end{equation}
as desired. \qed

\subsection{An upper bound on the volume of the outer-pivotal ball: proof of Lemma \ref{lem:pivball}}\label{sec:pivball}
\label{sec-num-pivs}
In this section, we prove Lemma \ref{lem:pivball}. 
Write 
\begin{equation}
	G(n) = \sum_{y \in \Zd} \Ppc(0 \conn y \text{ with $\le n$ pivotal edges}).
\end{equation}
We will prove Lemma \ref{lem:pivball} by induction using recursive upper and lower bounds on $G(n)$ (stated in Lemmas \ref{lemma-M1} and \ref{lemma-UpperBoundG2n} below). 
\begin{lemma}\label{lemma-M1}
    If $d>3(\alpha\wedge2)$, then exists $c_1>0$ such that for $n\in\N$,
    \begin{equation}\label{eqM3}
        G(2n)\ge \frac{c_1}{n} G(n)^2.
    \end{equation}
\end{lemma}
\proof
The proof is inspired by the proof of \cite[Lemma 3.1]{KozNac09}.

Write $N_n$ for the number of pairs $(b,y)$, such that the following events all occur:
\begin{itemize}
	\item $0 \conn y$, 
	\item the edge $b=(\bb,\tb)$ is pivotal for $0 \conn y$,
	\item $0\conn\bb$ with $\le n$ pivotal edges,
	\item $\tb\conn y$ with $\le n$ pivotal edges.
\end{itemize}
 If a pair $(b,y)$ contributes towards $N_n$,
then we must have that $\{0\conn y$ with $\le 2n$ pivotal edges$\}$, and there are at most $2n$ choices for $b$,
so
    \begin{equation}\label{eqM4}
    N_n\le2n\,\sum_{y\in \Z^d} \indic{0\conn y \text{ with $\leq 2n$ pivotal edges}}.
    \end{equation}
It follows that $\expec_{p_c}[N_n]\le2n\,G(2n)$.

We now bound $\expec_{p_c}[N_n]$ from below.
By the Factorization Lemma \cite[Lemma 2.2]{HofHolSla07b} we get
    \begin{align}
    \expec_{p_c}[N_n]=
    \sum_{\text{edge } b}\;\sum_{y\in\Zd}\expec_{p_c}\bigg(&
        \indi_{\{0\conn\bb \text{ with $\le n$ piv's on }\Ct\}}\;
        \indi_{\{b\text{ occupied}\}}\;\nnb\label{eqM5}
        &{}\times\P_{p_c}^{\Ct}\left(\tb\conn y\text{ with $\le n$ piv's off }\Ct\right)\bigg),
    \end{align}
where again the set $\Ct$ is to be considered as a \emph{fixed} set w.r.t.\ the probability measure $\Ppc^\Ct$, but it is a random set w.r.t.\ the overall expectation $\expec_{p_c}$.

Next we replace the indicator function in the first line of \eqref{eqM5} by $\indi_{\{0\conn\bb \text{ with $\le n$ piv's}\}}$. The reason for this is as follows: Suppose, $\{0\conn\bb \text{ with $\le n$ piv's}\}$ occurs, but not $\{0\conn\bb \text{ with $\le n$ piv's on }\Ct\}$. 
Then the (restricted) probability on the second line of \eqref{eqM5} equals 0, so these configurations do not contribute to the expectation. 

We expand the probability in the last line of \eqref{eqM5} as
\begin{multline}\label{eqM6}
    \P_{p_c}^{\Ct} \left(\tb\conn y\text{ with $\le n$ piv's off }\Ct\right)\\
    = \Ppc \left(\tb\conn y\text{ with $\le n$ piv's}\right) -\P_{p_c}^{\Ct} \left(\tb\conn y\text{ with $\le n$ piv's through }\Ct\right),
\end{multline}
where we used the definition of ``through'' (see Definition \ref{def:def}(iv)).

The first term in \eqref{eqM6} is equal to
\begin{equation}\label{eqM10}
  \sum_{\text{edge } b}\;\sum_{y\in\Zd}\expec_{p_c}\bigg[
        \indi_{\{0\conn\bb \text{ with $\le n$ piv's}\}}\;
        \indi_{\{b\text{ occupied}\}}\;\Ppc(\tb\conn y\text{ with $\le n$ piv's})\bigg]  = p_c G(n)^2.
\end{equation}
We bound the last term in \eqref{eqM6} from above by
\begin{multline}\label{eqM8}
    \sum_{\text{edge } b}\;\sum_{y\in\Zd}\expec_{p_c}\bigg[ \indi_{\{0\conn\bb \text{ with $\le n$ piv's}\}}\indi_{\{b\text{ occupied}\}} \Ppc^{\Ct}\big(\tb\conn y\text{ with $\le n$ piv's through }\Ct\big)\bigg]\\
        \le  \sum_{\text{edge } b}\;\sum_{y,v,w\in\Zd}\Ppc\Big( \{0\conn v \text{ with $\le n$ piv's}\} \circ\{b\text{ occupied}\} \circ\{v\conn \bb\}\circ\{v\conn w\}\Big)\\
            \times\Ppc\Big(\{\tb\conn w\}\circ\{w\conn y \text{ with $\le n$ piv's}\}\Big).
\end{multline}

We apply the BK-inequality to both probabilities in \eqref{eqM8} and take sums in the right order to obtain that \eqref{eqM8} is bounded above by
    \begin{equation}\label{eqM9}
      p_c\,G(n)^2\,\tilde T_{p_c}.
    \end{equation}
 where 
 \begin{equation}
 	\tilde T_{p_c} = \sum_{a,b,c} \Ppc(0 \conn a) p_c D(b-a) \Ppc(b \conn c) \Ppc(c \conn 0).
\end{equation}
It is proved in \cite{HarSla90a,HeyHofSak08} that $\tilde T_{p_c} \le C \beta$ for both high-dimensional finite-range and long-range percolation, respectively.
We insert \eqref{eqM6} into \eqref{eqM5}, and bound the two summands with \eqref{eqM10} and \eqref{eqM8}--\eqref{eqM9}, respectively,
to arrive at
    \begin{equation}\label{eqM11}
    \expec_{p_c}[N_n]\ge p_c\,(1-C \beta)\,G(n)^2.
    \end{equation}
Since $\expec_{p_c}[N_n]\le2n\,G(2n)$ from \eqref{eqM4},
this proves the claim for $c_1=p_c\,(1- C \beta)/2$, which is positive when $\beta$ is sufficiently small.
\qed

A complementary bound to \eqref{eqM3} is given by the next lemma.
\begin{lemma}\label{lemma-UpperBoundG2n}
There is a constant $C$ such that
\begin{equation}\label{eqM11}
    G(2n)\le C\;\frac{G(n)}{n^{\twa}}.
\end{equation}
\end{lemma}
\proof[Proof of Lemma \ref{lemma-UpperBoundG2n}]
We fix $\eps>0$, and let
\begin{equation}\label{eqM23}
    N=\frac\eps n\,G(2n)^{1/\twa}.
\end{equation}
Write $e_1$, \dots,$e_{\ell}$ for the ordered sequence of pivotal edges for the event $\{0 \conn y\}$ on the path from $0$ to $y$, write $\bbe_0=0$, and write $y = \bbe_{\ell+1}$. If the event $\{0\conn y \text{ with $\le 2n$ pivotals}\}$ occurs then there are two possibilities:
\begin{enumerate}
  \item[(i)]
    $|\bbe_i - \bbe_{i-1}|\le N$ for all $i=1,\dots,\ell+1$, or
  \item[(ii)]
    there exists an $i\in\{1,\dots,\ell+1 \}$ such that $|\bbe_i - \bbe_{i-1}|>N$.
\end{enumerate}
In case (i), we must have that $|y|\le 2nN$, so the contribution towards $G(2n)$ from this case is bounded above by
\begin{equation}\label{eqM12}
    \sum_{|y|\le 2nN}\prob_{p_c}(0\conn y)
    \:\le\: C (2nN)^{\twa}
    \:\le\: C_1\,\eps^{\twa}\,G(2n),
\end{equation}
where we used Theorem \ref{ExpectationBounds}.

The contribution towards $G(2n)$ from case (ii) is bounded above by
\begin{equation}\label{eqM13}
    \sum_{u\colon|u-v|>N}
    \Big(\sum_{v\in\Zd}p_cD(v)\,\prob_{p_c}(v\Conn  u)\Big)\,G(n)\;G(2n),
\end{equation}
and the sum of \eqref{eqM12} and \eqref{eqM13} forms an upper bound on $G(2n)$.
Here we denote by $\{v\Conn u\}\equiv\{v\conn u\}\circ\{v\conn u\}$ a \emph{double connection} between $u$ and $v$.

It follows from \cite[Proposition 2.5 and Remark 2.6]{HeyHofHul12a} that there exists a constant $C_2$ such that
\begin{equation}
	\sum_{|u-v| > N} \sum_{v \in \Zd}  p_c D(v) \P(v \Conn u) \le \sum_{u,v \in \Zd} \left| \frac{u-v}{N}\right|^\twa  p_c D(v) \P(v \Conn u) \le C_2 N^{-\twa}.
\end{equation}

Together with \eqref{eqM23} and \eqref{eqM12}, this implies
\begin{equation}\label{eqM15}
    G(2n)\le C_1\,\eps^{\twa}\,G(2n)+C_2\,\frac{n^{\twa}}{\eps^{\twa}}\,G(n),
\end{equation}
and choosing $\eps$ such that $C_1\eps^{\twa}=1/2$ implies the claim. \qed

Indeed, Lemmas \ref{lemma-M1} and \ref{lemma-UpperBoundG2n} imply Lemma \ref{prop-LinearBoundGn}:
\proof[Proof of Lemma \ref{prop-LinearBoundGn}.]
The proof is similar to the proof of \cite[Theorem 1.2(i)]{KozNac09}, though our infinite-range setting requires Lemma \ref{lemma-UpperBoundG2n} as an additional ingredient.

We give a proof by contradiction:

Assume that
$G(n_0)\ge \bar C\,n_0$
for a constant $\bar C$ that satisfies $\bar C>\max\{2,\,2/c_1,\,2^{\twa+1}\}$ and $n_0\in\N$.
We claim that this implies
\begin{equation}\label{eqM17}
    G(2^kn_0)\ge \bar C^{k+1}\,n_0.
\end{equation}
The claim is proved by induction.
The case $k=0$ is our assumption, and we advance the induction by using  \eqref{eqM3}, the induction hypothesis \eqref{eqM17}, and finally $\bar C>\max\{2,\,2/c_1\}$, to obtain
\begin{equation}\label{eqM18}
    G(2^{k+1}n_0)
    \ge \frac{c_1\,G(2^kn_0)^2}{2^kn_0}
    \ge \frac {c_1}{2^k}\,\bar C^{2k+2}\,n_0
    \ge \bar C^{k+2}\,n_0.
\end{equation}
The combination of \eqref{eqM3} and \eqref{eqM11} implies
\begin{equation}
	 \frac {c_1} n\,G(n)^2\le G(2n) \le C\,n^{\twa}G(n),
\end{equation}
so $G(n)\le \left(C/c_1\right)\,n^{\twa+1}$.
Together with \eqref{eqM18},
\begin{equation}\label{eqM19}
    \bar C^{k+1}n_0\le G(2^kn_0)\le \frac C c_1 2^{(\twa+1)k}\,n_0^{\twa+1}.
\end{equation}
For large $k$ this causes a contradiction because we chose $\bar C$ such that $\bar C>2^{\twa+1}$.
\qed

\subsection{Tail estimates for the cluster size of long-range models: proof of Lemma \ref{lem:Evolscaling}}\label{sec:proofEvolscaling}
In this section we prove Lemma \ref{lem:Evolscaling}. 
\proof[Proof of Lemma \ref{lem:Evolscaling}.]
Combining results from \cite{BarAiz91} and \cite{HeyHofSak08} we know that there exist $0<c'\le C'$ such that
\begin{equation}\label{e:Vvolsc}
    \frac{c'}{\sqrt{n}} \le \Ppc(\lvert \Ccal(0) \rvert \ge n) \le \frac{C'}{\sqrt{n}}.
\end{equation}
Recall that $V(\Ccal(0))$ is the edge volume of $\Ccal(0)$.
Since $V(\Ccal(0)) \ge \lvert \Ccal(0) \rvert-1$, the lower bound in \refeq{Evolsc} immediately follows. For models where each vertex has bounded degree, i.e., for any finite range model, the upper bound also follows immediately , since then $V(\Ccal(0)) \le \Delta_{\sss \text{max}} \lvert \Ccal(0) \rvert$, where $\Delta_{\sss \text{max}}$ is the maximal degree a vertex can have. But establishing this bound for models with unbounded degree requires a bit more work.

Fixing a constant $\gamma >0$ (to be determined later) we can bound
\begin{equation}\label{e:VCsplit}
    \begin{split}
        \Ppc(V(\Ccal(0)) \ge n)  &\le \Ppc(\lvert \Ccal(0) \rvert \ge \gamma n) + \Ppc(V(\Ccal(0)) \ge n \mid \lvert \Ccal(0) \rvert < \gamma n)\\
                                &\le \frac{C'}{\sqrt{\gamma n}} + \Ppc(V(\Ccal(0)) \ge n \mid \lvert \Ccal(0) \rvert < \gamma n).
    \end{split}
\end{equation}

Now, consider a fixed $\Ccal(0)$ such that $\lvert \Ccal(0) \rvert = m$. We do a standard \emph{breadth-first exploration} of $\Ccal(0)$, starting at $0$ and labeling the vertices as we go along. Let $v_i$ be the $i$th explored vertex of $\Ccal(0)$. After the exploration is completed, we know that there are $m-1$ explored open edges.

Given that we know the exploration tree, the only way that the unexplored edge $\{x, y\}$, with $x,y \in \Ccal(0)$, could be open is if both the topology and the labeling of the exploration tree are not affected by information on the status of $\{x, y\}$.
Clearly, any open edge in $\Ccal(0)$ that has not been explored has to be an edge between two explored vertices, so suppose that $x$ and $y$ have been explored. If $x$ and $y$ are more than two generations apart in the exploration tree and the edge $\{x, y\}$ has not been explored, then it is closed, because the exploration would have explored the edge if it was open. The edge $\{x, y\}$ is also closed when having the edge be open would be inconsistent with the order of the exploration. When the status of the edge affects neither the topology nor the labeling, then the status of the edge under $\Ppc$ is \emph{independent} of the status of other edges, since $\Ppc$ is a product measure. So the number of unexplored edges at vertex $v$, say, is stochastically bounded from above by $\Delta_v$, a random variable whose law $\P_{\Delta}$ is the degree of vertex $v$ under $\Ppc$ (in the setting without conditioning).

These observations imply that we have the following bound: for all $x$,
\begin{equation}
    \Ppc(V(\Ccal(0)) \ge x \mid \lvert \Ccal(0) \rvert < \gamma n) \le \P_{\Delta} \left( \gamma n + \sum_{i=1}^{\lceil \gamma n \rceil} X_i \ge x\right) \qquad \text{ with } X_i \sim \Delta_0 \text{ for all $i$}.
\end{equation}
Furthermore, by our choice of $D$ as described in the introduction we have that $\Epc[\Delta_0]  = p_c$ and it is easy to show that $\text{Var}_{p_c} (\Delta_0) = \sigma^2 < \infty$. Thus, choosing $\gamma < (2 p_c)^{-1}$ we can use Cantelli's inequality to bound
\begin{equation}\label{e:Cheby}
    \begin{split}
        \Ppc(V(\Ccal(0)) \ge n \mid \lvert \Ccal(0) \rvert < \gamma n) \le & \P_{\Delta}\left(\sum_{i=1}^{\lceil \gamma n \rceil} X_i -p_c \gamma n\ge (1- 2 p_c \gamma)n \right)\\
         \le & \left(1+\frac{(1-2 p_c \gamma)^2 n^2}{n \sigma^2}\right)^{-1} \le \frac{C}{n}.
    \end{split}
\end{equation}
Using \refeq{Cheby} with \refeq{VCsplit} completes the proof. \qed

\appendix
\section{Proof of Theorem \ref{th:Eucexittimes}(a), (b) and (c)}\label{app:Pfth14}
This appendix contains a proof of Theorem \ref{th:Eucexittimes}(a), (b) and (c) that closely follows the proof that Kumagai and Misumi \cite{KumMis08} give for a similar, more general results. We have made a few small modifications to their proof to make it work for the extrinsic distance metric.
Before we prove Theorem \ref{th:Eucexittimes} (a), (b) and (c) we state a lemma that gives some bounds on exit times in terms of bounds on volume and effective resistance:
\bl[Parts of Proposition 3.3 and 3.5 from \cite{KumMis08}]\label{prop:33}
Let $\lambda > 0$.
\begin{itemize}
    \item[(1)] Suppose $r \in \jel$. Then, for $z \in \bqr$,
    \begin{equation}
        \label{e:RVI}   E_{\omega}^z \tqr \le 2 \lambda^2 r^{6}.
    \end{equation}
    \item[(2)] Let $\vep^{2} =\vep(\lambda)^{2} =1/(8 \lambda^2)$. If $r, \vep r \in \jel$, then
    \begin{equation}
        \label{e:A2}    E_{\omega}^x \tqr \ge  \frac{r^{6}}{2^{11} \lambda^6}  \text{ for } x \in \bqhr.
    \end{equation}
\end{itemize}
\el
\proof We start by noting that by Definition \ref{def:properties}, for any $z \in \bqr$,
\begin{equation}
    E_{\omega}^z \tqr = \sum_{y \in \bqr} G_{Q_r}(z,y) \mu_y
\end{equation}

\emph{Proof of (1).} $\Reff$ is a metric, so it satisfies the triangle condition, i.e.,\ for any $z \in \bqr$,
\begin{equation}
    \Reff(z,Q_r^c) \le \Reff(0,z) + \Reff(0,Q_r^c) \le 2 \lambda r^{2}
\end{equation}
where the second inequality follows from the fact that $r \in \jel$ and Definition \ref{def:radiussets}(i). Since for any $y,z$, we have $G_{Q_r}(z,y) \le G_{Q_r}(z,z)$ and $G_{Q_r}(z,z) = \Reff(z, Q_r^c)$, it follows that
\begin{equation}\label{e:A1}
    E_{\omega}^z \tqr = \sum_{y \in \bqr} G_{Q_r} (z,y) \mu_y \le \sum_{y \in Q_r} G_{Q_r}(z,z) \mu_y = \Reff(z, Q_r^c)\Vr \le 2 \lambda^2 r^{6},
\end{equation}
where we have used Definition \ref{def:radiussets}(i) again for the final inequality.\\
\\
\emph{Proof of (2).} Since $r \in \jel$, we have the following bound for any $x \in \bqi{\vep r}$:
\begin{equation}
    \frac{r^{2}}{\lambda} \le \ReffB \le \Reff(0,x) + \Reff(x, Q_r^c) \le \lambda (\vep r)^{2} + \Reff(x, Q_r^c).
\end{equation}
So if we take $\vep >0$ sufficiently small, we get
\begin{equation}\label{e:Ass1}
     \Reff(x, Q_r^c) \ge \frac{r^{2}}{2 \lambda} \qquad \text{ for all }x \in \bqi{\vep r}.
\end{equation}

Let $\pqr (y) = \gqr (x,y)/\gqr(x,x)$. Since
\begin{equation}
    \lvert f(x) -f(y) \rvert^2 \le \Reff(x,y) \Ecal (f,f) \qquad \text{ for all } f \in L^2 (\Gamma, \mu)
\end{equation}
and furthermore $\Ecal(\pqr, \pqr) = \Reff(x, Q_r^c)^{-1} = \gqr (x,x)^{-1}$, we have, for $x,y \in \bqhr$,
\begin{equation}
    \lvert 1 - \pqr (y) \rvert^2  \le \frac{\Reff(x,y)}{\Reff(x, Q_r^c)} \le \frac{4 \lambda^2 (\vep r)^{2}}{r^{2}} = \frac12.
\end{equation}
Hence, $\pqr (y) \ge 1 - 1/\sqrt{2} \ge 1/4$, so that
\begin{equation}\label{e:UB1}
    E_{\omega}^x \tqr \ge \sum_{y \in \bqi{\vep r}} \gqr (x,x) \pqr(y) \mu_y \ge \frac14 \mu(\bqhr) \Reff (x, Q_r^c) \ge \frac{r^{2} V(\bqhr)}{8 \lambda} \ge \frac{r^{6}}{2^{11} \lambda^6}.
\end{equation}\qed

\proof[Proof of Theorem \ref{th:Eucexittimes} (a), (b) and (c)] We start with \refeq{RXI}. Choose $\lambda \ge 1$ such that $2 c_E \lambda^{-q_E} < \delta$. Let $r \ge r^\star$ and set $F_1 = \{r, \vep r \in \jel\}$. Suppose first that $\vep r \ge 1$. Then, by Theorem \ref{th:Goodballs}(a), $\Piic(F_1) \ge 1 - 2 c_E \lambda^{-q_E}$. For $\omega \in F_1$, by Lemma \ref{prop:33}, there exist $c_1 < \infty$ and $q_1 \ge 0$ such that
\begin{equation}\label{e:A6}
    (c_1 \lambda^{q_1})^{-1} \le \frac{E_{\omega}^x \tqr}{r^{6}} \le c_1 \lambda^{q_1} \qquad \text{ for } x \in \bqi{\vep r}.
\end{equation}
Thus, if $\theta_0 = c_1 \lambda^{q_1}$, then for $\theta = \theta_0$ it follows that
\begin{equation}
    \Piic \left(\theta^{-1} r^{6} \le \Eotr \le \theta r^{6}\right) \ge \Piic(F_1) \ge 1 -\delta.
\end{equation}

Now consider the case where $r \le 1/\vep$. For each graph $\Gamma_{{\sss \mathsf{IIC}}}(\omega)$, let
\begin{equation}
    Y(\omega) = \sup_{1 \le s \le 1/\vep}\frac{E_{\omega}^{0} \tau_{Q_s}}{s^{6}}.
\end{equation}
Then $Y(\omega) < \infty$ for each $\omega$, so there exists $\theta_1$ such that
\begin{equation}
    \Piic\left(\Eotr > \theta_1 r^{6}\right) \le \Piic(Y > \theta_1) \le \delta.
\end{equation}
If we take $\theta_1 > \vep(\lambda)^{-6}$, then $\Eotr \ge \theta_1^{-1} r^{6}$, since $\Eotr \ge 1$. Thus, for $\theta \ge \theta_1$, we also have
\begin{equation}
    \Piic \left(\theta^{-1} r^{6} \le \Eotr \le \theta r^{6}\right) \ge 1-\delta,
\end{equation}
which completes the proof of \refeq{RXI}.\\
\\
Now we prove \refeq{RXVI}. We begin with the upper bound. By \refeq{A1} and Theorem \ref{th:Goodballs}(a),
\begin{equation}
    \Eiic[\Eotr] \le \Eiic[\ReffB \Vr] \le c r^{6}.
\end{equation}

For the lower bounds, it suffices to find a set $F \subset \Omega$ of `nice' graphs with $\Piic(F) \ge c > 0$ such that, for all $\omega \in F$, we have a suitable lower bound on $\Eotr$. Assume that $r \ge r^\star$ is large enough so that $\vep (\lambda_0) r \ge 1$, where $\lambda_0$ is chosen large enough so that $c_E \lambda_0^{-q_E} < 1/8$. We can then get results for all $n$ (chosen below to depend on $r$) and $r$ by adjusting the constant $c_1$ in \refeq{RXVI}.

Let $F = \{r, \vep (\lambda_0) r \in \jel\}$. Then $\Piic(F) \ge 3/4$, and for $\omega \in F$, by \refeq{A2}, $\Eotr \ge \cso{1} r^{6}$, so that
\begin{equation}
    \Eiic [\Eotr] \ge \Eiic[\Eotr \indi_{\{F\}}] \ge \cso{1} r^{6} \Piic(F) \ge \cso{2} r^{6}.
\end{equation}
\\
Finally we prove \refeq{RXX}. Let $r_n = e^n$ and $\lambda_n = n^{2/q_E}$. Let $F_n = \{r_n, \vep(\lambda_n) r_n \in J_{{\rm Euc}}(\lambda_n)\}$. Then $\Piic(F_n^c) \le 2 n^{-2}$ (provided $\vep(\lambda_n) r_n \ge 1$). Therefore, by Borel-Cantelli, if $\Omega_b = \liminf F_n$, then $\Piic(\Omega_b)=1$. Hence there exists $M_0$ with $M_0(\omega) <\infty$ on $\Omega_b$ and such that $\omega \in F_n$ for all $n \ge M_0 (\omega)$.

Choose a fixed $\omega \in \Omega_b$ and let $x \in \iic(\omega)$. By \refeq{A6} there exist constants $c_2, q_2$ such that
\begin{equation}\label{e:A7}
    (c_2 \lambda_n^{q_2})^{-1} \le \frac{E_{\omega}^{x} \tau_{Q_{r_n}}}{r_n^{6}} \le c_2 \lambda_n^{q_2},
\end{equation}
provided that $n \ge M_0 (\omega)$ and $n$ is also large enough so that $x \in \bqi{\vep(\lambda_n) r_n}$. Writing $M_x (\omega)$ for the smallest such $n$, we have
\begin{equation}
    c_2^{-1} (\log r_n)^{-2 q_2 / q_E} r_n^{6} \le E_{\omega}^{x} \tau_{Q_{r_n}} \le c_2 (\log r_n)^{2 q_2 / q_E} r_n^{6} \qquad \text{ for all }n \ge M_x (\omega).
\end{equation}
Note that by definition $E_{\omega}^{x} \tau_{Q_{r_n}}$ is increasing in $n$. If $r \ge R_x = 1+e^{M_x}$, then let $r$ be such that $r_{n-1} \le r \le r_n$. Then,
\begin{equation}
    E_{\omega}^{x} \tau_{Q_{r}} \le E_{\omega}^{x} \tau_{Q_{r_n}} \le c_2 (\log r_n)^{2 q_2 / q_E} r_n^{6} \le c'_2 (\log r)^{2 q_2/q_E} r^{6}.
\end{equation}
Similarly,
\begin{equation}
    E_{\omega}^{x} \tau_{Q_{r}} \ge E_{\omega}^{x} \tau_{Q_{r_{n-1}}} \ge c_3(\log r_{n-1})^{-2 q_2 /q_E} r_{n-1}^{6} \ge c'_3 (\log r)^{-2 q_2/q_E} r^{6}.
\end{equation}
Taking $\xi > 2 q_2 /q_E$ large enough to absorb the constants $c'_2$ and $c'_3$ in the $\log r$ term, we get \refeq{RXX}. \qed

\vspace{0.5cm}
{\bf Acknowledgements.}
The work of MH, RvdH and TH was supported in part by the Netherlands Organization for Scientific Research (NWO). MH and RvdH thank the Institut Henri Poincar\'e Paris for kind hospitality during their visit in October 2009. TH thanks Robert Fitzner for stimulating discussions.

\bibliographystyle{abbrv}

\begin{thebibliography}{10}

\bibitem{Aize97}
M.~Aizenman.
\newblock On the number of incipient spanning clusters.
\newblock {\em Nucl. Phys. B}, {\bf 485}:551--582, (1997).

\bibitem{AleOrb82}
S.~Alexander and R.~Orbach.
\newblock Density of states on fractals: ``fractons''.
\newblock {\em J. Physique (Paris) Lett.}, {\bf 43}:L625--L631, (1982).

\bibitem{Barl04}
M.~Barlow.
\newblock Random walks on supercritical percolation clusters.
\newblock {\em Ann. Probab.}, {\bf 32}(4):3024--3084, (2004).

\bibitem{BarJarKumSla08}
M.~Barlow, A.~J{\'a}rai, T.~Kumagai, and G.~Slade.
\newblock Random walk on the incipient infinite cluster for oriented
  percolation in high dimensions.
\newblock {\em Comm. Math. Phys.}, {\bf 278}(2):385--431, (2008).

\bibitem{BarKum06}
M.~T. Barlow and T.~Kumagai.
\newblock Random walk on the incipient infinite cluster on trees.
\newblock {\em Illinois J. Math.}, 50(1-4):33--65 (electronic), (2006).

\bibitem{BarAiz91}
D.~Barsky and M.~Aizenman.
\newblock Percolation critical exponents under the triangle condition.
\newblock {\em Ann. Probab.}, {\bf 19}:1520--1536, (1991).

\bibitem{BerKes85}
J.~v.~d. Berg and H.~Kesten.
\newblock Inequalities with applications to percolation and reliability.
\newblock {\em J. Appl. Prob.}, {\bf 22}:556--569, (1985).

\bibitem{BerGanPer03}
N.~Berger, N.~Gantert, and Y.~Peres.
\newblock The speed of biased random walk on percolation clusters.
\newblock {\em Probab. Th. Rel. Fields}, {\bf 126}(2):221--242, (2003).

\bibitem{Bill99}
P.~Billingsley.
\newblock {\em Convergence of Probability Measures}.
\newblock John Wiley and Sons, New York, 2nd edition, (1999).

\bibitem{BorChaHofSlaSpe05b}
C.~Borgs, J.~Chayes, R.~v.~d. Hofstad, G.~Slade, and J.~Spencer.
\newblock Random subgraphs of finite graphs. {II}. {T}he lace expansion and the
  triangle condition.
\newblock {\em Ann. Probab.}, {\bf 33}(5):1886--1944, (2005).

\bibitem{CheSak12}
L.-C. Chen and A.~Sakai.
\newblock Critical two-point functions for long-range statistical-mechanical
  models in high dimensions.
\newblock {\em Ann. Probab.}, to appear.

\bibitem{Croy08}
D.~A. Croydon.
\newblock Random walk on the range of random walk.
\newblock {\em J. Stat. Phys.}, 136(2):349--372, (2009).

\bibitem{DoySne84}
P.~Doyle and J.~Snell.
\newblock {\em Random walks and electric networks}, volume~{\bf 22} of {\em
  Carus Mathematical Monographs}.
\newblock Mathematical Association of America, Washington, DC, (1984).

\bibitem{Fitz13}
R.~Fitzner.
\newblock {\em Non-backtracking lace expansion}.
\newblock PhD thesis, Eindhoven University of Technology, (2013).

\bibitem{Grim99}
G.~Grimmett.
\newblock {\em Percolation}.
\newblock Springer, Berlin, 2nd edition, (1999).

\bibitem{Hara08}
T.~Hara.
\newblock Decay of correlations in nearest-neighbor self-avoiding walk,
  percolation, lattice trees and animals.
\newblock {\em Ann. Probab.}, {\bf 36}(2):530--593, (2008).

\bibitem{HarHofSla03}
T.~Hara, R.~v.~d. Hofstad, and G.~Slade.
\newblock Critical two-point functions and the lace expansion for spread-out
  high-di\-men\-sio\-nal percolation and related models.
\newblock {\em Ann. Probab.}, {\bf 31}(1):349--408, (2003).

\bibitem{HarSla90a}
T.~Hara and G.~Slade.
\newblock Mean-field critical behaviour for percolation in high dimensions.
\newblock {\em Commun. Math. Phys.}, {\bf 128}:333--391, (1990).

\bibitem{HarSla94}
T.~Hara and G.~Slade.
\newblock Mean-field behaviour and the lace expansion.
\newblock In G.~Grimmett, editor, {\em Probability and Phase Transition},
  Dordrecht, (1994). Kluwer.

\bibitem{HeyHofHul12a}
M.~Heydenreich, R.~v.~d. Hofstad, and T.~Hulshof.
\newblock High-di\-men\-sio\-nal incipient infinite clusters revisited.
\newblock arXiv:1108.4325v2, (2012).

\bibitem{HeyHofHulMie12a}
M.~Heydenreich, R.~v.~d. Hofstad, T.~Hulshof, and G.~Miermont.
\newblock Backbone scaling limit of the high-di\-men\-sio\-nal {IIC}.
\newblock In preparation.

\bibitem{HeyHofSak08}
M.~Heydenreich, R.~v.~d. Hofstad, and A.~Sakai.
\newblock Mean-field behavior for long- and finite range {I}sing model,
  percolation and self-avoiding walk.
\newblock {\em J. Statist. Phys.}, {\bf 132}(5):1001--1049, (2008).

\bibitem{HeyHof11}
M.~Heydenreich and R.~van~der Hofstad.
\newblock Random graph asymptotics on high-dimensional tori {II}: volume,
  diameter and mixing time.
\newblock {\em Probab. Theory Related Fields}, 149(3-4):397--415, (2011).

\bibitem{HofHolSla02}
R.~v.~d. Hofstad, F.~d. Hollander, and G.~Slade.
\newblock Construction of the incipient infinite cluster for spread-out
  oriented percolation above {$4+1$} dimensions.
\newblock {\em Comm. Math. Phys.}, {\bf 231}(3):435--461, (2002).

\bibitem{HofHolSla07b}
R.~v.~d. Hofstad, F.~d. Hollander, and G.~Slade.
\newblock The survival probability for critical spread-out oriented percolation
  above 4+1 dimensions. {II.} {E}xpansion.
\newblock {\em Ann. Inst. H. Poincar\'e Probab. Statist.}, {\bf 5}(5):509--570,
  (2007).

\bibitem{HofJar04}
R.~v.~d. Hofstad and A.~J{\'a}rai.
\newblock The incipient infinite cluster for high-di\-men\-sio\-nal unoriented
  percolation.
\newblock {\em J. Statist. Phys.}, {\bf 114}(3-4):625--663, (2004).

\bibitem{HofSap11}
R.~v.~d. Hofstad and A.~Sapozhnikov.
\newblock Cycle structure of percolation on high-dimensional tori.
\newblock {\it arXiv:1109.1233}, (2011).

\bibitem{Jara03b}
A.~J{\'a}rai.
\newblock Incipient infinite percolation clusters in 2{D}.
\newblock {\em Ann. Probab.}, {\bf 31}(1):444--485, (2003).

\bibitem{JarNac13}
A.~J{\'a}rai and A.~Nachmias.
\newblock Electrical resistance of the low dimensional critical branching
  random walk.
\newblock arXiv:1305.1092, (2013).

\bibitem{Kest80}
H.~Kesten.
\newblock The critical probability of bond percolation on the square lattice
  equals {$\frac12$}.
\newblock {\em Comm. Math. Phys.}, {\bf 74}(1):41--59, (1980).

\bibitem{Kest86a}
H.~Kesten.
\newblock The incipient infinite cluster in two-di\-men\-sio\-nal percolation.
\newblock {\em Probab. Th. Rel. Fields}, {\bf 73}(3):369--394, (1986).

\bibitem{Kest86b}
H.~Kesten.
\newblock Subdiffusive behavior of random walk on a random cluster.
\newblock {\em Ann. Inst. H. Poincar\'e Probab. Statist.}, 22(4):425--487,
  (1986).

\bibitem{KozNac09}
G.~Kozma and A.~Nachmias.
\newblock The {A}lexander-{O}rbach conjecture holds in high dimensions.
\newblock {\em Invent. Math.}, {\bf 178}:635--654, (2009).

\bibitem{KozNac11}
G.~Kozma and A.~Nachmias.
\newblock Arm exponents in high-di\-men\-sio\-nal percolation.
\newblock {\em J. Amer. Math. Soc.}, {\bf 24}:375--409, (2011).

\bibitem{KumMis08}
T.~Kumagai and J.~Misumi.
\newblock Heat kernel estimates for strongly recurrent random walk on random
  media.
\newblock {\em J. Theoret. Probab.}, {\bf 21}:910--935, (2008).

\bibitem{LeySta83}
F.~Leyvraz and H.~Stanley.
\newblock To what class of fractals does the {A}lexander-{O}rbach conjecture
  apply?
\newblock {\em Phys. Rev. Lett.}, {\bf 51}:2048--2051, (1983).

\bibitem{NacPer08}
A.~Nachmias and Y.~Peres.
\newblock Critical random graphs: diameter and mixing time.
\newblock {\em Ann. Probab.}, 36(4):1267--1286, (2008).

\bibitem{Sapo10}
A.~Sapozhnikov.
\newblock Upper bound on the expected size of the intrinsic ball.
\newblock {\em Elect. Comm. in Probab.}, {\bf 15}:297--298, (2010).

\end{thebibliography}
\begin{small}

\end{small}
\end{document}